\documentclass[12pt]{amsart}
\usepackage{amssymb}
\usepackage{amsmath} 
\usepackage{pictex}
\newcommand\AND{\quad\text{and}\quad}

\newcommand\bd{\partial}
\newcommand\BS{\operatorname{\sf BS}}

\newcommand\cf{\curlywedge}

\newcommand\ep{\varepsilon}

\newcommand\geo[1]{\overline{#1}}

\newcommand\GM{\mathsf{\Gamma}\mathsf{M}}
\newcommand\GMa{\mathsf{\Gamma_1}\mathsf{M_1}}
\newcommand\GMb{\mathsf{\Gamma_2}\mathsf{M_2}}
\newcommand\GMQ{\mathsf{\Gamma_0}\mathsf{M_0}}
\newcommand\Hb{{\mathbb H}}
\newcommand\hor{\mathfrak{h}}
\newcommand\HT{\operatorname{\sf HT}}
\newcommand\im{\mathfrak{i}\,}
\newcommand\Id{\text{\sl Id}}
\newcommand\IM{\operatorname{\text{\sl Im}}}
\newcommand\jm{\mathfrak{j}}
\newcommand\Lf{\mathsf{L}}
\newcommand\Lap{\mathfrak A}

\newcommand\Ms{M}
\newcommand\N{\mathbb N}
\newcommand\nnu{\mathsf{n}}
\newcommand\olA{\,\overline{\Lap }}
\newcommand\om{\varpi}
\newcommand\qq{\mathsf{q}}

\newcommand\pp{\mathsf{p}}
\newcommand\R{\mathbb{R}}
\newcommand\RE{\operatorname{\text{\sl Re}}}
\newcommand\scs{\scriptstyle}
\newcommand\si{\sigma}

\newcommand\T{\mathbb T}
\newcommand\uf{\mathfrak{u}}
\newcommand\uno{\mathbf{1}}
\newcommand\uH{\underline{H}}
\newcommand\uX{\underline{X}}
\newcommand\ux{\underline{x}}
\newcommand\uxi{\underline{\xi}}
\newcommand\vf{\mathfrak{v}}

\newcommand\wh{\widehat}

\newcommand\Z{\mathbb Z}

\numberwithin{equation}{section}

\newtheoremstyle{mythm}
  {9pt}
  {9pt}
  {\itshape}
  {0pt}
  {\bfseries}
  {}
  { }
  {\thmnumber{(#2)}\thmname{ #1}\thmnote{ #3}}

\newtheoremstyle{mydef}
  {9pt}
  {9pt}
  {\normalfont}
  {0pt}
  {\bfseries}
  {}
  { }
  {\thmnumber{(#2)}\thmname{ #1}\thmnote{ #3}}

\theoremstyle{mythm}
\newtheorem{thm}[equation]{Theorem.}
\newtheorem{pro}[equation]{Proposition.}
\newtheorem{lem}[equation]{Lemma.}
\newtheorem{cor}[equation]{Corollary.}

\theoremstyle{mydef}
\newtheorem{dfn}[equation]{Definition.}
\newtheorem{exa}[equation]{Example.}
\newtheorem{exas}[equation]{Examples.}

\newtheorem{rmk}[equation]{Remark.}

\newtheorem{rmks}[equation]{Remarks.}

\begin{document}$\,$ \vspace{-1truecm}
\title{The heat semigroup and Brownian motion\\ on strip complexes}
\author{\bf Alexander BENDIKOV, Laurent SALOFF-COSTE \\
Maura SALVATORI, and Wolfgang WOESS}
\address{\parbox{.8\linewidth}{Insitute of Mathematics, Wroclaw University\\
Pl. Grundwaldzki 2/4, 50-384 Wroclaw, Poland\\}}
\email{bendikov@math.uni.wroc.pl}
\address{\parbox{.8\linewidth}{Department of Mathematics, Cornell University\\
Ithaca, NY 14853, USA\\}}
\email{lsc@paris.math.cornell.edu}
\address{\parbox{.8\linewidth}{Dipartimento di Matematica, Universit\`a di Milano\\
Via Saldini 50, 20133 Milano, Italy\\}}
\email{maura.salvatori@unimi.it}
\address{\parbox{.8\linewidth}{Institut f\"ur Mathematische Strukturtheorie,
Technische Universit\"at Graz,
Steyrergasse 30, A-8010 Graz, Austria\\}}
\email{woess@TUGraz.at}
\date{\today}
\thanks{S. Bendikov was partially supported by SFB 701 at the 
University of Bielefeld. L. Saloff-Coste was partially supported by 
NSF grant DMS--063886.
W. Woess was partially supported by
the Autrian Science Fund (FWF) project no. P19115-N18.}
\begin{abstract}
We introduce the notion of strip complex. A strip complex
is a special type of
complex obtained by gluing ``strips'' along their natural boundaries
according to a given graph structure. The most familiar example is
the one dimensional complex classically associated with a graph, 
in which case
the strips are simply copies of the unit interval (our setup actually
allows for variable edge length). A leading key  example
is treebolic space, a geometric object studied in a number of recent 
articles, which arises as a horocyclic product of a metric 
tree with the hyperbolic plane.  In this case, the graph is a regular  tree, the strips are 
$[0\,,\,1]\times \R$, and each strip is  equipped with the hyperbolic geometry of a specific strip in
upper half plane.
We consider natural families of Dirichlet forms on a general strip complex 
and show that the associated heat kernels and harmonic functions
have very strong smoothness 
properties. We study questions such as essential selfadjointness of 
the underlying differential operator acting on a suitable space of
smooth functions satisfying a Kirchoff type condition at points where the
strip complex bifurcates. Compatibility with projections that arise from
proper group actions is also considered. 
\end{abstract}

\maketitle

\setcounter{tocdepth}{1}
\tableofcontents

\maketitle

\markboth{{\sf A. Bendikov, L. Saloff-Coste, M. Salvatori, and W. Woess}}
{{\sf Strip complexes}}
\baselineskip 15pt

\section{Introduction}\label{intro}

\subsection*{A. The treebolic spaces $\HT(\pp,\qq)\,$}
Let $\Hb = \{ x + \im y : x \in \R\,,\; y > 0 \}$
be the hyperbolic upper half space, and
$\T=\T_{\pp}$ be the homogeneous tree with degree $\pp+1$, where $\pp \in \N$.
The \emph{treebolic space} is a Riemannian 2-complex which can be viewed as  a
\emph{horocyclic product} of  $\Hb$ and $\T$. Let us start with a picture
and an informal description.
$$
\beginpicture

\setcoordinatesystem units <1.3mm,1.3mm>

\setplotarea x from -16 to 54, y from 0 to 48

\plot 0 0  40 20  48 36  8 16  0 0  -8 16  -2 28  38 48  36.666 45.333 /

\plot 48 36  54 48  14 28  8 16  2  28  42 48   44.4 43.2 /

\plot 4.8 22.4  -8 16  -14 28  26 48  28.4  43.2 /

\put{$\leftarrow$ copies of $S_{k-1}$} [l] at 46 28
\put{$\leftarrow$ copies of $S_{k}$} [l] at 53 42

\endpicture
$$
\begin{center}
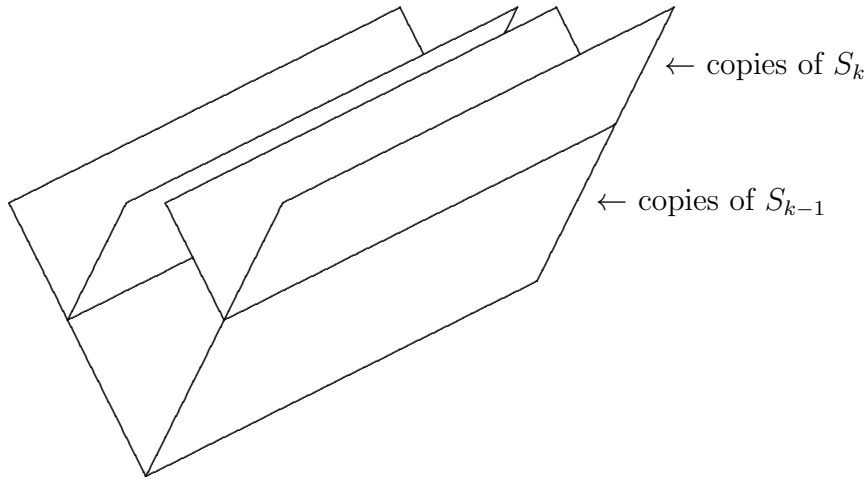
{\sl Figure 1.} A finite portion of treebolic space, with $\pp=2$.
\end{center}

\medskip

Let $1 < \qq\in \R$. Subdivide $\Hb$ into the strips
$S_k = \{ x+\im y : x \in \R\,,\; \qq^{k-1} \le y \le \qq^{k} \}$,
where $k \in \Z$.
Each strip is bounded by two horizontal lines of the form
$L_k= \{ x+\im \qq^k : x \in \R \}$, which, in hyperbolic geometry, are
horospheres with respect to the boundary point at $\infty$ (or rather $\im\infty$).
In the treebolic space $\HT(\pp,\qq)$, infinitely many copies of those strips
are glued together in a tree-like fashion: for each $k \in \Z$, the bottom lines of
$\pp$ copies of $S_{k}$ are identified among each other and with the top line
of $S_{k-1}$. Each strip is equipped with the standard hyperbolic length
element and, in this way,  one obtains a natural metric on $\HT(\pp,\qq)$
as well as a natural measure.

This space admits interesting isometric group actions.
On the one hand, when $\qq=\pp$, the amenable Baumslag-Solitar group
$\BS(\pp) = \langle a,b \mid ab = b^{\pp}a \rangle$ acts on $\HT(\pp,\pp)$
by isometries and with compact quotient. This fact has been exploited by
{\sc Farb and Mosher}~\cite{FaMo} in order to classify
the Baumslag-Solitar groups up to quasi-isometry. See also the nice picture
in {\sc Meier}~\cite[p. 118]{Me}.
On the other hand, for $\pp \ne \qq$, no discrete group can act in such
a way on $\HT(\pp,\qq)$ and its
isometry group is a non-unimodular locally compact group. 
This isometry group admits various subgroups that act with compact quotients,
see our forthcoming paper \cite{BSSW}.

This article is motivated by the following questions. What is Brownian
motion on the treebolic space $\HT(\pp,\qq)$~?
What is the concrete description of the Laplacian, i.e., the
generator of Brownian motion? Can one prove some  essential self-adjoindness
results for this Laplacian? How smooth is the associated heat kernel?
Can one describe explicitly the cone of positive harmonic functions?
The last question, which is at the origin of this work, will be discussed in 
detail in \cite{BSSW}. 
Answers to the other questions are described in theorems 
\ref{thm-sg}-\ref{thm-saHT}.

\bigskip

\subsection*{B. General strip complexes}
The treebolic spaces $\HT(\pp,\qq)$ form one family of examples of what we 
call a \emph{strip complex}, and this work is devoted to the study of 
the heat equation and heat kernel on strip complexes. The  simplest family 
of strip complexes are metric graphs (``quantum graphs''). 
In fact, as a topological space, 
a strip complex is simply the direct product of a (connected) metric graph and
a topological space $\Ms$, e.g., $\{0\}$, $\R$, or a fixed manifold. 
In particular, strip complexes are typically not smooth as they 
bifurcate along the \emph{bifurcation manifolds} 
at the  vertices of the underlying graph structure. See, e.g.,  Figure 1.
We will equip those spaces with certain adapted geometries and adapted  
measures which will give rise to specific Laplacians and heat semigroups. 
Our aim is to show that, because of the specific structure of strip complexes, 
harmonic functions and solutions of the heat equation on such spaces have 
very strong global smoothness properties. Namely, these solutions have locally 
bounded derivatives  of all orders \emph{up to} the bifurcation manifolds 
even though these derivatives are typically not continuous \emph{across} 
the bifurcation manifolds.

In order to carry this out in spite of the singularities of the underlying 
strip complex structure, we build the theory ``from scratch'', using
the theory of strictly local regular Dirichlet forms. See, e.g., 
{\sc Fukushima, Oshima and Takeda}~\cite[Cor.1.3.1]{FOT}
and {\sc Sturm}~\cite{Stu1}, \cite{Stu2}, \cite{Stu3}. The Laplace operators
constructed by this approach are somewhat esoteric objects and one of our goals 
is to describe them in a more concrete way as the closure of operators that
are classical second order 
elliptic differential operators in the smooth part of the complex  
and whose domains of 
definition involve Kirchhoff type laws along  bifurcation manifolds. 

Our material and results should be compared with some previous work.
First, the theory of the Laplacian, heat kernel, etc., on metric graphs is 
quite well understood. See, e.g., {\sc Baxter and Chacon}~\cite{BaCh},
{\sc Cattaneo}~\cite{Cat}, {\sc Enriquez and Kifer}~\cite{EnKi} and
{\sc Kuchment}~\cite{Ku1}, \cite{Ku2}. Note however that, 
even in this simple setting, the exact smoothness of the heat kernel 
is not entirely understood. See {\sc Bendikov and 
Saloff-Coste}~\cite{BSCg}.

Second, {\sc Brin and Kifer}~\cite{BK} introduced Brownian motion
on $2$-dimensional Euclidean complexes
(strongly connected simplicial complexes, where each simplex carries the
Euclidean structure) via a local probabilistic construction. 
The Dirichlet form approach on more general Riemannian
complexes is discussed by {\sc Eells and Fuglede}~\cite{EF} and 
{\sc Pivarski and Saloff-Coste}~\cite{PSC}. None of these references 
provide the type of regularity results proved below for strip complexes.

It is worth emphasizing that, despite the existence of very many different approaches to
the definition of Brownian motion on complexes such as $\HT(\pp,\qq)$, 
the basic problem of uniqueness is not adequately discussed in the literature. 
From this perspective, we view Theorem \ref{thm-saHT} (and its much more general version
Theorem \ref{th-SA}) as an important result.

Many of our results are local in nature. We note that, locally, the simplest 
strip complex structure (a star of finitely many Euclidean half spaces, 
glued along their boundaries) is the model for the neighbourhood
of any generic
singular point in a general $n$-dimensional Euclidean polytopal complex,
that is, any point $\xi$ where the $n$-dimensional closed faces 
containing $\xi$ meet along an $(n-1)$-face. 
The strong regularity results that we obtain thus apply to
small neighbourhoods of such points in any Euclidean polytopal complex.

This paper is organized as follows. 
In Section \ref{geometry}, we exhibit our main  results in  
the key  example of the treebolic space. We describe
a two-parameter family of Dirichlet forms on $\HT$ whose associated Laplacians
and heat semigroups satisfy all regularity and smoothness properties that
one would  wish to have (Theorem \ref{thm-sg}). In each case, the Laplacian is 
the unique self-adjoint extension of a naturally defined, essentially 
self-adjoint operator that is elliptic inside the strips of $\HT$ and acts on 
a space of smooth functions which satisfy a Kirchhoff condition along the
bifurcation lines in $\HT$ (Theorem \ref{thm-saHT}). To the best of our 
knowledge, this is the first time that essential self-adjointness is discussed 
is such a setting.  This construction gives rise to
a Hunt process (``Brownian motion'') on $\HT$ with  natural
projections from $\HT$ onto the underlying (metric) tree and 
onto the hyperbolic
plane (``sliced'' into a strip complex by the lines $L_v$). On each of those
objects, there is a corresponding Dirichlet form and associated 
Laplacian which is the infinitesimal generator of the respective 
projection of the process on $\HT$ (Theorem \ref{thm-projections}).
Uniqueness properties are used here to identify 
the projections with the natural processes
intrinsically defined on the quotient spaces.

In Section \ref{strip}, we introduce 
the notion of strip complex as the product of a metric graph with
a manifold. In a series of definitions, we introduce several function spaces
that are needed to do analysis on such a complex. The geometry of a strip 
complex is obtained through the following data: a length function describing 
the length of the edges of the graph, a Riemannian structure on the 
manifold $\Ms$, and a positive function $\phi$ on the metric graph 
that serves as a conformal factor to define the metric on each strip.
We also introduce a second positive function $\psi$ on the metric graph 
that serves as a weight function to define the underlying measure. 
These data turn the strip complex from a topological space into a geodesic
metric measure space. This structure is used to define a Dirichlet form
whose basic properties are discussed
(Theorems \ref{thm-strilo}--\ref{thm-dense}). This Dirichlet form gives rise
to the associated Laplacian, harmonic functions and heat equation. 

Basic properties of the heat semigroup are  derived in Section \ref{heat}.
Crucial geometric-analytic ingredients are  the local doubling property and
local Poincar\'e inequality (Theorem \ref{th-locDP}). Via the work of
{\sc Sturm}~\cite{Stu1}, \cite{Stu2}, \cite{Stu3} and 
{\sc Saloff-Coste}~\cite{SalB}, this has far reaching consequences for
weak solutions of the heat equation and for the heat diffusion semigroup
(Theorems \ref{th-locHR}--\ref{th-UR} plus corollaries). 

In Section \ref{smoothweak}, we consider weak solutions of the
Laplace and heat equations. We show that these weak solutions are smooth 
up to (but not across) the bifurcation manifolds and satisfy
Kirchhoff type \emph{bifurcation conditions} (Theorems \ref{th-Hsm}
, \ref{th-Heatsm} and  \ref{th-HKsmth}). These results are the most 
significant technical results contained in the present paper. 
 
Section \ref{projs} studies how Dirichlet forms and the associated
heat semigroups are compatible with  natural projections of one strip
complex onto another induced by a proper,
continuous group action (Theorem \ref{th-proj}). 

Uniqueness of the heat semigroup is studied in Section \ref{uniq}. 
First, this question is dealt with on the space of continuous functions that
vanish at infinity, where besides completeness, a uniform local doubling 
property plus uniform local Poincar\'e inequality is needed 
(Theorem \ref{th-UC0}). Second, a very precise  
essential self-adjointness result is obtained
provided completeness and the existence of a \emph{strip-adapted sequence of
functions approximating} $\uno$ (Theorem \ref{th-SA}). The proof of this 
uses in an essential way the heat kernel regularity results proved earlier.
Since we require the existence of an adapted approximation of
$\uno$, this question is briefly dealt with in Section
\ref{adapted}. 

Finally, the appendix contains a hypoellipticity result for the 
operator $\sqrt{-\Delta_{\Ms}}$ on an arbitrary Riemannian manifold which is 
a key element for the proof of the regularity results in Section 
\ref{smoothweak}.

\section{More on $\HT(\pp,\qq)$}\label{geometry}
\subsection*{A. First construction}
We start with a rapid review of some relevant features of the homogeneous
tree $\T=\T_{\pp}\,$. Consider $\T$ as a one-complex, where each edge
is a copy of the
unit interval $[0\,,\,1]$. Let $\T^0$ be the vertex
set ($0$-skeleton) of $\T$. This space is equipped with its natural metric.
A geodesic in $\T$ is the image of an isometric embedding
$t \to w_t \in \T$ of
an interval $I \subset \R\,$.

An \emph{end} of $\T$ is an equivalence class of geodesic rays
(parametrized by $[0\,,\,\infty)$), where
two rays $(w_t)$ and $(\bar w_t)$ are equivalent if they coincide except
perhaps on bounded  initial pieces, i.e., there are $s_0, t_0 \ge 0$ such that
$w_{s_0+t} = \bar w_{t_0+t}$ for all $t \ge 0$.
We write $\bd \T$ for the space of ends, and $\wh \T = \T \cup \bd \T$.
For all $\uf, \vf \in \wh \T$ there is a unique geodesic
$\geo{\uf\,\vf}$ (parametrized by $(-\infty\,,\,\infty)$)
that connects the two.
We choose and fix a reference vertex $o \in \T^0$ and a reference end
$\om \in \bd \T$.
For $v_1,v_2 \in \wh \T \setminus \{ \om \}$, their confluent
$b = v_1 \cf v_2$ with respect to $\om$ is defined by
$\geo{v_1\,\om} \cap \geo{v_2\,\om} = \geo{b\,\om}$.
The \emph{Busemann function} $\hor: \T \to \R$ and the \emph{horocycles} $H_t$
with respect to $\om$ are defined as
$
\hor(w) = d(w,w \cf o) - d(o,w \cf o)$ and
$H_t = \{ w \in \T : \hor(w) = t \}\,.
$
Every horocycle is infinite and denumerable. The vertex set $\T^0$ is the 
union of all $H_k$ with $k \in \Z$. Every vertex $v$ in $H_k$ has one neighbour
$v^-$ (its predecessor) in $H_{k-1}$ and $\pp$ neighbours (its successors)
in $H_{k+1}$.
We set $\bd^* \T = \bd \T \setminus \{\om\}$.

\begin{center}

$$
\beginpicture  \label{fig2}

\setcoordinatesystem units <.7mm,1mm>

\setplotarea x from -10 to 104, y from -84 to -3

\arrow <6pt> [.2,.67] from 0 0 to 80 -80

\plot 32 -32 62 -2 /

 \plot 16 -16 30 -2 /

 \plot 48 -16 34 -2 /

 \plot 8 -8 14 -2 /

 \plot 24 -8 18 -2 /

 \plot 40 -8 46 -2 /

 \plot 56 -8 50 -2 /

 \plot 4 -4 6 -2 /

 \plot 12 -4 10 -2 /

 \plot 20 -4 22 -2 /

 \plot 28 -4 26 -2 /

 \plot 36 -4 38 -2 /

 \plot 44 -4 42 -2 /

 \plot 52 -4 54 -2 /

 \plot 60 -4 58 -2 /



 \plot 99 -29 64 -64 /

 \plot 66 -2 96 -32 /

 \plot 70 -2 68 -4 /

 \plot 74 -2 76 -4 /

 \plot 78 -2 72 -8 /

 \plot 82 -2 88 -8 /

 \plot 86 -2 84 -4 /

 \plot 90 -2 92 -4 /

 \plot 94 -2 80 -16 /


\setdots <3pt>
\putrule from -4.8 -4 to 102 -4
\putrule from -4.5 -8 to 102 -8
\putrule from -1.3 -16 to 102 -16
\putrule from -1.0 -32 to 102 -32
\putrule from  2.2 -42 to 102 -42
\putrule from -1.0 -64 to 102 -64
\setdashes <3pt>
\linethickness =.7pt
\putrule from -3 6 to 102 6
\setlinear

\put {$\vdots$} at 32 3
\put {$\vdots$} at 64 3

\put {$\dots$} [l] at 103 -6
\put {$\dots$} [l] at 103 -48

\put {$H_{-3}$} [l] at -12.5 -64
\put {$H_{-2.3}$} [l] at -12.5 -42
\put {$H_{-2}$} [l] at -12.5 -32
\put {$H_{-1}$} [l] at -12.5 -16
\put {$H_0$} [l] at -12.5 -8
\put {$H_1$} [l] at -12.5 -4

\put {$\vdots$} at -10 3
\put {$\vdots$} [B] at -10 -70
\put {$o$} [rt] at 7.6 -8.6
\put {$\om$} at 82 -82

\put {$\bd^*\T$} [l] at -14 6

\put {\scriptsize $\bullet$} at 8 -8

\endpicture
$$

\vspace{.1cm}

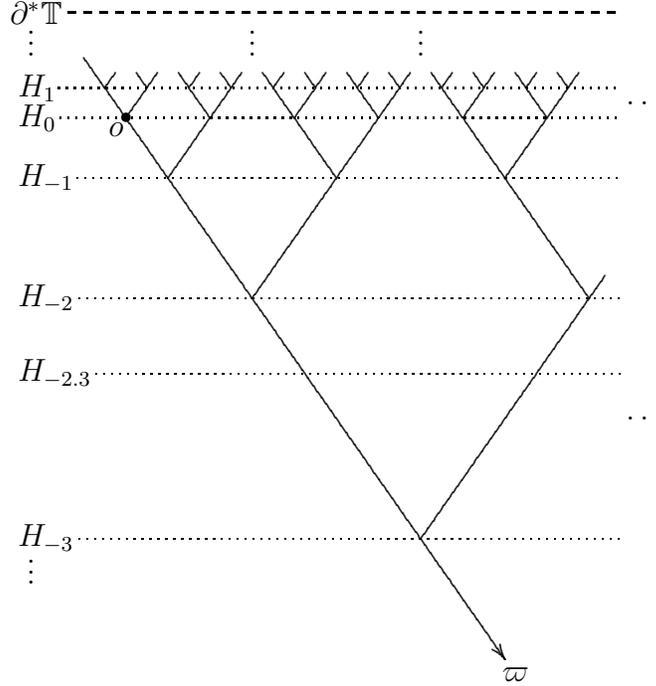
{\sl Figure 2.} The ``upper half plane'' drawing of $\T_2$\\ (top down,
edge lengths are not meaningful in this picture)

\end{center}

\medskip

Fix $\qq>1$ and consider the hyperbolic plane $\Hb$ in its upper-half 
space representation. The horocycles (with respect $\im \infty$) are 
horizontal lines. Recall that  $\T$ is subdivided horizontally by the 
horocycles $H_k$, $k \in \Z$. Similarly, subdivide $\Hb$ in the horizontal 
strips $S_k$ delimited by the lines $y=\qq^k$, 
$k \in \Z$, see Figure~3. Note that all $S_k$ are hyperbolically
isometric.

\medskip

As outlined in the Introduction, the treebolic space with parameters $\qq$ and
$\pp$ is
\begin{equation}\label{eq:treebolicdef}
\HT(\qq,\pp) = \{ (z,w) \in \Hb \times \T_{\pp} :
\hor(w) = \log_{\qq}(\IM z)
\}\,,
\end{equation}
where $\IM z$ is the imaginary part of $z$. Thus, Figures~2 and~3 are the
``side'' and ``front'' views of $\HT$, that is, the images of $\HT$
under the projections $\pi_{\T}: (z,w) \mapsto w$ and 
$\pi_{\Hb}: (z,w) \mapsto z$, respectively.

For each end $\uf \in \bd^*\T$, treebolic space contains the isometric copy
$$
\Hb_{\uf} = \{ (z,w) \in \Hb\times \T_{\pp} : \hor(w) = \log_{\qq}(\IM z)\,,\;
w \in \geo{\uf\,\om}\}
$$
of $\Hb$, and if $\uf  ,\vf \in \bd^*\T$ are distinct and
$v=\uf \cf \vf$ (a vertex),
then $\Hb_{\uf}$ and $\Hb_{\vf}$ bifurcate along the line
$$
\Lf_v = \{ (z,v) \in \Hb\times \T_{\pp} : \IM z = \qq^{\hor(v)} \} 
= \R \times \{v\},
$$ that is,
$\Hb_{\uf} \cap \Hb_{\vf} = \{ (z,w) \in \HT : w \in \geo{v\,\om}\}$.
The metric of $\HT$ is induced by the hyperbolic length element in
the interior of each $\Hb_{\uf}$.

$$
\beginpicture

\setcoordinatesystem units <1.4mm,1.00mm>

\setplotarea x from -40 to 40, y from 0 to 80

\arrow <6pt> [.2,.67] from 0 0 to 0 77.5



\plot -40 4  40 4 /

\plot -40 8  40 8 /

\plot -40 16  40 16 /

\plot -40 32  40 32 /

\plot -40 64  40 64 /

\put {$\im$} [rb] at -0.3 8.6
\put {$\scs \bullet$} at 0 8




\put {$y = \qq^{-1}$} [r] at -40 4
\put {$y = 1$} [r] at -43 8
\put {$y = \qq$} [r] at -43 16
\put {$y = \qq^2$} [r] at -42 32
\put {$y = \qq^3$} [r] at -42 64

\put {$\infty$} [t] at 0 80

\put {$\R$} [r] at -45 0

\put {$S_1$} at 30 12
\put {$S_2$} at 30 24
\put {$S_3$} at 30 48

\setdashes <3pt>
\linethickness =.7pt
\putrule from -43.4 0 to 40.4 0
\setlinear

\endpicture
$$
\vspace{.1cm}

\begin{center}
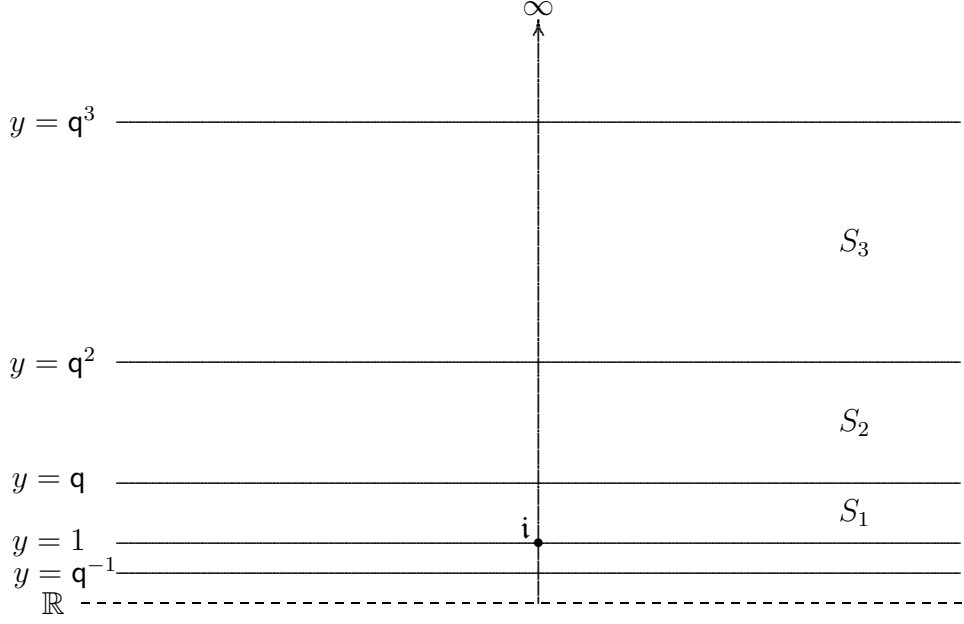
{\sl Figure 3.} Hyperbolic upper half plane $\Hb$ subdivided in isometric strips
\end{center}

\medskip

\subsection*{B. Second construction}
We now present an alternative construction of $\HT=\HT(\pp,\qq)$ which
leads to further generalizations. It is clear that, as a topological space,
$\HT$ is simply
$$
\HT=\T_{\pp}\times \R\,.
$$
Note that topologically, $\qq$ plays no role.
Now, let us view $\T_{\pp}$ as a metric tree $\T_{\pp,\qq}$
by setting the length of all edges between the horocycles
$H_{k-1}$ and $H_{k}$ to be $\qq^{k-1}(\qq-1)$.
Hence,  $\T_{\pp,\qq}\times \R$ comes equipped with a natural geometry.
Namely, given any edge  $e=[v^-,v]$, parametrized by $s\in [q^{k-1}\,,\,q^{k}]$,
$k=\hor(v)$, we can view
$[v^-,v]\times \R$ as a manifold with global coordinates
$(s,x)\in [q^{k-1}\,,\,q^{k}]\times \R\,$.
We can equip this manifold with the length element 
$s^{-2}\bigl((ds)^2+(dx)^2\bigr)$.
Doing this for all edges yields a new metric structure on $\HT$ which is
isometric to its treebolic structure described earlier. Indeed, any
doubly infinite geodesic joining $\om$ to another end of $\T$ determines
an upper-half plane in $\T_{\pp,\qq}\times \R\,$,
and the construction outlined  above yields the hyperbolic metric on  
any of these upper-half planes (with $s=y$, $z=x+\im y$).
The natural measure on $\T_{\pp,\qq}\times \R$ is given on
a strip $[v^-,v]\times \R\,$, viewed as a manifold with global coordinates
$(s,x)\in [q^{k-1},q^{k}]\times \R\,$, by $s^{-2}\,ds\,dx\,$.

\bigskip

\subsection*{C. The two parameters family of Dirichlet forms $\mathcal
E_{\alpha,\beta}\,$}
Recall that the Riemannian metric and 
measure of the
hyperbolic plane  $\Hb={\R}^2_+$ (upper half plane model)
are  given by $y^{-2}(dx^2+dy^2)$ and
$d\mu=y^{-2}\,dx\,dy\,$, respectively. The natural Dirichlet form on $\Hb$ is
$$
\int_{\Hb}|\nabla f|^2\,d\mu =
\int_{\Hb} (|\partial_xf|^2+|\partial_y f|^2) \,dx\,dy\,.
$$
The Laplacian is $y^2(\partial_x^2+\partial_y^2)$.
See, e.g., {\sc Chavel}~\cite[p. 263--265]{Cha}.

Any element $\xi$ in $\HT$ is described uniquely by
a pair $(z,v)$ with $v\in \T^0$ and $z=x+\im y \in \Hb$ with $x\in \R$,
$\qq^{k-1}<y\le \qq^k$ and $k=\hor(v)$. In this case, we write
$y = y(\xi)$ and $v=v(\xi)$.

Thus, for each $v\in \T^0$, we consider
\begin{eqnarray*}
S_v&=&\{(z,v):z=x+\im y\in \Hb\,,\; x\in \R\,,\;\qq^{k-1}\le y\le \qq^k\}\\
S^o_v&=&\{(z,v):z=x+\im y\in \Hb\,, \; x\in \R\,,\;
\qq^{k-1}<y\ <\qq^k\}
\end{eqnarray*}
where $k=\hor(v)$.
The lines
$$
L_v=\{(z,v):z=x+\im \qq^{\hor(v)}, x\in \R\}
$$
are called bifurcation lines. With this notation,
we have
$$
\HT=\bigcup_{v\in \T^0} (S_v \setminus L_{v^-})
\qquad\text{(a disjoint union)}.
$$
Note that all the strips $S^o_v$ are isometric and have hyperbolic 
width $\log \qq$. However, above we have kept the Euclidean coordinates, 
taking into account the ``height'' of the strip $S_v\,$, i.e., $k=\hor(v)\,$.

As mentioned, the space $\HT$ carries a natural measure 
(again coming from $\Hb$) that we denote by  $d\xi$.
Namely,
\begin{equation}\label{eq-dxi}
\int_{\HT}f(\xi)\,d\xi=\sum_{v\in \T^0}\,\int_{S^o_v}f(x+\im y,v)\,y^{-2}\,dx\,dy\,.
\end{equation}
For $\alpha\in \R\,$, $\beta>0$, set
\begin{equation}\label{muab}
d\mu_{\alpha,\beta}(\xi)= \beta^{\hor(v)}\,y^\alpha\, d\xi
=\beta^{\hor(v)}\,y^{\alpha-2}\, dx\,dy\,.
\end{equation}
This means that
\begin{equation}\label{eq-intHT}
\int_{\HT}f(\xi)\,d\mu_{\alpha,\beta}(\xi)
=\sum_{v\in \T^0}\beta^{\hor(v)}\int_{S^o_v}f(x+\im y,v)\,y^{-2+\alpha}\,dx\,dy\,.
\end{equation}
For any open strip $S^o_v$ equipped with the $(x,y)$-coordinates as above,
let $\mathcal W^1(S^o_v)$ be the Sobolev space of those functions
$f$ in ${\mathcal L}^2(S^o_v)$ whose distributional first order partial derivatives
$\partial_xf,\partial_y f$ can be represented by functions in ${\mathcal L}^2(S^o_v)$
(with respect to the measure $dx\,dy$, say). By a fundamental theorem 
concerning Sobolev spaces, such functions admit a trace 
$\mbox{Tr}^{S^o_v}_{L}(f)$ on
each of the lines bordering the strip. This trace is in fact in the
fractional Sobolev space $\mathcal W^{1/2}(L)$ of the lines $L$. Namely,
the trace theorem asserts that $\mbox{Tr}^{S^o_v}_L$ defined on
$\mathcal C^\infty(S_v)$
extends as a bounded operator
$$
\mbox{Tr}^{S^o_v}_L: \mathcal W^1(S^o_v)\rightarrow \mathcal W^{1/2}(L).
$$
We can now describe a two parameters family of function spaces and
Dirichlet forms on $\HT$ which all share the same underlying geometry.

\begin{dfn}\label{def-Dirforms}
Fix $\alpha\in \R\,$, $\beta>0$. Let $\Omega$ be an open set in $\HT$.
We define $\mathcal W^1_{\alpha,\beta}(\Omega)$ as
the space of all functions $f$ in ${\mathcal L}^2(\Omega,\mu_{\alpha,\beta})$ such that the
following two properties hold.
\begin{enumerate}
\item For each $v\in \T^0$, the function $f$, restricted to 
$S^o_v\cap \Omega\,$, is in $\mathcal W^1(S^o_v\cap \Omega)$, and
\begin{eqnarray*}
\|f\|_{\mathcal W^1_{\alpha,\beta}(\Omega)}^2&=&\sum_{v\in \T^0}\beta^{\hor(v)}
\int_{S^o_v\cap \Omega}
\Bigl(|f(z,v)|^2\,y^{-2}+|\partial_x f(z,v)|^2+|\partial_y f(z,v)|^2\Bigr) 
\,y^\alpha\,dx\,dy \\
&=& \int_{\Omega}\Bigl(|f(\xi)|^2 + |\nabla f(\xi)|^2\Bigr) \,d\mu_{\alpha,\beta}(\xi)
<\infty\,,
\end{eqnarray*}
where, for $\xi=(z,v)$, we have set
$\nabla f(\xi)=\bigl(y^2 \,\partial_x f(z,v), y^2\,\partial_y f(z,v)\bigr)$ 
and
$$
|\nabla f (\xi)|^2= \bigl\langle \nabla f(\xi),\nabla f(\xi)\bigr\rangle_{z}=
y^2\bigl(|\partial_x f(z,v)|^2+|\partial_y f(z,v)|^2\bigr).
$$
(The inner product is with respect to the hyperbolic metric in the
$z$-variable.)
\item For any pair of neighbours $u,v\in \T^0$ such that 
$S_v\cap S_u=L$, one has 
$\mbox{Tr}^{S^o_v}_Lf=\mbox{Tr}^{S^o_u}_Lf$
along $L\cap \Omega$.
\end{enumerate}
Let $\mathcal W^1_{\alpha,\beta,0}(\Omega)$ be the completion of
$\mathcal W^1_{\alpha,\beta}(\Omega) \cap \mathcal C_c(\Omega)$ with respect to the 
norm  $\|\cdot\|_{W^1_{\alpha,\beta}(\Omega)}\,$.
\end{dfn}

\begin{dfn}\label{def-Eab}
Let $\mathcal E_{\alpha,\beta}$ be the bilinear form
\begin{eqnarray}
{\mathcal E}_{\alpha,\beta}(f,g)&=&\sum_{v\in \T^0}\beta^{\hor(v)}
\int_{S^o_v}
\Bigl(\partial_xf(z,v)\,\partial_xg(z,v)+
\partial_yf(z,v)\,\partial_yg(z,v)\Bigr)\,y^\alpha \,dx\,dy \nonumber \\
&=&
\int_{\HT} \bigl\langle  \nabla f(\xi),\nabla g(\xi)\bigr\rangle_{z(\xi)}
\,\, d\mu_{\alpha,\beta}(\xi)\,.
\end{eqnarray}
with domain $\mathcal{D}(\mathcal{E}_{\alpha,\beta})=
\mathcal W^1_{\alpha,\beta}(\HT)\subset {\mathcal L}^2(\HT,\mu_{\alpha,\beta})$.
Here, $z(\xi)=z$ if $\xi =(z,v)\in \HT$.
\end{dfn}

Note that for $f\in \mathcal W^1_{\alpha,\beta}(\HT)$,
the function $\xi\mapsto |\nabla f(\xi)|$ is well defined
as an element of ${\mathcal L}^2(\HT)$. In the present context, $|\nabla f|^2$
is the \emph{carr\'e du champ,} also often denoted by
$$
|\nabla f|^2=\Gamma(f,f)=
\frac{d\,\Gamma_{\alpha,\beta}(f,f)}{d\mu_{\alpha,\beta}}\,,
$$
where $d\,\Gamma_{\alpha,\beta}(f,f)$ is the \emph{energy measure}
associated to $f\in \mathcal W^1_{\alpha,\beta}(\HT)$.
Observe that the \emph{carr\'e du champ} does not depend on the parameters
$\alpha,\beta$. This explains why we say that these Dirichlet forms all 
share the same geometry.

\begin{dfn} We let  $\mathcal {C}^\infty(\HT)$ be the set of those continuous
functions $f$ on $\HT$ such that, for each $v\in \T^0$, the restriction
$f_v=f(\cdot,v)$ of $f$ to the  
closed strip $S_v$ 
has continous derivatives
$\partial_x^m\partial_y^n f(z,v)$ of all orders in the interior $S^o_v$
which satisfy, for all $R>0$,
$$
\sup\bigl\{|\partial_x ^m\partial_y^n f(z,v)|:
(z,v) \in S^o_v\,,\; |\RE z| \le R\bigr\}
<\infty\,.
$$
Given an open set $\Omega\subset \HT$, we let $\mathcal C^\infty_c(\Omega)$
be the space of those functions in $\mathcal C^\infty(\HT)$ that have compact 
support in $\Omega\,$.
\end{dfn}
\begin{rmk}\label{rem-derivs}
The condition implies that each partial derivative 
$\partial_x^m\partial_y^n f(z,v)$ extends continuously to the boundary
of $S_v$. We write $\partial_x^m\partial_y^n f_v$ for this extension. 

Note however that only the function $f\in \mathcal {C}^\infty(\HT)$
itself has to be continuous at the bifurcation lines, not its derivatives.
That is, if $w^-=v$ then it is in general \emph{not} true that 
$\partial_x^m\partial_y^n f_w = \partial_x^m\partial_y^n f_v$
on $L_v = S_v \cap S_w\,$, unless $m=n=0$.
\end{rmk}

\begin{pro} 
For each $\alpha\in \R$ and $\beta>0$,
the form $\bigl(\mathcal E_{\alpha,\beta}\,,\mathcal W^1_{\alpha,\beta}(\HT)\bigr)$
is a strictly local regular Dirichlet form, and  $\mathcal C^\infty_c(\HT)$
is a core for this Dirichlet form.

For any open set $\Omega$, the space 
$\mathcal C^\infty_c(\Omega)$ is dense in $\mathcal W^1_{\alpha,\beta,0}(\Omega)$.
\end{pro}

Note that the regularity of these Dirichlet forms is not obvious at all.
We will prove this result in a more general setting below.

\bigskip

\subsection*{D. The heat semigroup and Brownian motion}
For each $\alpha\in \R$, $\beta>0$,  the Dirichlet form
$\bigl(\mathcal E_{\alpha,\beta}\,,\mathcal W^1_{\alpha,\beta}(\HT)\bigr)$
induces  a self-adjoint contraction semigroup $e^{t\Delta_{\alpha,\beta}}$
with infinitesimal generator (``Laplacian'') $\Delta_{\alpha,\beta}$
on ${\mathcal L}^2(\HT,\mu_{\alpha,\beta})$.  The domain
$\mbox{Dom}(\Delta_{\alpha,\beta})$ of $\Delta_{\alpha,\beta}$ is the
set of functions $f\in \mathcal W^1_{\alpha,\beta}(\HT)$ for which there exists
a constant $C_f$ such that
$$
\mathcal E_{\alpha,\beta}(f,g)=
\int_{\HT} \bigl\langle \nabla f(\xi),\nabla g(\xi)\bigr\rangle_{z(\xi)}\,
d\mu_{\alpha,\beta}(\xi)
\le C_f \, \|g\|_{{\mathcal L}^2(\HT,\mu_{\alpha,\beta})}
$$
for all $g\in \mathcal W^1_{\alpha,\beta}(\HT)$. As $\mathcal W^1_{\alpha,\beta}(\HT)$
is dense in ${\mathcal L}^2(\HT,d\mu_{\alpha,\beta})$,
this condition and the Riesz representation theorem imply that there exists 
a (unique) function $h\in {\mathcal L}^2(\HT,d\mu_{\alpha,\beta})$ such that 
$\mathcal E_{\alpha,\beta}(f,g)=-\int_{\HT} h\,g\,d\mu_{\alpha,\beta}\,$.
By definition, $\Delta_{\alpha,\beta}f=h\,$\; see, e.g., 
\cite[Cor.1.3.1]{FOT}. If $f$ is in
$\mbox{Dom}(\Delta_{\alpha,\beta})\cap C^\infty(\HT)$ then,
in each open strip,
\begin{equation}\label{lap}
\Delta_{\alpha,\beta}f 
= \bigl[y^2(\partial_x^2+\partial_y^2)+ \alpha y\partial_y \bigr]\,f\,,
\end{equation}
but $f$ must also satisfy
the \emph{bifurcation} or \emph{Kirchhoff condition}
\begin{equation}\label{bif}
\partial_y f_v=\beta\sum_{w\,:\,  w^- = v}\partial_y f_w \quad
\text{on $L_v$ for each $v \in \T^0\,$.}
\end{equation}
Note that the parameter $\beta$ comes into play only at the bifurcation lines
where it appears in the bifurcation condition \eqref{bif}
relating the different vertical partial derivatives in the $\pp+1$ strips 
meeting along any given bifurcation line. This will be discussed in detail later
on.

\begin{thm} \label{thm-sg}
The semigroup $e^{t\Delta_{\alpha,\beta}}$, $t>0$,
acting on ${\mathcal L}^2(\HT,\mu_{\alpha,\beta})$ has the following properties:
\begin{enumerate}
\item[(a)] It admits a continuous positive symmetric transition kernel
$$
(0,\infty)\times \HT\times\HT\ni (t,\xi,\zeta)\mapsto
h_{\alpha,\beta}(t,\xi,\zeta)
$$
such that for all $f\in \mathcal C_c(\HT)\,$,
$$
e^{t\Delta_{\alpha,\beta}}f(\xi)=\int_{\HT}
h_{\alpha,\beta}(t,\xi,\zeta)\,f(\zeta)\,d\mu_{\alpha,\beta}(\zeta)\,.
$$
\item[(b)]
For each fixed $(t,\xi)$, the function 
$\zeta\mapsto h_{\alpha,\beta}(t,\xi,\zeta)$  is
in $\mathcal C^\infty(\HT)$ and satisfies {\em (\ref{bif})}.
\item[(c)] For each $k\in \N$, the function
$(0,\infty)\times \HT\times\HT\ni (t,\xi,\zeta)\mapsto 
\partial_t^k h_{\alpha,\beta}(t,\xi,\zeta)$ is H\"older continuous,
and for each $\xi\in \HT$, the function
$\zeta\mapsto \partial_t^k h_{\alpha,\beta}(t,\xi,\zeta)$ is 
in $\mathcal C^\infty(\HT)$ and satisfies {\em (\ref{bif})}.
\item[(d)] For any fixed $\epsilon\in (0,1)$ and $k\in \N$, there is a constant
$C=C(\alpha,\beta,\pp,\qq,k,\epsilon)$ such that for
all $(t,\xi,\zeta)\in (0,\infty) \times \HT\times \HT$, 
\begin{equation}\label{heatest'}
|\partial_t^kh_{\alpha,\beta}(t,\xi,\zeta)|
\le
\frac{C}{\beta^{\hor(v(\xi))}\,y(\xi)^\alpha \,\min\{1,t\}\,t^{k}}\,
\exp\left(-\frac{d(\xi,\zeta)^2}{4(1+\epsilon)t}\right).
\end{equation}
\item[(e)] It is conservative, that is, 
$e^{t\Delta_{\alpha,\beta}} \mathbf 1=\mathbf 1$. Equivalently, 
$\int_{\HT} h_{\alpha,\beta}(t,\xi,\cdot)\,d\mu_{\alpha,\beta}=1$.
\item[(f)] It sends ${\mathcal L}^\infty(\HT)$ into $\mathcal C(\HT)\cap {\mathcal L}^\infty(\HT)$.
\item[(g)] It sends $\mathcal C_0(\HT)$ into itself.
\item[(h)] The associated  Hunt process is transient, that is, for all pairs of
distinct points $\xi, \zeta\in \HT$,
$$
G_{\alpha,\beta}(\xi,\zeta)=\int_0^\infty
h_{\alpha,\beta}(t,\xi,\zeta)\,dt<\infty\,.
$$
\item[(i)] The bottom $\lambda=\lambda(\alpha,\beta,\pp,\qq)$
of the ${\mathcal L}^2(\HT,\mu_{\alpha,\beta})$-spectrum  of $-\Delta_{\alpha,\beta}$ is 
strictly positive if and only if $ \qq^{1-\alpha} \ne \beta\pp $.

In particular, in addition to \eqref{heatest'} the following holds.
     
For any fixed $\epsilon\in (0,1)$ and $k\in \N$, there is a constant
$C=C(\alpha,\beta,\pp,\qq,k,\epsilon)$ such that for all 
$(t,\xi,\zeta)\in (0,\infty) \times \HT\times \HT$, 
\begin{equation}\label{heatest2}
|\partial_t^kh_{\alpha,\beta}(t,\xi,\zeta)|\le
\frac{C}{
\beta^{\hor(v(\xi))}y(\xi)^\alpha (\min\{1,t\})^{1+k}}
\exp\left(-\lambda t-\frac{d(\xi,\zeta)^2}{4(1+\epsilon)t}\right).
\end{equation}
\end{enumerate}
\end{thm}
\begin{proof} Statements (a) through (g) follow from more general results
proved in this paper. That $\lambda$ is positive if and only if
$ q^{1-\alpha}/(\beta\pp)\neq 1$ can be obtained by the techniques
and results of {\sc Saloff-Coste and Woess}~\cite{SCW-S} which also provides 
an explicit formula for $\lambda$ in terms of the parameters.
Transience is explained below after Theorem \ref{thm-projections}. 
\end{proof}

\begin{dfn}\label{dfn-ADHT} Let $\HT^o=\bigcup_v S^o_v\,$ be the treebolic space
without the bifurcation lines.  For $f\in \mathcal C^\infty(\HT^o)$, set
$$
\Lap_{\alpha}f(\xi)= y^{2}(\partial_x^2+\partial_y^2)f(\xi)
+\alpha y\partial_yf(\xi),\;\;\xi=(x+\im y, v)\in \HT^o\,.
$$
Let $\mathcal D^\infty_{\alpha,\beta,c}$ be the
space of those functions in $\mathcal C^\infty_c(\HT)$ such that:
\begin{itemize}
\item For any $k$, the function $\Lap_\alpha^k f$, originally defined on
$\HT^o,$  admits a continuous extension to all of $\HT$. (Here, 
$\Lap_\alpha^k$ is the $k$-th iterate of $\Lap_\alpha$.)
This implies that $\Lap_\alpha ^kf \in \mathcal C^\infty_c(\HT)$ for each $k$.
\item Using the same notation as in Remark \ref{rem-derivs} and formula
\eqref{bif},
$$ 
\partial_y \Lap_\alpha^k f_v=\beta\sum_{w\,:\,w^-=v}
\partial_y \Lap_\alpha^k f_w \quad
\text{on $L_v$ for each $v \in \T^0\,$.}
$$
\end{itemize}
\end{dfn}
The following statement
yields a clear and fundamental  uniqueness result concerning the  Laplacian
$\Delta_{\alpha,\beta}$ introduced  above. For the proof, 
see Theorem \ref{th-SA} and Proposition \ref{pro-tbexhaust}.

\begin{thm} \label{thm-saHT}
The operator $(\Lap_\alpha,\mathcal D^\infty_{\alpha,\beta,c})$ is symmetric
on ${\mathcal L}^2(\HT,\mu_{\alpha,\beta})$. It is essentially self-adjoint and its
unique self-adjoint extension is the infinitesimal generator
$\bigl(\Delta_{\alpha,\beta},\mbox{\em Dom}(\Delta_{\alpha,\beta})\bigr)$ 
associated with the Dirichlet form
$\bigl(\mathcal E_{\alpha,\beta}\,,\mathcal W^1_{\alpha,\beta}(\HT)\bigr)$ on
${\mathcal L}^2(\HT,\mu_{\alpha,\beta})\bigr)$.
\end{thm}

\begin{rmk} Let $X$ be a topological space equipped with a Borel measure
$\mu$ with full support. A densely defined operator
$\bigl(\Lap,\mbox{Dom}(\Lap)\bigr)$ on ${\mathcal L}^1(X,\mu)$ is called strongly Markov-unique
if and only if there is at most one sub-Markovian $\mathcal C^0$-semigroup
on ${\mathcal L}^1(X,\mu)$
whose infinitesimal generator extends $\bigl(\Lap,\mbox{Dom}(\Lap)\bigr)$.
It is not hard to see that a symmetric essentially self-adjoint operator
is strongly Markov-unique. See, e.g., {\sc Eberle}~\cite{Eb}.
\end{rmk}

\bigskip

\subsection*{E. The $(\alpha,\beta)$-Markov process}
By the general theory of Markov processes,
there is a Hunt process associated with the
conservative semigroup $H^{\alpha,\beta}_t=e^{t\Delta_{\alpha,\beta}}:
\mathcal C_0(\HT)\mapsto \mathcal C_0(\HT)$. It is defined for every 
starting point $\xi\in\HT$, has infinite life time and continuous sample paths.
Its family of distributions $(\mathbb P^{\alpha,\beta}_\xi )_{\xi\in \HT}$
on $\boldsymbol{\Omega}=\mathcal C ([0,\infty]\to \HT)$ is determined by
the one-dimensional distributions
$$
\mathbb P^{\alpha,\beta}_\xi[X_t\in U]
= \int_U h_{\alpha,\beta}(t,\xi,\zeta) \,d\mu_{\alpha,\beta}(\zeta)
=H^{\alpha,\beta}_t\mathbf 1_U(\xi)$$
where $U$ is any Borel subset of $\HT$.

Setting $\tau_U=\inf\{t: X_t\not\in U\}$,
we  can define the exit distribution from a bounded Borel set $U$ by
$$
\pi^{\alpha,\beta}_U(\xi,B)=\mathbb P^{\alpha,\beta}_\xi[X_{\tau_U}\in B]
$$
for any Borel set $B \subset U$ and set
$$
\pi^{\alpha,\beta}_U(\xi,f)
=\mathbb E^{\alpha,\beta}_\xi\bigl(f(X_{\tau_U})\bigr)
$$
for any bounded Borel measurable function $f$. Since the process has continuous
sample paths, the exit distribution is supported by $\partial U$ for any
starting point $\xi \in U$.

As outlined at the beginning of this section, the treebolic space $
\HT(\qq,\pp) = \{ (z,w) \in \Hb \times \T_{\pp} :\hor(w) = \log_{\qq}(\IM z)\}$
(here written in terms of the first construction)
admits natural projections, $\pi_{\Hb}: (z,w) \mapsto z$ and
$\pi_{\T}: (z,w) \mapsto w$, corresponding
respectively to the
``side'' and ``front'' views of $\HT$ depicted in  Figures~2 and~3.

By the general theory of transformations of the  
state space, it is plain that the images of the Hunt process 
$(X_t,\mathbb P^{\alpha,\beta}_\xi, t\ge 0, \xi\in \HT)$ by the projections 
$\pi_{\Hb}$ and $\pi_{\T}$
are Markov processes. What is not entirely obvious, \emph{a priori,} is to
describe what these processes are in intrinsic terms in $\Hb$ and $\T$.
One of the multiple motivations behind this work was indeed to obtain 
an intrinsic description of each of these processes.

Analogously to $\HT$, we can describe the metric tree
$\T=\T_{\pp,\qq}$ as
$$
\T=\{(s,v): v\in \T^0\,,\;
s\in (\qq^{\hor(v)-1}\,,\,\qq^{\hor(v)}]\}\,,
$$
where $\{ v \} \times (\qq^{\hor(v)-1}\,,\,\qq^{\hor(v)}]$ parametrizes
the ``metric edge'' $(v^-,v]$ as a left-open interval. On $\T$
we consider the measure
$\mu^{\T}_{\alpha,\beta}$ defined by
$d\mu^{\T}_{\alpha,\beta}(s,v)= \beta^{\hor(v)}\,s^{-2+\alpha}\, ds$, that is,
for all $f\in \mathcal C_c(\T)$
\begin{equation}\label{eq-treemeasure} 
\int_{\T}f\, d\mu^{\T}_{\alpha,\beta}
=\sum_{v\in \T^0}\beta^{\hor(v)} \int_{q^{\hor(v)-1}}^{q^{\hor(v)}}
f(s,v)\, s^{-2+\alpha}\, ds\,,
\end{equation}
and the Dirichlet form
\begin{equation}\label{eq-treeform}
\mathcal E^{\T}_{\alpha,\beta}(f,f)=
\int_{\T}s^2\,|\partial_s f|^2\,d\mu^{\T}_{\alpha,\beta}
=\sum_{v\in \T^0}\beta^{\hor(v)}\int_{q^{\hor(v)-1}}^{q^{\hor(v)}}
|\partial_s f(s,v)|^2\,s^{\alpha}\,ds\,,
\end{equation}
with domain
$$
\mathcal W^1_{\alpha,\beta}(\T)=\{
f\in \mathcal C(\T)\cap {\mathcal L}^2(\T,\mu^\T_{\alpha,\beta}): s\,\partial_s f
\in {\mathcal L}^2(\T,\mu^\T_{\alpha,\beta})\}.
$$
Here $\partial_s\, f$ denotes the distributional derivative of $f$
along any open edge $(v^-,v)= \{v\}\times (q^{\hor(v)-1}\,,\,q^{\hor(v)})$ 
of $\T$.
Let $h^\T_{\alpha,\beta}(t,\cdot,\cdot)$, $t>0$, be the heat kernel
associated with this Dirichlet form.

On the hyperbolic space $\Hb$, subdivided by the horocycle lines
$L_k=\{z=x+\im y: y=q^k\}$, consider the measure
$\mu^{\Hb}_{\alpha,\beta}$ which is defined 
for all $f\in \mathcal C_0(\Hb)$ by
\begin{equation}\label{eq-planemeasure}
\int_{\Hb}f\, d\mu^{\Hb}_{\alpha,\beta}
=\sum_{k\in \mathbb Z}\beta^{k}\int_{q^{k-1}} ^{q^{k}}
\int_{-\infty}^\infty
f(x+\im y)\, y^{-2+\alpha}\,dx\,dy\,,
\end{equation}
and the Dirichlet form
\begin{equation}\label{eq-planeform}
\begin{aligned}
\mathcal E^{\Hb}_{\alpha,\beta}(f,f)&=
\int_{\Hb} |\nabla f|^2\,d\mu^{\Hb}_{\alpha,\beta}\\
&=\sum_{k\in \mathbb Z}\beta^{k}\int_{q^{k-1}}^{q^{k}}\int_{-\infty}^\infty
\bigl(|\partial_x f(x+\im y)|^2 + |\partial_y f(x+\im y)|^2\bigr)\, 
y^{\alpha}\,dx\,dy\,,
\end{aligned}
\end{equation}
where $|\nabla f|$ denotes the hyperbolic gradient length of $f$.
The domain of this form is the space $\mathcal W^1_{\alpha,\beta}(\Hb)$
of those functions in ${\mathcal L}^2(\Hb,\mu^\Hb_{\alpha,\beta})$ which admit
locally integrable first order partial  derivatives in the sense
of distributions and such that $|\nabla f|$ is in
${\mathcal L}^2(\Hb,\mu^\Hb_{\alpha,\beta})$. Let $h^\Hb_{\alpha,\beta}(t,\cdot,\cdot)$,
$t>0$, be the heat kernel associated with this Dirichlet form on $\Hb$.
(All this coincides precisely whith what we have considered in the 
previous subsections on $\HT(\pp,\qq)$, but now we are in the ``degenerate''
case when $\pp=1$ and the tree is a two-way-infinite linear graph.)

\begin{thm}\label{thm-projections} Fix $\pp\in \{2,3,\dots\}$, $\qq>1$ and 
$\alpha\in \R$, $\beta>0$. Let $(X_t)$ be the process on $\HT(\pp,\qq)$ 
associated with the Dirichlet form 
$\bigl(\mathcal E_{\alpha,\beta}\,,\mathcal W^1_{\alpha,\beta}(\HT)\bigr)$.
Let $Y_t=\pi^\T(X_t)$, $Z_t=\pi^\Hb(X_t)$, $t>0$, be the projections on
$\T$ and $\Hb$, respectively.
\begin{itemize}
\item[(a)]  The process $(Y_t)$ is a Markov process on $\T$ and, for any $t>0$
and $y\in \T$, the law of  $Y_t$ given $Y_0=y_0$ has probability density
$h^{\T}_{\alpha,\beta}(t,y_0,\cdot)$ with respect to 
$\mu^{\T}_{\alpha,\beta}\,$.

In other words, $(Y_t)$ is  a version of the Hunt process associated with
the strictly local regular Dirichlet form 
$\bigl(\mathcal E^\T_{\alpha,\beta},\mathcal  W^1_{\alpha,\beta}(\T)\bigr)$.
\item[(b)]
The process $(Z_t)$ is a Markov process on $\Hb$ and, for any $t>0$
and $z\in \Hb$, the law of $Z_t$ given $Z_0=z_0$ has probability density
$h^{\Hb}_{\alpha,\beta\pp}(z_0,\cdot)$ with respect to
$\mu^{\Hb}_{\alpha,\beta\pp}\,$.

In other words, $(Z_t)$ is  a version of the Hunt process associated with
the strictly local regular Dirichlet form
$\bigl(\mathcal E^{\Hb}_{\alpha,\beta\pp}, \mathcal W^1_{\alpha,\beta\pp}(\Hb)\bigr)$.
\end{itemize}
\end{thm}

See Proposition \ref{pro-goodact} and Example \ref{ex-goodact}(C) at the end of 
Section \ref{projs}.

\begin{pro}\label{pro-trans}
Each of the processes $(X_t)$, $(Y_t)$ and $(Z_t)$ appearing in Theorem
\ref{thm-projections} is transient.
\end{pro}
 
\begin{proof}
Via the projections, transience of $(X_t)$ will follow from
transience of $(Z_t)$. 

This amounts to showing that for every choice of $\alpha \in \R$ and 
$\beta > 0$, the process on $\HT(1,\qq) = \Hb$ associated with
$\bigl(\mathcal E^{\Hb}_{\alpha,\beta}, \mathcal W^1_{\alpha,\beta}(\Hb)\bigr)$
is transient. Now, the associated measure $\mu^{\Hb}_{\alpha,\beta}\,$
can be compared below and above, up to multiplying with positive constants,
with the measure $\mu^{\Hb}_{\bar\alpha,\bar \beta}\,$,
where $\bar \alpha = \alpha + \log \beta / \log q$ and $\bar \beta = 1$.
Hence, the associated metric measure spaces are (measure) quasi-isometric
(i.e., quasi-isometric with adapted measures, see {\sc Coulhon and
Saloff-Coste}~\cite{CSC}). 
This implies that the corresponding processes are
either both transient or both recurrent. Hence, it thus suffices to study
the transience of the process on $\Hb$ associated with 
$\bigl(\mathcal E^{\Hb}_{\bar \alpha,1}, \mathcal W^1_{\bar \alpha,1}(\Hb)\bigr)$.
This process does not ``see'' the separating lines bounding the strips.
Indeed, the associated infinitesimal generator on the whole upper half plane is
$$
\Delta_{\bar \alpha,1}=y^2(\partial_x^2+\partial _y^2) 
+ \bar \alpha \, y \,\partial_y\,.
$$
The process is just standard hyperbolic Brownian motion on $\Hb$ 
with an additional vertical drift term. 
It is very well known to be transient. For example, one finds
nonconstant positive harmonic functions that are expressed in terms
of the Poisson kernel. Another way is to identify $\Hb$ with the \emph{affine
group} of all transformations $x \mapsto ax+b$, where
$a > 0$ and $b \in \R$, via $(a,b) \leftrightarrow b + \im a \in \Hb$.
Then the law of our process is invariant under the action of the affine
group on itself, whence it must be transient, compare e.g. with
{\sc Guivarc'h, Keane and Roynette}~\cite{GKR}. Namely, when we consider the
process at integer times, we obtain a random walk on the affine group,
which must be transient since that group is non--unimodular.

\smallskip

Also transience of $(Y_t)$ can be shown by 
constructing non-constant positive harmonic functions. 
More details are deferred to forthcoming work \cite{BSSW}, where among other we 
shall describe \emph{all} harmonic functions associated with 
$\bigl(\mathcal E^\T_{\alpha\beta}, \mathcal W^1_{\alpha,\beta}(\T)\bigr)$.
\end{proof}

\begin{rmk}Theorems \ref{thm-sg} and \ref{thm-saHT}, which describe some 
basic properties of the $(\alpha,\beta)$-heat semigroup and Laplacian on 
$\HT$ have obvious versions that apply to the heat semigroups and Laplacians 
on  $\T$ and $\Hb$ (respectively) that appear in the above result on
projections.
All these results illustrate the more general theory developed below in
the setting of what we call strip complexes. In fact, the introduction of 
the notion of strip complex is motivated in part by the justification of
the projections described above and the need to treat all these objects 
and their properties in a unified way.
\end{rmk}

\section{Strip complexes}\label{strip}

\subsection*{A. The basic structure of strip complexes}
Let $V,E$ be countable sets equipped with a map
$$
E \rightarrow V \times V\,,\; e \mapsto (e^-,e^+)\,.
$$
This defines an oriented graph $\Gamma$ with vertex set $V$ and edge set
$E$. We will assume throughout that $e^-\neq e^+$.
Hence multiple edges are allowed, but there
are no loops. The ``no loops'' convention will  simplify our considerations.
 Moreover, this is no real lack of generality for our purpose: loops can be 
handled by adding a virtual vertex in the middle of any existing loop.

The vertices $e^-,e^+$ are the extremities of the edge $e$.
We set $V_e=\{e^-,e^+\}\subset V$ 
and $E_v=\{e:v\in V_e\}$.
We let $\Gamma^1$ be the associated
$1$-dimensional complex. In $\Gamma^1$, the edge $e$ is realized by
a subset $I_e$ of $\Gamma^1$, homeomorphic to the closed interval $[0\,,\,1]$.
We will also use the notation $I_e=[e^-\,,\,e^+]$ and $I^o_e=(e^-\,,\,e^+)$
for the closed and open intervals corresponding to edge $e$, respectively.
Similarly, we write $\Gamma^o=\Gamma^1\setminus V$.
We assume throughout that $\Gamma^1$ is connected
and that each vertex has only finitely many neighbours, that is,
$E_v$ is a finite set.
For reasons that will become clear later, we refer to $\deg(v)= |E_v|$
as the \emph{bifurcation number} at $v$.

Although the edges are oriented, this orientation will
not play an important role for us. In particular, the notion of neighbours
introduced above does not take the orientation into account.
Observe also that we can view
$\Gamma^1$ as the union of all the edges $I_e$, $e\in E$, with the appropriate
identification  at the vertices where several edges meet.

Given a topological space $\Ms$ (we will be mostly  interested here
in the case where $\Ms$ is $\{o\}$, a line, a circle, or more generally,
a Riemannian manifold), the \emph{strip complex}
(more precisely, the $\Ms$-strip complex) associated to $\Gamma$ and
$\Ms$ is simply the direct product
$$
\GM =\Gamma^1 \times \Ms.
$$
This is a topological space with a simple ``coordinate system'' $\GM\ni
\xi=(\gamma,m)$. However, this viewpoint  is not entirely well suited
to capture the additional structure that these spaces
have in the cases of interest to us.

Instead, it will be essential to view $\GM$ as the union
of the \emph{strips}
$$
\bigcup_{e\in E} S_e\,,\quad\text{where}\;\; S_e= I_e\times \Ms.
$$
This is not a disjoint union,
as the strips $S_e=I_e\times \Ms$, $e\in E_v$, $v\in V$,
all meet along $\Ms_v=\{v\}\times \Ms$.
We call $\Ms_v$ the \emph{bifurcation manifold} at $v$.
This is simply the copy of $\Ms$ passing through $v$ in $\GM$.

(In Section \ref{geometry}, $\Ms = \R$, and the strips were labeld by the 
vertices of the tree, because there is a one-to-one correspondence between 
vertices $v$ and edges $[v^-,v]$.)

We let 
$$
S^o_e=(e^-,e^+)\times \Ms
$$
be the interior of the strip $S_e$ and set 
$$
\GM^o = \bigcup_{e \in E}S^o_e\,,
$$ 
the union of all open strips in $\GM$ (this is an open dense set in $\GM$).
For any function $f$ defined on $\GM^o$, we let 
$$
f_e= f|_{S^o_e}
$$
be the restriction of $f$ to the open strip  $S^o_e\,$. This notation
plays an important role and will be used throughout. In addition, we make the following 
natural convention. Whenever $f_e$ admits a continuous extension to the closed strip $S_e$,
we (abusively) use the same notation $f_e$ to denote this continuous extension.  Note that
if $f_e$ and $f_{e'}$ are defined on $S_e$ and $S_{e'}$  with $M_v=S_e\cap S_{e'}$ then
it may well be that $f_e$ and $f_{e'}$ take different values along $M_v$.

We also set
$$
X_v=\Ms_v \cup\biggl(\bigcup_{e\in E_v} S^o_e\biggr).
$$
The set $X_v$ is called the \emph{star of strips} at $v$. It is an open set 
in $\GM$.

\begin{rmk} Note that the definition of a strip complex given above 
is of a global nature and corresponds to what could be called ``untwisted'' 
strip complexes in the context of the following more general definition
which yields the same  local structure. In this more general definition,
the graph $\Gamma$ is decorated at each vertex by a collection 
$\{g^e_v:e\in E_v\}$ of homeomorphisms  $g^e_v: \Ms\rightarrow \Ms$
(when $\Ms$ is equipped with a Riemannian structure, these maps 
are required to be isometries).
Then, the boundaries $\Ms_v^e$ of different strips  $S^o_e$, $e\in E_v$,  
meeting at a vertex $v$, are  identified with a unique copy $\Ms_v$ of 
$\Ms$ through the homeomorphisms $g_v^e$. 
For instance, if $\Ms=(0\,,\,1)$, and the 
graph $\Gamma$ has  two vertices $a,a'$  and two edges $e,e'$
joining $a$ and $a'$, the strip complex $\GM=\Gamma^1\times \Ms$  
is a cylinder with two marked lines corresponding to $a,a'$.
However, we could identify  the two intervals $(0\,,\,1)$ 
at $a$ through the identity map  and at $a'$ through the flip  
$x\mapsto 1-x$. In this case, we get a ``twisted strip complex'' 
which is a Moebius 
band with two marked lines. Note that this ``twisted strip complex'' is not 
globally the direct product of $\Gamma^1$ and $\Ms$ although, locally, it 
has the same structure.  We will not discuss twisted strip complexes 
in this paper. But we note that all of our results (properly interpreted) 
will  hold as well  for such more general structures. In particular, 
our local smoothness results will apply to these twisted structures in an
obvious way.
\end{rmk}

\begin{rmk} The treebolic space (see Figure 1) gives a good illustration 
of a strip complex structure, but it may be useful for the reader to think 
of the case when $\Ms$ is the circle $\mathbb T=\R/(2\pi\Z)$
and $\Gamma$ is some finite graph. Although one can easily draw sketches 
of such examples, in most cases, 
these circle strip complexes cannot be embedded (without crossings)
in three-space.
\end{rmk}

\bigskip

\subsection*{B. Smooth functions on strip complexes}

Fix a graph $\Gamma$ as defined above. Let  $\Ms$ be an $n$-dimensional 
manifold and consider the associated strip complex $\GM$.
Let $\mathcal C_c(\GM)$,
$\mathcal C_0(\GM)$ and $\mathcal C_b(\GM)$
be the spaces of continuous functions on $\GM$ that
are, respectively,  compactly supported, vanishing at infinity, bounded.

Without further comments, we will assume that $\Ms$ is equipped with a Radon
measure which, in any coordinate chart on $\Ms$, admits a
smooth positive density with respect to the Riemannian measure. The strip complex
$\GM$ is then equipped with the product measure of one-dimensional 
Lebesgue measure on $\Gamma^1$ and the given Radon measure on $M$. 
Later we will make a more precise choice of such a measure. 
For the time being, this measure is used only
for the definition of negligeable sets (sets of measure zero) and the
particular choice made is irrelevant.

\begin{dfn}\label{def-chart} 
A \emph{relatively compact coordinate chart} in $\GM$ is an open,
relatively compact set of the form $I\times U\subset \GM$ where
$I\subset (e^-\,,\,e^+)\subset \Gamma^1$ for some $e\in E$ is an open interval
and  $(U;x_1,\cdots,x_n)$ is a relatively compact coordinate chart in $\Ms$.
The associated \emph{local coordinate system} on the open subset $I\times U$
is denoted by $\xi=(s,x_1,\dots,x_n)$, $s\in I$, $(x_1,\dots,x_n)\in U$.
For any $(n+1)$-tuple $\kappa=(\kappa_0,\kappa_1,\dots, \kappa_n)$
of integers and any smooth enough function $f$ defined over $I\times U$,
we set
$$
\partial^\kappa_\xi f(\xi)=\partial_s^{\kappa_0}
\partial^{\kappa_1}_{x_1}\dots
\partial^{\kappa_n}_{x_n}f(s,x_1,\dots,x_n).
$$
If necessary, we can also consider $\partial^\kappa_\xi f$
to be defined in the sense of distributions in $I\times U$.
\end{dfn}

\begin{rmk}
The above definition never involves the bifurcation manifolds,
except possibly at the boundary of $I\times U$. Hence, smoothness of a function
in a relatively compact chart $I\times U$
as defined above is a classical notion.
\end{rmk}

\begin{dfn}\label{def-smooth}
(a) The space of \emph{strip-wise smooth functions} 
on $\GM^o$, denoted  $\mathcal S^\infty(\GM^o)$, is
the set of those locally bounded functions $f$  on $\GM^o$ such that,
for any open edge $I^o_e=(e^-\,,\,e^+)$, $e\in E$, and any precompact
coordinate chart $(U;x_1,\dots,x_n)$ in $\Ms$, the function 
$f|_{I^o_e\times U}$ is a bounded continuous function with
bounded continuous derivatives of all orders with respect
to the coordinates $(s,x_1,\dots,x_n)$ in $I^o_e\times U$.
The vector space $\mathcal S^\infty(\GM^o)$
is equipped with the family of seminorms
\begin{equation}\label{seminorms}
\begin{aligned}
N^k_{K,I\times U}(f)&=\sup\{|f(\xi)| : \xi \in K \cap \GM^o\} \\
&\quad +
\sup\Bigl\{\, |\partial^\kappa_\xi f(\xi)|: \xi\in I\times U,
\kappa=(\kappa_0,\kappa_1,\dots,\kappa_n), 
\textstyle{\sum_0^n}\kappa_i\le k\Bigr\}\,,
\end{aligned}
\end{equation}
where $k$ is an integer, $K$ a compact subset of $\GM$ and $I\times U$
a relatively compact coordinate chart in $\GM$.

Abusing notation, we will also consider any function $f$ in 
$\mathcal S^\infty(\GM^o)$ as a function on $\GM$ that is defined almost 
everywhere (a representative of a class of functions under the 
usual equivalence of coinciding almost everywhere).

\smallskip

(b) The space of \emph{continuous strip-wise smooth functions} on $\GM$, denoted
$\mathcal C^\infty(\GM)$ is defined as 
$$
\mathcal C(\GM)\cap \mathcal S^\infty(\GM^o) 
= \{ f \in  C(\GM) : f|_{\GM^o} \in S^\infty(\GM^o) \}\,.
$$ 
We also let
$$
\mathcal C^\infty_c(\GM)=\mathcal C^\infty(\GM) \cap \mathcal C_c(\GM).
$$
The vector space $\mathcal C^\infty(\GM)$ is equipped with the same family of
seminorms $N^k_{K,U}$ as $\mathcal S^\infty(\GM^o)$.
\end{dfn}

\begin{rmks} 
(i) A function $f\in \mathcal S^\infty(\GM^o)$ is not necessarily
continuous across bifurcation manifolds  (it need not even be defined on the 
latter). However, the functions $f_e$ are bounded
continuous with bounded continuous derivatives on $I_e^o\times U$
for any relatively compact set $U\subset \Ms$. This implies that each $f_e$
can be extended as a smooth continuous function to the closed strip $S_e$.
According to our earlier convention, we still denote this extension by $f_e\,$. 
In particular, for any vertex $v$,
a function $f\in \mathcal S^\infty(\GM)$, yields $\deg(v)$ smooth functions
$$
\Ms \ni x\mapsto f_e(v,x) \in \mathcal C^\infty(\Ms).
$$

\smallskip

\noindent
(ii) Note that a function in $\mathcal S^\infty(\GM^o)$ may not have a 
continuous  extension to $\GM$ but is always (essentially) bounded 
on compact sets.

\smallskip

\noindent
(iii) The space $\mathcal C^\infty(\GM)$ is a complete seminormed space.
In view of (i), a function $f\in \mathcal C^\infty(\GM)$
is a continuous function on $\GM$ such that its restriction $f_e$
to any closed strip $S_e$ is a smooth function in the usual sense on
the manifold $S_e$. 

Since $f$ is continuous it follows that the
partial derivatives $\partial_x^\kappa f$, $\kappa=(\kappa_1,\dots,\kappa_n)$
in the direction of $\Ms$  have to be continuous
across bifurcation manifolds. That is, for any fixed coordinate chart
$(U;x)$ in $\Ms$, with $x=(x_1,\dots,x_n)$,
$$
\partial_x^\kappa f_{e_1}(v,x)=
\partial_x^\kappa f_{e_2}(v,x)\,, \quad\text{if}\;\;e_1,e_2\in E_v\,.
$$
Note, however,
that the partial derivatives $\partial_s^k\partial^\kappa_xf_e$ with $k\ge 1$
and computed in different strips meeting along a bifurcation manifold $\Ms_v$
do not have to match along $\Ms_v$.
\end{rmks}

\begin{rmk} We will sometimes consider functions $f$ of space and 
time variables,
such as for example $(0,T)\times \GM\ni(t,\xi)\mapsto f(t,\xi)$. 
Since $(0,T)\times \GM$ is also a
strip complex, with $\Ms$ replaced by $\Ms\times(0,T)$,
(\ref{def-smooth}.b) also defines
$\mathcal C^\infty\bigl((0,T)\times \GM\bigr)$.
\end{rmk}

The following subspace of $\mathcal C^\infty_{c}(\GM)$ will be useful
for our purpose. It is the subspace of those functions in
$\mathcal C^\infty_{c}(\GM)$ which are locally constant along $\Gamma^1$
near each bifurcation manifold $\Ms_v\,$.

\begin{dfn}
Let $\mathcal C^\infty_{c,c}(\GM)$ be the subspace of
$\mathcal C^\infty_{c}(\GM)$ of those functions whose partial derivative
$\partial_s f_e$ in any strip $S_e=I_e\times \Ms$, $s\in I_e\,$,
has compact support in $S^o_e\,$.
\end{dfn}

\begin{lem}\label{lem-dens0}
The space $\mathcal C^\infty_{c,c}(\GM)$ is dense
in $\mathcal C_{0}(\GM)$ for the uniform norm. 
\end{lem}
\begin{proof}Since $\mathcal C_c(\GM)$ is dense in $\mathcal C_0(\GM)$
for the uniform norm, it suffices to show that for any $f\in
\mathcal C_c(\GM)$ and $\epsilon>0$ there is
$f_\epsilon\in \mathcal C^\infty_{c,c}(\GM)$ such that
$\|f-f_\epsilon\|_\infty\le \epsilon.$

Let $K$ be the support of $f$ and $\{U_n,n\le N\}$ be a finite covering of $K$
by open precompact subsets which are so small that for each $n$,
$\sup\{|f(\xi)-f(\zeta)|: \xi,\zeta\in U_n\}<\epsilon$
(uniform continuity of $f$) and $U_n$ is either of the form
$J_n\times V_n$ where $V_n$ is a small coordinate chart in $\Ms$
and $J_n$ is relatively compact in $(e^-,e^+)$ for some $e$, or
$U_n=\bigcup_{e\in E_v} J^e_n\times V_n$
where $V_n$ is a small coordinate chart in $\Ms$ and each $J_n^e$ 
is a semi-open interval in $I_e$ with closed extremity at $v$. 
By standard arguments adapted to the present situation,
we can construct a family of functions $\omega_n\in \mathcal C^\infty_{c,c}(\GM)$
such that $\omega_n$ is supported in $U_n$ and $\sum_{n\le N}\omega_n=1$
on $K$. For each $n\le N$, pick $\xi_n\in U_n$ and set
$$
f_\epsilon=\sum_{n\le N}f(\xi_n)\,\omega_n\,.
$$
By construction, $f_\epsilon\in \mathcal C^\infty_{c,c}(\GM)$ and,
for any $\xi\in K$,
$$
|f-f_{\epsilon}|(\xi)\le \sum_{n\le N}
|f(\xi)-f(\xi_n)|\omega_n(\xi)\le \epsilon.
$$
This provides the desired approximation.
\end{proof}

The next definition introduces smoothness (of various orders)
in an open subset $\Omega$ of $\GM$.

\begin{dfn}\label{smoothOmega}
Let $\Omega$ be an open subset of $\GM$. 

\smallskip

(a) A function $f$ is in $\mathcal C^k(\Omega)$, where $k \ge 1$, 
if it is continuous in $\Omega$, and for any relatively compact 
coordinate chart $I\times U$
with $\overline{I\times U}\subset \Omega$,
$f$ has continuous partial derivatives of order up to $k$ in $I\times U$.
This space is equipped with the family of seminorms
$N^k_{K,I\times U}$ defined as in \eqref{seminorms},
where $K$ runs over  compact subsets of $\Omega$ and $I\times U$
over all  relatively compact coordinate charts with
$\overline{I\times U}\subset \Omega$.

\smallskip

(b) A function $f$ is in
$\mathcal C^\infty(\Omega)$ if it is continuous in $\Omega$
and for any relatively compact coordinate chart $I\times U$
with $\overline{I\times U}\subset \Omega$,
$f$ has continuous partial derivatives of all orders in $I\times U$.
This space is equipped with the family of seminorms
$N^k_{K,I\times U}$ defined as in \eqref{seminorms},
where $k$ runs over the positive integers,
$K$ runs over  compact subsets of $\Omega$ and $I\times U$
over all  relatively compact coordinate charts with
$\overline{I\times U}\subset \Omega$.
\end{dfn}

The spaces $\mathcal C^k(\Omega)$ and $\mathcal C^{\infty}(\Omega)$ are
complete seminormed spaces.

\bigskip

\subsection*{C. Diffeomorphisms}
Let $\GMa$, $\GMb$ be two strip complexes. Since these spaces are equipped 
with a natural topology, the notion of homeomorphism is well defined. 
Observe that bifurcation manifolds $\Ms_v$ with bifurcation number 
$\deg(v)=2$ may be ignored by a homeomorphism. Otherwise, by definition,
a homeomorphism must send strips to strips and
send  any bifurcation manifold with bifurcation number $\deg(v)>2$
to a bifurcation manifold with the same bifurcation number.

\begin{dfn}
Let $\GMa$, $\GMb$ be two strip complexes. 
A homeomorphism $\jm:\GMa\rightarrow \GMb$ is called a \emph{diffeomorphism} 
if $\jm$ and $\jm^{-1}$ send any bifurcation manifold
to a bifurcation manifold and, for any pair of closed strips
$S_1\subset \GMa$, $S_2\subset \GMb$ such that $\jm(S_1)=S_2\,$, the 
restriction $\jm|_{S_1}:S_1\mapsto S_2$ is a diffeomorphism.

A \emph{local diffeomorphism between} open sets $\Omega_1,\Omega_2$
is a map $\jm:\Omega_1\rightarrow \Omega_2$ which is a homeomorphism,
sends any trace of a bifurcation manifold
to a trace of a bifurcation manifold and is a diffeomorphism between traces
of closed strips.
\end{dfn}

\begin{rmks} (1) Diffeomorphisms must respect the bifurcation structure,
even for bifurcation manifolds with bifurcation number $\deg(v)=2$.

\smallskip

(2) If $\jm:\GMa\rightarrow \GMb$ is a diffeomorphism
then for any $f\in\mathcal C^\infty(\GMb)$, resp. $\mathcal S^\infty(\GM^o)$,
the function $f\circ \jm$ is in $\mathcal C^\infty(\GMa)$, resp. 
$\mathcal S^\infty(\GM^o)$.
If $\jm:\Omega_1\rightarrow \Omega_2$ is a local  diffeomorphism then for any 
$f\in\mathcal C^\infty(\Omega_2)$, the function $f\circ \jm$ is in 
$\mathcal C^\infty(\Omega_1)$. The same holds for functions that are smooth
up to order $k$.
\end{rmks}

\bigskip

\subsection*{D. Geometric structures on strip complexes}

We now introduce a rather specific class of geometric structures
on the strip complex $\GM$. This is done in two stages. The special features
of these structures will play a central role in our analysis.

In the first stage, we introduce a product geometric structure
on $\GM$ associated with given geometric structures on $\Gamma$ and $M$
as follows.

First, we assume that the edge map contains an additional information,
namely, the \emph{length} of the edge $e$. More precisely, we have a map
$$
E \to V\times V\times(0,\infty)\,,\quad e \mapsto (e^-,e^+,l_e).
$$
Thus, with this additional information, the edge $I_e=[e^-,e^+]$ is isometric
to the real interval $[0\,,\,l_e]$.
We can view $\Gamma^1=(\Gamma^1,l)$ as a metric space in the obvious way.
We will always use the letter $s$ to refer to an
arc length parameter on $\Gamma^1$ or connected pieces of $\Gamma^1$.
From now on, we always assume that $\Gamma$ comes equipped with a
specific edge length map $l$.

Second, we assume that $(\Ms,g)$ is a Riemannian manifold with gradient
$\nabla_{\Ms}$. Given these two geometric inputs
(length of edges, Riemanian metric on $\Ms$), we immediately obtain a natural
metric on $\GM$ by equipping  each strip $S_e=I_e\times \Ms$ with the
Riemannian metric $(ds)^2+g_x(\cdot,\cdot)$, where $(s,x)\in I_e\times \Ms$.

Here and elsewhere, the subscript $x$ in $g_x$ indicates that
$g$ is considered with respect to the $x$-variable of $(s,x)$.

The second stage of our construction depends on the choice of
a function $\phi$, positive and strip-wise smooth on $\Gamma^o = 
\Gamma^o$,  that is, $\phi\in \mathcal S^\infty(\Gamma^o)$.
On each  strip $S_e=I_e\times M$, we consider the smooth Riemannian structure
\begin{equation}\label{Riem-structure}
\phi_e(s)\cdot\bigl[(ds)^2+g_x(\cdot,\cdot)\bigr]
\end{equation}
obtained from the product structure by multiplication by $\phi$.
The associated Riemannian measure is 
$\phi_e(s)^{(1+n)/2}\,ds\,dx\,$, 
where $dx$ is the volume element of $\Ms$ (resp. area or length element,
according to the dimension of $\Ms$). 
This induces our reference measure on $\GM$ that reflects the 
underlying geometry, given by
\begin{equation}\label{basic-measure}
\sum_{e \in E} \phi_e(s)^{(1+n)/2}\,\uno_{S_e^o}\,\,ds\,dx\,.
\end{equation}
Note that $\GM \setminus \GM^o$, the union of all the bifurcation manifolds,
is a negligeable set. (Below we shall consider
a larger class of measures, associated forms and processes.)
We are led to the following.

\begin{dfn}\label{def-gradsquare} 
Let $f,h$ be  functions in $\mathcal S^\infty(\GM^o)$.
The gradient $\nabla f$ and its length square are given
at $(s,x)\in S^o_e$ by
$$
\nabla f(s,x)=\frac{1}{\phi_e(s)}
\bigl(\partial_sf_e(s,x),\nabla_{\Ms} f_e(s,x)\bigr)
$$
and
$$
|\nabla f(s,x)|^2= \frac{1}{\phi_e(s)}
\Bigl(|\partial_s f_e(s,x)|^2
  + g_x\bigl(\nabla_{\Ms} f_e(s,x),\nabla_{\Ms} f_e(s,x)\bigr)\Bigr),
$$
that is, 
$
|\nabla f|^2= \sum_{e\in E} \frac{1}{\phi_e(s)}|\nabla f_e|^2 \,
{\mathbf 1}_{S^o_e}\,.
$
Correspondingly, the inner product of the gradients at $(s,x) \in \GM^o$ is
$$
(\nabla f,\nabla h)(s,x) =
\sum_{e\in E} \frac{1}{\phi_e(s)}
\Bigl(\partial_s f_e(s,x)\,\partial_s h_e(s,x)
  + g_x\bigl(\nabla_{\Ms} f_e(s,x),\nabla_{\Ms} h_e(s,x)\bigr)\Bigr)\,. 
$$
\end{dfn}
Note that these definitions involve the edge
length function $l$, the metric $g$ on $\Ms$ and the function
$\phi\,$, but these are omitted in our notation. 

Now, if we have a continuous path in $\GM$
which is rectifiable (i.e., is rectifiable in each strip), then we can
compute its length by adding the lengths of the parts of the path within 
each strip.

\begin{dfn} For any two points $\xi,\zeta\in \GM$, let
$\rho(\xi,\zeta)$ be the infimum of the lengths of continuous rectifiable
paths in $\GM$ joining $\xi$ to $\zeta$.
\end{dfn}

One easily checks that this defines a distance function on $\GM$ which 
defines the original topology of this space. We set
$$
B(\xi,r)=\{\zeta\in \GM: \rho(\xi,\zeta)<r\},
$$
the open ball with radius $r$ around $\xi$.

The (easy) proof of the following lemma is left to the reader.

\begin{lem}\label{lem-GMcomp} Assume that $(M,g)$ is a complete Riemannian
manifold. Then the metric space $(\GM,\rho)$ is complete if and only if
the metric space $(\Gamma^1,\rho)$ is complete. This is the case if and only if,
for any infinite family $F\subset E$ of edges such that
$\bigcup_{e\in F} I_e$ is connected in $\Gamma^1$, we have
\begin{equation}\label{GMcomp}
\sum_{e\in F}  \int_{I_e}\sqrt{\phi_e(s)}\,ds=\infty\,.
\end{equation}
\end{lem}

\begin{dfn}\label{def-iso}
Given two strip complexes $\GMa,\GMb$, each equipped
with respective geometric structures $(\phi_1,\rho_1)$ and $(\phi_2,\rho_2)$ 
as above, we say that a diffeomorphism $\jm:\GMa\rightarrow \GMb$ is an 
\emph{isometry} if it satisfies
$$
\rho_2\bigl(\jm(\xi),\jm(\zeta)\bigr)=\rho_1(\xi,\zeta)\quad
\text{for all}\;\;\xi\,,\;\zeta\in \GMa\,.
$$
A \emph{local isometry} between two open sets
$\Omega_1\,,\;\Omega_2$ is defined analogously.
\end{dfn}

\begin{rmk}
If $\,\jm$ is an isometry then for any
$f\in \mathcal C^\infty(\GMb)$ and any $\xi$ in the interior of a strip in $\GMa$, we have
$$
\bigl(\nabla_1 f\circ \jm\,,\nabla_1 h\circ \jm\bigr)_1 (\xi)
=\bigl(\nabla_2 f,\nabla_2h\bigr)_2\bigl(\jm(\xi)\bigr).
$$
Indeed the differential map $d\jm|_\xi$ is an isometry between the tangent spaces
at $\xi$ and $\jm(\xi)$, when $\xi \in \GMa^{\!\! o}$.
\end{rmk}

\bigskip

\subsection*{E. Dirichlet forms on $\GM$}

We now equip $\GM$ with a measure $d\mu$ which will serve as
our basic underlying measure to define ${\mathcal L}^p$ spaces on $\GM$, in particular,
${\mathcal L}^2(\GM,\mu)$. This measure  $\mu$ is described by its density
$\psi\in \mathcal S^\infty(\Gamma^o)$  with respect to the basic
measure of \eqref{basic-measure}.

\begin{dfn} \label{dfn-mu}
(a) Given the positive function $\psi\in \mathcal S^\infty(\Gamma^o)$,
let $\mu = \mu_{\psi}$ be the positive Radon measure on $\GM$  such that, 
for any $f\in  \mathcal C_c(\GM)$,
\begin{eqnarray*}
\int_{\GM}f\, d\mu &=& \int_{\GM} f(s,x)\, \psi(s)\,\phi(s)^{(1+n)/2}\,
ds\, dx\\
&=&\sum_{e\in E}\int_{S^o_e} f_e(s,x)\,\psi_e(s)\,\phi_e(s)^{(1+n)/2}\,ds\,dx\,,\\
\end{eqnarray*}
where $ds$ is Lebesgue measure on $(\Gamma^1,l)$ and $dx$ is the
Riemannian measure on $(M,g)$.
\smallskip

(b) For each $\xi\in \GM$ and $r>0$ set
$$
V(\xi,r)=\mu\bigl(B(\xi,r)\bigr)
=\sum_{e\in E} \int_{B(\xi,r)\cap S^o_e}\psi_e(s) \, \phi_e(s)^{(1+n)/2} \,ds\, dx\,.
$$
\end{dfn}

Above, Lebesgue measure on $(\Gamma^1,l)$ is of course the measure which
restricted to each edge is Lebesgue measure assigning length $l_e$
to $I_e\,$, while the vertex set has measure $0$.

To construct Dirichlet forms on $\GM$, we need to recall a version of the 
classical trace theorem for Sobolev spaces. 
For any strip $S_e$,
consider the set $\mathcal W^1(S_e,\mu)=\mathcal W^1(S^o_e,\mu)$ 
of all functions $f$ in 
${\mathcal L}^2(S^o_e,\mu)$ whose  distributional first derivatives 
in $S^o_e$ can be represented by functions in ${\mathcal L}^2(S^o_e,\mu)$.
Note that, by definition, $\mathcal W^1(S_e,\mu)=\mathcal W^1(S^o_e,\mu)$.
However choosing $S_e$ or $S^o_e$ makes a difference when considering the 
local versions of this space since compact subsets of $S^o_E$ and $S_e$ are differeent. 
We let $\mathcal W^1_{\mbox{\tiny loc}}(S_e,\mu)$ be the space of  
of all functions $f$ in 
${\mathcal L}^2_{\mbox{\tiny loc}}(S_e,\mu)$ whose  distributional first derivatives 
in $S^o_e$ can be represented by functions in ${\mathcal L}^2_{\mbox{\tiny loc}}(S_e,\mu)$.

For any $f$ in $\mathcal W^1_{\mbox{\tiny loc}}(S_e,\mu)$, using the global coordinates
$(s,x)$ on $S_e=I_e\times M$, we have that the derivative $\partial_s f$, the 
$\Ms$  gradient $\nabla_M f$ and the global gradient $\nabla f$ are 
well defined locally square integrable functions on $S_e$.
In particular, for such functions, the length square and inner product
of the gradient(s) are well defined as  locally integrable functions
in the sense of Definition \ref{def-gradsquare}.

By the classical trace theorem, those functions admit a trace on 
each of the copies $\Ms_{e^-}$ and $\Ms_{e^+}$ of $\Ms$ bounding the strip 
$S_e\,$. More precisely, there exist two continuous  linear operators
\begin{equation}\label{traces}
\mbox{Tr}^{S_e}_{M_{e^\pm}}: \mathcal W^1_{\mbox{\tiny loc}}(S_e,\mu) \rightarrow
{\mathcal L}^2_{\mbox{\tiny loc}}(M_{e^\pm},dx)
\end{equation}
which extend the natural restriction operators
defined from  $\mathcal C^\infty(S_e)$ to $\mathcal C^\infty(M_{e^\pm})$.

\begin{dfn} Given $\Gamma$, $(M,g)$ and $\phi,\psi\in \mathcal S^\infty(\Gamma^o)$, 
as above, let $\mathcal W^1(\GM,\mu)$ be the space of those functions $f$ 
in ${\mathcal L}^2(\GM,\mu)$ whose restrictions $f_e$, $e\in E$, are all in
$\mathcal W^1_{\mbox{\tiny loc}}(S_e)$ and satisfy:
\begin{itemize}
\item
$\int_{\GM} |\nabla f|^2d\mu<\infty$
\item For any vertex $v$ and any two edges $e,e'\in E_v$,
$\mbox{Tr}^{S_e}_{M_v}f=\mbox{Tr}^{S_{e'}}_{M_v}f$.
\end{itemize}
\end{dfn}

\begin{dfn}For $f,h\in \mathcal W^1(\GM,\mu)$, set
$$
\mathcal E(f,h)= \int_{\GM}(\nabla f, \nabla h)\, d\mu\,.
$$
Let $\mathcal W^1_0(\GM,\mu)$ be the closure of
$\mathcal C^\infty_c(\GM)$ in $\mathcal W^1(\GM,\mu)$.
\end{dfn}

\begin{exa} \label{main-ex-HT}
Let $\Gamma=\T = \T_\pp$ be a $\pp$-regular tree equipped with 
an origin $o$, a reference end $\om$, and the associated horocycle 
function $\hor$.
Edges are oriented away from $\om$ so that $\hor(e^+)=\hor(e^-)+1$. See 
Figure \ref{fig2}. Turn $\T$ into a metric tree by given length 
$q^{k-1}(q-1)$ to all edges $e$ with $\hor(e^-)=k-1$. Define 
$\phi\in \mathcal S^\infty(\T^1 \setminus V)$ by
$\phi_e(s)=  s^{-2}$ on $I_e\cong[q^{k-1},q^k]$ if $\hor(e^-)=k-1$. Setting 
$\Ms=\R$, the corresponding structure on $\GM$ is isometric to that of the 
treebolic space $\HT(\pp,\qq)$. Next, for any fixed $\alpha\in \R$, define
$\psi\in \mathcal S^\infty(\T^1 \setminus V)$
by $\psi_e(s)=  s^{\alpha}$ on $I_e\cong[q^{k-1},q^k]$ if $\hor(e^-)=k-1$. 
Then the corresponding measure $\mu$ on $\GM$ is the measure 
$\mu_{\alpha,\beta}$ on $\HT(\pp,\qq)$ (modulo the isometry mentioned 
earlier between these two spaces)
and the associated Dirichlet form $\bigl(\mathcal E, \mathcal W^1(\GM,\mu)\bigr)$
is the form $\bigl(\mathcal E_{\alpha,\beta}\,,\mathcal W^1_{\alpha,\beta}(\HT)\bigr)$ 
from Definition \ref{def-Eab}. 
\end{exa}

\begin{thm} \label{thm-strilo}
The quadratic form $\bigl(\mathcal E, \mathcal W^1(\GM,\mu)\bigr)$ is a
strictly local   Dirichlet form and the quadratic form
$\bigl(\mathcal E, \mathcal W^1_0(\GM,\mu)\bigr)$ is a
strictly local regular Dirichlet form.
\end{thm}

\begin{proof} The Markov character and strict locality of these forms are clear
from the definitions. See \cite{FOT}.
The fact that the first form is closed follows from the
fact that the corresponding forms on all strips are closed and from
the continuity of the trace operators. The fact that the second form is closed
and regular is obvious from the definition and the fact that
$\mathcal C^\infty_c(\GM)$ is dense in $\mathcal C_0(\GM)$
for the uniform norm (see Lemma \ref{lem-dens0}).
\end{proof}

\begin{thm} \label{th-reg}
Assume that $(\GM,\rho)$ is a complete metric space
(see {\em Lemma \ref{lem-GMcomp}} for a necessary and sufficient condition).
Then the forms $\bigl(\mathcal E, \mathcal W^1(\GM,\mu)\bigr)$ and
$\bigl((\mathcal E, \mathcal W^1_0(\GM,\mu)\bigr)$ coincide. In particular,
$\bigl((\mathcal E, \mathcal W^1(\GM,\mu)\bigr)$ is a strictly local regular Dirichlet 
form.
\end{thm}

\begin{proof} To prove this, we simply need to show that
$\mathcal C^\infty_c (\GM)$ is dense in $\mathcal W^1(\GM,\mu)$. First, we show that
any $f\in \mathcal W^1(\GM,\mu)$ can be approximated in $\mathcal W^1(\GM,\mu)$ by functions
with compact support. Consider the distance function $\rho$ on $\GM$.
Observe that, for any set $U$,
the function $\xi \mapsto \rho(\xi,U)$ is a contraction in each strip $S_e$.
Therefore this function is in $\mathcal W^1_{\mbox{\tiny loc }}(\GM)$
with $|\nabla \rho_U|\le 1$. If follows that the functions
$$
\theta_n = \max\{1-\rho(\cdot,B(o,n))/n,0\}\,,
$$
where $o$ is a fixed point in $\GM$, are in $\mathcal W^1(\GM,\mu)$ and satisfy
$|\nabla \theta_n|\le 1/n$. The function $\theta_n$ is supported in $B(o,2n)$,
which is precompact since
$(\GM,\rho)$ is a complete locally compact metric space.
This yields that the compactly supported functions
$\theta_n f$ converge to $f$ in $\mathcal W^1(\GM,\mu)$.

Next, we show that any compactly supported function $f$ in $\mathcal W^1(\GM,\mu)$
can be approximated in $\mathcal W^1(\GM,\mu)$ by compactly supported functions in
$\mathcal C^\infty_c(\GM)$. By compactness of the support of $f$,
we can find a finite collection of functions $\omega_i$ in
$\mathcal C^\infty_c(\GM)$ such that $\sum\omega_i=1$ on the support of $f$
and each $\omega_i$ either has its compact support in
 an open strip $(e^-,e^+)\times M$ or $\omega_i$ has its compact support in
a star of strips $X_v$ at vertex $v$. At this point, it suffices to approximate
each $f\, \omega_i$ by functions in $\mathcal C_c^\infty(\GM)$.
If $\omega_i$ has compact support within one open strip, this follows from
a classical procedure. 

The interesting case is when $\omega_i$ is compactly supported in a star 
$X_v\,$. In this case we can assume that the support of $\omega_i$ is so 
small that it is contained in an open  set of the form 
$\bigcup_{e\in E_v} U_e\,$, where
the $U_e$ meet on $\Ms_v$ along an open set $U_v=\{v \} \times U \subset \Ms_v$ 
and each $U_e$ is of the form
$J_e\times U$ where the $J_e\subset I_e$ are semi-open intervals of the same length
all containing $v$. 

Now pick one of the edges $\tilde e\in E_v\,$, and
let $\tilde f$ be the function  which,  on each $U_e$, equals 
$f\,\omega_i|_{U_{\tilde e}}$, and is zero outside of $\bigcup_{e\in E_v} U_e\,$.
That is, we copy the values of $f\,\omega_i$ from $U_{\tilde e}$ to all the 
other $U_e\,$, $e \in E_v\,$, via the obvious coordinate-wise correspondences 
between those sets, taking into account the identification between 
$U_{\tilde e}$ and the other sets $U_e$ along $U_v$.
On each strip $S^o_e$, the  function $f\,\omega_i- \tilde f$ is in 
$\mathcal W^1_0(S^o_e,\mu)$ because, by construction, the functions $f\,\omega_i$ and 
$\tilde f$ coincide on $U_v\subset  M_v$.  Hence we can approximate 
$f\,\omega_i-\tilde f$ in $\mathcal W^1(\GM,\mu)$
by functions $g_n$  whose restrictions to each $S^o_e$, $e\in E_v$, are smooth 
and compactly supported in the respective set $U_e\,$. Those $g_n$ are in
$\mathcal C^\infty_{c,c}(\GM)$. 

Next, the function  $f\,\omega_i|_{U_{\tilde e}}$
is in $\mathcal W^1(S_{\tilde e},\mu)$ with compact support in $U_{\tilde e}$. 
Recall that $U_{\tilde e}$ contains $U_v$ as part of its boundary.
By classical constructions, $f\,\omega_i|_{U_{\tilde e}}$
can be approximated by functions $h_n \in \mathcal C^\infty_c(U_{\tilde e})$.
We now use $h_n$ to define $\tilde h_n$ on $\bigcup_{e\in E_v} U_e$ by setting,
for each $e\in E_v$,  $\tilde h_n|_{U_{e}}=h_n$ in the same way as above
via the natural correspondence between $U_e$ and $U_{\tilde e}\,$.
Obviously, $\tilde h_n\in \mathcal C^\infty_c(\GM)$
and it  approximates $\tilde f$ in $\mathcal W^1(\GM,\mu)$. This implies that 
$g_n+\tilde h_n$, which is in $\mathcal C^\infty_c(\GM)$, approximates 
$f\omega_i$ in $\mathcal W^1(\GM,\mu)$.
\end{proof}

In fact, the smaller space $\mathcal C^\infty_{c,c}(\GM)$
is already dense in $\mathcal W^1_0(\GM,\mu)$, and thus in $\mathcal W^1(\GM,\mu)$ when
$(\GM,\rho)$ is complete. Recall that $\mathcal C^\infty_{c,c}(\GM)$
is the set of those functions in $\mathcal C^\infty_{c}(\GM)$
such that in any strip $S_e=I_e\times \Ms$, the partial derivative
$\partial_s f|_{s_e}$ has compact support contained in the open strip $S^o_e$
(as usual, $s$ is the variable in the interval $I_e$).

\begin{thm}\label{thm-dense}
 The subspace $\mathcal C^\infty_{c,c}(\GM)$ of $\mathcal W^1_0(\GM,\mu)$
is dense in $\mathcal W^1_0(\GM,\mu)$, and thus in $\mathcal W^1(\GM,\mu)$, when $(\GM,\rho)$
is complete.
\end{thm}
\begin{proof}
To see that this is the case, we return to the end of the argument
in the proof of Theorem \ref{th-reg}. We claim we can approximate
$f\,\omega_i|_{U_{\tilde e}}\in \mathcal W^1(S_{\tilde e},\mu)$
by functions $h_n \in \mathcal C^\infty_{c,c}(U_{\tilde e})$.
If that is the case, we use $h_n$ to define $\tilde h_n$ on 
$\bigcup_{e\in E_v}U_e$ by copying the values of $\tilde h_n$ from 
$U_{\tilde e}$ to $U_e$ for each $e\in E_v\,$. Obviously, 
$\tilde h_n\in \mathcal C^\infty_{c,c}(\GM)$ and it  approximates $\tilde f$ 
in $\mathcal W^1(\GM,\mu)$. Then, as before, $g_n+\tilde h_n$,
approximates $f\,\omega_i$ in $\mathcal W^1(\GM,\mu)$ as desired. The function
$g_n+\tilde h_n$ is in $\mathcal C^\infty_{c,c}(\GM)$ because $\tilde h_n$ 
is in that space by construction and $g_n$ has compact support in the union 
$\bigcup_{e\in E_v}S^o_e$
of the open strips surrounding $M_v$ and thus is also in
$\mathcal C^\infty_{c,c}(\GM)$.

Thus, the only thing left to prove is that a function $f\in \mathcal W^1(R)$,
$R=[e^-\,,\,e^+) \times U$, with compact support in $R$ can be approximated 
in $\mathcal W^1(R)$ 
by a sequence of functions $h_n\in \mathcal C^\infty(R)$ with compact support 
in $R$ and such that $\partial_s h_n$ has compact support in
$I_e^o\times U$. Note that, by definition, $R$ contains the bottom 
$\{e^-\}\times U$.

Since this is a local problem, we can  regard $U$ as a small open set in
$\R^n$ that contains the origin, and ignore completely the role of the 
functions $\phi,\psi$. Modifying notation in this sense, we use coordinates 
$(s,x)\in R = [0,l)\times U$ (instead of $s \in [e^-\,,\,e^+)$,
with $l=l_e$) and write $d\mu=ds\,dx$. For $n=1,2,\dots$, set
$$
f_n(s,x)=n\int_{s}^{s+1/n}f(\tau,x)\,d\tau \AND
\tilde f_n(s,x) =\begin{cases} f_n(s,x)\,,&\text{if}\; s\in (1/n,l)\\
                               f_n(1/n,x)\,,& \text{if}\; s\in [0,1/n].
                 \end{cases}
$$
We can assume that the support of $f$ in $R$ is small enough so that
$f_n$ and $\tilde f_n$ are still supported in $R$. It is plain that $f_n$ 
tends to $f$ in $\mathcal W^1(R)$, and we claim that the same is true for 
$\tilde f_n\,$. It is clear that $\tilde f_n$ tends to $f$ in ${\mathcal L}^2(R)$ and 
we only need to check that $|\nabla (f_n-\tilde f_n)|$ tends to $0$ in ${\mathcal L}^2(R)$.
Setting $R_n=[0,1/n]\times U$, we write
\begin{eqnarray*}
\int_{R}|\nabla (\tilde f_n-f_n)|^2\, d\mu
&=&\int_{R_n} \Bigl(|\partial_sf_n|^2 + |\nabla_M (\tilde f_n-f_n)|^2 
                                                \Bigr)\,d\mu\\
&\le& C \int_{R_n}\Bigl(|\partial_sf_n|^2+|\nabla_M f_n\Bigr)|^2 \,d\mu
      + \frac{C}{n}\int_U|\nabla_M \tilde f_n(1/n,x)|^2\,dx .
\end{eqnarray*}
It is plain that
$$
\int_{R_n} \Bigl(|\partial_sf_n|^2+|\nabla_M f_n)|^2 \Bigr)\,d\mu
\rightarrow 0\,.
$$
Moreover
\begin{eqnarray*}
\frac{1}{n}\int_U|\nabla_M\tilde f_n(1/n,x)|^2\,dx
&\le & \frac{1}{n}\int_U
         \Bigl|\,n\int_{1/n}^{2/n} |\nabla_M f(s,x)| \, ds\Bigr|^2\,dx\\
&\le & \int_{1/n}^{2/n}\int_U |\nabla_M f(s,x)|^2 \,ds\,dx \rightarrow 0\,.
\end{eqnarray*}
The functions $\tilde f_n$ satisfy $\partial_s\tilde f_n=0$ in 
$[0,1/n)\times U$ but are not smooth. To obtain smooth functions approximating 
$f$ with the desired property,  extend $f$ and $\tilde f_n$ by symmetry
to $R^*=(-l,l)\times U$, that is, $f(-s,x)=f(s,x)$ and 
$\tilde f_n(-s,x)=\tilde f_n(s,x)$. Obviously, 
$\|\tilde f_n-f\|_{\mathcal W^1(R^*)}\rightarrow 0$. For each $n$, let $\theta_n$ be a 
smooth non-negative function with integral $1$ and support in the ball of 
radius less than $1/(5n)$ around $(0,0)$ in $(-l,l)\times U$.
Consider $h_n=\theta_n*\tilde f_n$ ($*$ denoting convolution).
Now, the restriction of $h_n$ to $[0,l)\times U$ is a smooth function
which satisfies $\partial_s h_n=0$ in a neighbourhood of $\{0\}\times U$
and approximates $f$ in $\mathcal W^1(R)$. Indeed,
$\theta_n*f\rightarrow f$ in $\mathcal W^1(R^*)$, and
$\|\theta_n*(\tilde f_n-f)\|_{\mathcal W^1(R^*)}\le 
\|\tilde f_n-f\|_{\mathcal W^1(R^*)}\rightarrow 0\,.$
\end{proof}

The Dirichlet form structure on a strip complex $\GM$
is based on the choice of
\begin{itemize}
\item[(a)] the geometry
determined by $l,\phi\,$, and
\item[(b)] the measure $\mu$ determined by $\psi$.
\end{itemize}
The following definition takes this into account to introduce isometries
that are compatible with this additional structure.
\begin{dfn} Let $\GMa$ and $\GMb$ be two strip complexes equipped respectively 
with $\phi_i,\psi_i$ and the associated measures $\mu_i$, $i=1,2$ as above.
An isometry (or local isometry, with obvious modifications)
$\jm:\GMa\rightarrow \GMb$ is called \emph{measure-adapted} if there is a 
positive constant $c(\jm)$ such that, for any compact set $A\subset \GMb\,$,
$$
\mu_1(\jm^{-1}(A))=c(\jm)\,\mu_2(A).
$$
\end{dfn}

\begin{rmk} If $\jm$ is a measure-adapted isometry and $f_1= f_2\circ \jm$, 
where $f_2\in \mathcal W^1(\GMb)$, then $f_1\in  \mathcal W^1(\GMa)$ and
$$
\mathcal E_1(f_1,f_1)= c(\jm)\,\mathcal E_2(f_2,f_2).$$
\end{rmk}

\begin{exa} For any $\pp\in \{1,2\dots\}$ and $\qq>1$, the treebolic space 
$\HT(\pp,\qq)$ (equipped with its stripwise hyperbolic geometry, as described
in Section \ref{geometry}) admits a large group of isometries 
(see the introduction). These isometries are all measure-adapted
whenever $\HT(\pp,\qq)$ is equipped with any of the measures 
$\mu_{\alpha,\beta}$ ($\alpha\in \R$, $\beta>0$) defined in \eqref{muab}.
\end{exa}

\subsection*{F. The Laplacian and the heat equation on strip complexes}
Consider a srip complex $\GM$ where $(\Ms,g)$ is a Riemannian manifold, and 
equip $\GM$ with the data $(l,\phi,\psi)$, 
where $\phi,\psi\in \mathcal S^\infty(\Gamma^o)$ as in \S3.D and \S3.E. 
Let $\mu$ the associated measure. For simplicity, we write 
$\mathcal L^p(\GM)=\mathcal L^p(\GM,\mu)$, 
$\mathcal W^1_0(\GM)=\mathcal W^1_0(\GM,\mu)$
and $\mathcal W^1(\GM)=\mathcal W^1(\GM,\mu)$.

By the general theory of Dirichlet forms, there is a self-adjoint operator
$$
\bigl(\Delta,\mbox{Dom}(\Delta)\bigr)
$$
on ${\mathcal L}^2(\GM)$ which we call the Laplacian on $\GM$ and which is
defined as follows.

\begin{dfn} Set 
$$
\mbox{Dom}(\Delta)=
\bigl\{f\in \mathcal W^1_0(\GM): \text{there is} \;C_f\;\text{such that}\;
\mathcal E(f,h)\le C_f\|h\|_2 \;  \text{for all}\;h\in \mathcal W^1_0(\GM)
\bigr\}.
$$
For $f\in \mbox{Dom}(\Delta)$, there exists
a unique $u\in {\mathcal L}^2(\GM)$ such that $\mathcal E(f,h)=-\int u \,h\,d\mu$ 
for all $h\in {\mathcal L}^2(\GM)$ and we set
$$
\Delta f=u.
$$
\end{dfn}

Since the measure $\mu$ will be fixed most of the time, we will often
omit it in our notation.
The operator $\Delta$ with domain $\mbox{Dom}(\Delta)$ is the
infinitesimal generator of a strongly continuous
semigroup of self-adjoint contractions  $\{ H_t=e^{t\Delta} : t\ge 0 \}$,
on ${\mathcal L}^2(\GM)$ which has the Markov property:
$$
f\in {\mathcal L}^2(\GM)\,, \;\;0 \le f\le 1 
\;\;\Longrightarrow \;\;0\le H_t f\le 1.
$$
It follows that $H_t$ extends to a contraction on each  space ${\mathcal L}^p(\GM)$, 
$1\le p\le \infty$. For $1\le p<\infty$, the family
$\{H_t : t\ge 0\}$ is a strongly continuous semigroup on ${\mathcal L}^p(\GM)$.
We call $\{H_t : t > 0\}$ the \emph{heat semigroup} on $\GM$ 
(more precisely, on $(\GM;l,\phi,\psi)$). 

The following is immediate by inspection.

\begin{pro}Let $\GMa$ and $\GMb$ be two strip complexes, each  equipped 
with data $l_i\,,\phi_i\,,\psi_i\,$ ($i=1,2$) as above.
Let $\mu_i$ and $\bigl(\Delta_i,\mbox{\em Dom}(\Delta_i)\bigr)$, $i=1,2$,
be the associated measures and Laplacians.  If $\jm:\GMa\rightarrow \GMb$
is a measure-adapted isometry then
$$
\text{for all}\;\;f_2\in \mbox{\em Dom}(\Delta_2)\,,\quad
f_1=f_2\circ \jm \in \mbox{\em Dom}(\Delta_1)\;
\mbox{ and }\; \Delta_1 f_1 = (\Delta_2 f_2)\circ \jm\,.
$$
Also,
$$
\text{for all}\;t>0\;\text{and}\;f_2\in {\mathcal L}^2(\GMb)\,,\quad
f_1=f_2\circ \jm\in {\mathcal L}^2(\GMa)\;\mbox{ and }\; 
H_{1,t} f_1 = (H_{2,t} f_2)\circ \jm\,.
$$
\end{pro}

In the general theory of regular strictly local Dirichlet forms
$\bigl(\mathcal E,\mbox{Dom}(\mathcal E)\bigr)$, one introduces a notion of 
intrinsic distance. In the present setting, this definition reads
$$
\widetilde{\rho}(x,y)=\sup\left\{f(x)-f(y): f\in
\mathcal C(\GM)\cap \mathcal W^1(\GM), \;\;|\nabla f|\le 1\right\}.
$$
It is not hard to see that this intrinsic distance coincides with
the distance $\rho$ introduced earlier.

\begin{dfn} Let $\Omega$ be an open set in $\GM$. Set
$$
\mathcal W^1_c(\Omega)=\{f\in \mathcal W^1(\GM):
 f\; \mbox{is compactly supported in}\; \Omega\}
$$
and
$$
\mathcal W^1_{\mbox{\tiny loc}}(\Omega)
=\left \{f\in {\mathcal L}^2_{\mbox{\tiny loc}}(\Omega):
{\displaystyle \text{for every compact}\; K \subset \Omega\;\text{there is}\quad
\atop
\tilde f\in \mathcal W^1(\GM)\; \text{such that}\; f|_K=\tilde{f}|_K\; \mbox{a.e} }
\right\}.
$$
\end{dfn}

Fix an open set $\Omega$ and consider the topological vector spaces
$\mathcal W^1_c(\Omega)\subset \mathcal W^1_0(\GM)\subset {\mathcal L}^2(\GM)$ 
and their duals
${\mathcal L}^2(\GM)\subset \mathcal W^1_0(\GM)^*\subset \mathcal W^1_c(\Omega)^*$.

\begin{dfn} Let $\Omega$ be an open set in $\GM$. Let $f\in \mathcal W^1_c(\Omega)^*$.
We say  that a function $u$ is a weak solution of the  equation
$\Delta u=f$ in $\Omega$ if
\begin{itemize} \item $u\in \mathcal W^1_{\mbox{\tiny loc}}(\Omega)$, and 

\smallskip

\item $\mathcal E(u,h) = -f(h)$ for all
$h\in \mathcal W^1_c(\Omega)$.
\end{itemize}
\end{dfn}

Observe that $f(h) $ above is well defined since $f\in \mathcal W^1_c(\Omega)^*$ and
$h\in \mathcal W^1_c(\Omega)$. Observe also that if $f$ is represented by a 
locally integrable function on $\Omega$ (again called $f$) and if $u$ is such 
that there exists $\tilde{u}\in \mbox{Dom}(\Delta)$ satisfying
$u=\tilde{u}|_\Omega$ then  $u$ is a weak solution of $\Delta u=f$ in
$\Omega$ if and only if $(\Delta \tilde{u})|_\Omega=f$.

Given a Hilbert space $H$ and an interval $I$, let ${\mathcal L}^2(I\rightarrow H)$ be 
the Hilbert space of those functions $f: I\to H$ such that
$$
\|f\|_{{\mathcal L}^2(I\rightarrow H)}=\left(\int_I\|f(t)\|_H^2\,dt\right)^{1/2}<\infty\,.
$$
Let $\mathcal W^1(I\rightarrow H)\subset {\mathcal L}^2(I\rightarrow H)$
be the Hilbert space of those functions
$f: I\to H$ in ${\mathcal L}^2(I\rightarrow H)$ whose distributional time derivative
$f'$  can be represented by functions in ${\mathcal L}^2(I\rightarrow H)$,
equipped  with the norm
$$
\|f\|_{\mathcal W^1(I\rightarrow H)}=
\left(\int_I\bigl(\|f(t)\|_H^2+\|f'(t)\|_H^2\bigr)\,dt\right)^{1/2}<\infty\,.
$$
Given an open  time interval $I$, set
$$
\mathcal F(I\times \GM)={\mathcal L}^2 \bigl(I\rightarrow \mathcal W^1_0(\GM)\bigr)
\cap
\mathcal W^1\bigl(I\rightarrow  \mathcal W^1_0(\GM)^*\bigr).
$$
This notation is justified by the inclusions
$\mathcal W^1(\GM) \subset {\mathcal L}^2(\GM) 
= {\mathcal L}^2(\GM)^* \subset \mathcal W^1(\GM)^*$, compare with 
 \cite{Stu1}, \cite{Stu2}, \cite{Stu3}. While in these definitions
it was convenient to consider $f(t)$ as a function on $\GM$ for each
$t \in I$, we shall usually prefer the notation $f(t,\cdot)$, where we think
of $f$ as a function on $I \times \GM$. 

Given an open time interval $I$ and an open set $\Omega\subset \GM$
(both nonempty), let
\begin{equation}\label{Floc}
\mathcal F_{\mbox{\tiny loc}}(I\times \Omega)
\end{equation}
be the set of all functions $f:I\times \Omega\rightarrow \R$
such that, for any open interval $I'\subset I$ relatively compact in $I$
and any open subset $\Omega'$ relatively compact in $U$ there exists a function
$f^{\#}\in \mathcal F(I\times \GM)$  satisfying $f=f^{\#}$ a.e. in
$I'\times \Omega'$.
Finally, let
$$
\mathcal F_c(I\times \Omega)=\{f\in \mathcal F(I\times \GM): f(t,\cdot)
\; \mbox{has compact support in}\; \Omega\; \mbox{for a.e.}\;t\in I\}.
$$

\begin{dfn}
Let $I$ be an open time interval.
Let $\Omega$ be an open subset in $\GM$ and set $Q=I\times \Omega$.
A function $u:Q\mapsto \R$ is a \emph{weak (local) solution of the heat 
equation} $(\partial_t-\Delta)u=0$ in $Q$ if
\begin{enumerate}
\item $u\in \mathcal F_{\mbox{\tiny loc}}(Q)$, and
\item for any open interval $J$ relatively compact in $I$
and any  $f\in \mathcal F_c(Q)$,
$$ 
\int_J\int_U  f \,\partial_t u \,d\mu\, dt +
\int_{J}\mathcal E\bigl(f(t,\cdot),u(t,\cdot)\bigr)\,dt= 0.
$$
\end{enumerate}
\end{dfn}

The following proposition follows from the relevant definitions by inspection.

\begin{pro}
Let $\GMa$ and $\GMb$ be two strip complexes, each equipped with data
$\phi_i,\psi_i$, $i=1,2$, as above.
Let $\mu_i$ and $\bigl(\Delta_i,\mbox{\em Dom}(\Delta_i)\bigr)$, $i=1,2$,
be the associated measures and Laplacians.
Let $\,\jm:\Omega_1\rightarrow \Omega_2\,$, where $\Omega_1\subset \GMa$ and 
$\Omega_2\subset \GMb\,$, be a measure-adapted local isometry between 
the open sets $\Omega_1\,$ and $\Omega_2\,$. 
\begin{itemize}
\item If $f\in \mathcal W^1_c(\Omega_2)^*$  and $u_2$ is a weak solution of
$\Delta_2 u=f_2$ in $\Omega_2$ then 
$f_1(h)= f_2(h\circ \jm^{-1})\in \mathcal W^1_c(\Omega_1)^*$
and $u_1=u_2\circ \jm$ is  a weak solution of $\Delta_1 u_1=f_1\,$.

\smallskip

\item For any time interval $I$, if $u_2$ is a weak solution of the heat 
equation on $\GMb$
in $Q_2=I\times \Omega_2\,$, then $u_1=u_2\circ \jm$ is a weak solution of 
the heat equation on $\GMa$ in $Q_1=I\times \Omega_1\,$.
\end{itemize}
\end{pro}

\section{Basic properties of the heat semigroup}\label{heat}

In this section and the next, $\GM$  is a fixed strip complex based on graph 
$\Gamma$ and a Riemannian manifold 
$(\Ms,g)$. Furthermore, $\GM$ is equipped with data $(l,\phi,\psi)$, where 
$\phi,\psi\in \mathcal S^\infty(\Gamma^o)$,  
the associated distance $\rho$ and measure $\mu$, the Dirichlet form 
$\bigl(\mathcal E,\mathcal W^1_0(\GM)\bigr)$ and the corresponding 
Laplacian $\Delta$ and heat semigroup
$\{H_t=  e^{t \Delta} : t\ge 0 \}$. See \S3.D--\S3.F. 

Because of the singular nature of strip complexes, the local regularity 
properties of weak solutions of the Laplace
or heat equations are a non-trivial and crucial issue.

\begin{thm}\label{th-locDP}
For any compact set $K \subset \GM$, there exist 
$r_K > 0$ and constants $D_K\,$,  $P_K$ such that for all $\xi\in K$, 
$r\in (0,r_K)$ the following properties hold.
\begin{itemize}
\item $V(\xi,r)\le D_K\, V(\xi,2r)$, and 

\smallskip
 
\item for every $f\in \mathcal W^1(B)$,
$$
\int_{B}|f-f_B|^2\,d\mu \le P_K \,r^2\int_{B}|\nabla f|^2\,d\mu\,,
$$
where $B=B(\xi,r)$ and $f_B=\frac{1}{\mu(B)}\int_B f\,d\mu$.
\end{itemize}
\end{thm}

\begin{proof} The first property is clear by inspection because of the
continuity of $\phi,\psi$, the Riemannian nature of $\Ms$, and the fact that
the underlying graph $(V,E)$ is locally finite. The second property,
i.e., the Poincar\'e inequality, is more delicate to prove. First,
such a (localized) Poincar\'e inequality holds on $\Ms$
(i.e., on any Riemannian manifold). See, e.g.,  
\cite[5.6.3]{SalB}. This applies to any strip $S_e$ equipped with the 
$\phi$-structure. By continuity and
positivity of  $\psi$, the desired local Poincar\'e inequality holds
on balls that are contained in the interior of a strip.

The same is true when we have a vertex $v$ with $\deg(v)=1$, so that
$\Ms_v$ sits at the boundary of a unique strip $S_e\,$, and $B$ is contained
in the half-open strip $S_e^o \cup \Ms_v\,$.  

By classical arguments, it thus suffices to prove the stated result
assuming that the center $\xi$ belongs to a bifurcation manifold $\Ms_v\,$,
where $\deg(v) \ge 2$. 
We can further assume that $r$ is small enough so that the ball $B=B(\xi,r)$ is 
contained in
$X_v\,$, the star of strips around $v$.

The crucial observation is that  for any pair of edges $e, e'\in E_v\,$,
the open set $X_v^{e,e'}=\Ms_v\cup S^o_{e}\cup S^o_{e'}$ equipped with the 
$\phi$-structure is locally bi-Lipschitz equivalent to a smooth Riemannian 
manifold $I\times \Ms$, where the interval $I$ corresponds to
$\{v\} \cup I_e^o \cup I_{e'}^o\,$.

Therefore, setting $B=B(\xi,r)$ and $B_{e,e'}=B\cap X_v^{e,e'}\,$ 
the following Poincar\'e inequalities hold:
\begin{eqnarray*}
\int_{B_{e,e'}}|f-f_{B_{e,e'}}|^2\,d\mu
&\le& C_K\, r^2\int_{B_{e,e'}}|\nabla f|^2\,d\mu
\quad \text{for all}\;f\in \mathcal W^1(B)\,,\;e,e'\in E_v\,,\\
\text{where}\;\;
f_{B_{e,e'}}
&=&\frac{1}{\mu(B_{e,e'})}\int_{B_{e,e'}}f\,d\mu\,.
\end{eqnarray*}
Now choose and fix an edge $e\in E_v$ so that $\mu(B_e)$ is maximal
among all edges in $E_v\,$, where $B_{e}=B\cap S^o_e$. We set
$f_{B_e}=\frac{1}{\mu(B_e)}\int_{B_e}f\,d\mu$. Then
$$
\max_{e'\in E_v} \mu(B_{e,e'}) \le 2 \mu(B_e).
$$
Then we can estimate
\begin{eqnarray*}
|f_{B_e}-f_{B_{e,e'}}|
&=&
\left|\frac{1}{\mu(B_e)\mu(B_{e,e'})} \int_{B_e} \int_{B_{e,e'}}
\bigl(f(\eta)-f(\zeta)\bigr) \, d\mu(\eta)\, d\mu(\zeta)\right|\\
&\le&
\frac{2}{\mu(B_{e,e'})^2} \int_{B_{e,e'}}\int_{B_{e,e'}}
\bigl|f(\eta)- f(\zeta)\bigr|\, d\mu(\eta)\, d\mu(\zeta)\\
&\le&
 \left(\frac{4}{\mu(B_{e,e'})^2} \int_{B_{e,e'}}\int_{B_{e,e'}}
\bigl|\bigl(f(\eta)- f_{B_{e,e'}}\bigr)-\bigl(f(\zeta)- f_{B_{e,e'}}\bigr)\bigr|^2
\,d\mu(\eta)\, d\mu(\zeta)\right)^{\! 1/2}\\
&\le& 
\left(\frac{8}{\mu(B_{e,e'})}\int_{B_{e,e'}} |f-f_{B_{e,e'}}|^2\,d\mu
        \right)^{\!1/2}\\
&\le& 
 \left(\frac{8\,C_K \,r^2}{\mu(B_{e,e'})}\int_{B_{e,e'}}
|\nabla f|^2d\mu\right)^{\! 1/2}\,.
\end{eqnarray*}
In the last inequality, we have used the Poincar\'e inequality on 
$B_{e,e'}\,$. Next,
$$
\int_{B}|f-f_B|^2\,d\mu =\min_{c\in \R}
\int_{B}|f-c|^2 \, d\mu \le \int_{B}|f-f_{B_e}|^2\,d\mu,
$$
and
\begin{eqnarray*}
\int_{B}|f-f_{B_{e}}|^2\,d\mu 
&\le&
\sum_{e'\in E_v\setminus \{e\}} \int_{B_{e,e'}}|f-f_{B_e}|^2\, d\mu \\
&\le&
2\sum_{e'\in E_v\setminus \{e\}}
\left(\int_{B_{e,e'}}|f-f_{B_{e,e'}}|^2\,d\mu 
 + \mu(B_{e,e'})\,|f_{B_{e,e'}}-f_{B_e}|^2\right)\\
&\le& 
18\, C_K \, r^2\sum_{e'\in E_v\setminus \{e\}} 
\int_{B_{e,e'}}|\nabla f|^2\,d\mu \\
&\le& 
18\, C_K \, r^2 \bigl(\deg(v)-1\bigr)\int_{B}|\nabla f|^2\,d\mu.
\end{eqnarray*}
This is the desired Poincar\'e inequality when $\xi\in K\cap M_v$.
From this, for all $\xi\in K$ and $r\in (0\,,\,r_K)$,
elementary considerations give that for all $f\in \mathcal W^1(2B)$,
$$
\int_{B}|f-f_B|^2\, d\mu\le P_K \, r^2 \int_{2B}|\nabla f|^2\,d\mu\,,
$$
where $2B=B(\xi,2r)$. Now, it is well known (but not so elementary) that
this suffices to obtain the desired Poincar\'e inequality where
$f\in \mathcal W^1(2B)$ and $\int_{2B}|\nabla f|^2d\,\mu$ are replaced by
$f\in \mathcal W^1(B)$ and $\int_{B}|\nabla f|^2\,d\mu$. See \cite[5.3]{SalB}.
Compare also with \cite{EF} and \cite{PSC}.
\end{proof}

Theorem \ref{th-locDP} has far reaching consequences.
The next three theorems follow from the arguments of 
\cite{Stu1}, \cite{Stu2}, \cite{Stu3},
which are based on Moser iteration techniques and thus,
in the present situation, rely heavily
on Theorem \ref{th-locDP}. See also \cite{SalB} and {\sc Biroli and
Mosco}~\cite{BM}.

\begin{thm}\label{th-locHR}
Referring to the general setting of this section, the heat semigroup has
the following properties.
\begin{itemize}
\item 
For any open interval $I$ and compact intervals 
$J$, $J'$ of $I$ with $\max J < \min J'$ and for any  connected open set 
$\Omega \subset \GM$ and compact  $K\subset \Omega$,
there are positive constants $\alpha_1=\alpha_1(I,\Omega,J,K)$, 
$C_1=C_1(I,\Omega,J,K)$ and $C_2=C_2(I,\Omega,J,J',K)$ such that any
weak solution $u$ of the heat equation $(\partial_t-\Delta)u=0$ in 
$I\times \Omega$ admits a continuous version which satisfies
$$
\sup
\left\{ \frac{|u(t,\xi)-u(s,\zeta)|}{(|t-s|^{1/2}+\rho(\xi,\zeta))^{\alpha_1}} :
(t,\xi),(s,\zeta)\in J\times K
\right\} \le C_1\, \sup_{I\times \Omega} |u|\,,
$$
and, if $u$ is non-negative,
$$
\sup_{J\times K} u \le C_2\, \inf_{J'\times K} u.
$$
\item 
The heat diffusion semigroup $\{H_t=  e^{t \Delta} : t > 0 \}$ admits a
continuous kernel $(t,\xi,\zeta)\mapsto h(t,\xi,\zeta)$ --
which we call the \emph{heat kernel} of\/ $\GM$ -- so that
$$
H_tf(\xi)=\int_{\GM}h(t,\xi,\zeta)\,f(\zeta)\,d\mu(\xi).
$$
The heat kernel is symmetric in $\xi,\zeta$.
\item 
For each $\xi,\zeta\in \GM$, the function $t\mapsto h(t,\xi,\zeta)$ is in
$\mathcal C^\infty\bigl((0,\infty)\bigr)$, and for each $\zeta\in \GM$,
the function $(t,\xi)\mapsto \partial_t^k h(t,\xi,\zeta)$ is a weak solution of
the heat equation in  $(0,\infty)\times \GM$. Moreover,
$(t,\xi,\zeta)\mapsto \partial_t^k h(t,\xi,\zeta)$ is a continuous function
on $(0,\infty)\times \GM\times\GM$.
\item For any fixed compact $K\subset \GM$, $\zeta_0\in K$,
compact time interval $I=[a\,,\,b]\subset (0\,,\,\infty)$ and integer $k$, 
there are positive constants $\alpha_2=\alpha_2(I,K,k)$ and 
$C_3=C_3(I,K,k)$ such that, for all $\xi\in \GM$, we have
$$
\sup \bigl\{|\partial_t^k h(t,\zeta,\xi) : t\in I,\zeta\in K| \bigr\}
\le C_3\, h(2b,\zeta_0,\xi)
$$
and
$$
\sup
\left\{ \frac{|\partial^k_th(t,\zeta,\xi)-\partial^k_th(t,\zeta',\xi)|}
{\rho(\zeta',\zeta)^{\alpha_2}} : t\in I,\zeta,\zeta'\in K \right\}
\le C_3 \, h(2b,\zeta_0,\xi).
$$
\item Each operator $H_t$, $t>0$, sends bounded measurable functions to
continuous bounded functions, that is, $H_t {\mathcal L}^\infty(\GM)\subset
\mathcal C_b(\GM)$ for any $t>0$.
\end{itemize}
\end{thm}

Note that no global results can be obtained
under the present very general hypotheses. In particular, we have no bound on
the volume of large balls,  and stochastic completeness is not guaranteed. 
That is, it may very well occur that 
$\int h(t,\xi,\zeta)\,d\mu(\zeta)<1$ for some $t,\xi$.
Indeed, we have so far not even assumed the completeness 
of $(\GM,\rho)$, but will do so next.

\begin{thm}\label{th-BdUniq}
 Assume that $(\GM,\rho)$ is complete and that 
$$
\int_1^\infty \frac{r\,dr}{\ln V(\xi_0,r)}=\infty\,.
$$
Then uniqueness of the bounded Cauchy problem holds on $(0,T)\times \GM$
for the heat equation. More precisely, if $u: (0,T)\times \GM$ is a weak
solution of the heat equation on $(0,T)\times\GM$ which is bounded
and  satisfies $\lim_{t\rightarrow 0} u(t,\xi)=0$ $\mu$-almost everywhere,
then $u=0$ on $(0,T)\times \GM$. In particular,
the semigroup $\{H_t =  e^{t \Delta} : t > 0 \}$, is conservative, 
that is, $e^{t\Delta}\mathbf 1=\mathbf 1$.
\end{thm}

In the next theorem, we also assume that $(\GM,\rho)$ is complete,
and make uniform local assumptions on the geometry of $\GM$ that
allow us to obtain more quantitative results.

\begin{thm}\label{th-UR} 
Assume that $(\GM,\rho)$ is complete and that there are constants $D,P,r_0 > 0$
such that
\begin{enumerate}
\item[(i)] for any $\xi\in \GM$ and $r\in (0\,,\,r_0)$, we have the doubling property
$V(\xi,r)\le D \,V(\xi,2r)$, and
\item[(ii)] for any $\xi\in \GM$  and $r\in (0\,,\,r_0)$, setting $B=B(\xi,r)$,
$$
\int_B|f-f_B|^2\, d\mu\le P\, r^2\int_B|\nabla f|^2\,d\mu \quad
\text{for every}\;\;f\in \mathcal W^1(B)\,, \quad\text{where}\;\;
f_B=\frac{1}{\mu(B)}\int_B f\,d\mu\,.
$$
\end{enumerate}

\smallskip

Then the following properties hold.

\smallskip

\begin{enumerate}
\item 
For fixed $R>0$ there are positive constants
$\alpha$, $C_4$ and $C_5$ (depending only on $R$)
such that for all $\xi\in \GM$, $r\in (0\,,\,R)$, any
weak solution $u$ of the heat equation $(\partial_t-\Delta)u=0$ in
$Q=(0,4r^2)\times B(\xi,2r)$ satisfies
$$
\sup
\left\{ \frac{|u(t,\xi)-u(s,\zeta)|}{(|t-s|^{1/2}+\rho(\xi,\zeta))^\alpha} :
(t,\xi),(s,\zeta)\in Q'\right\}
\le \frac{C_4}{r^{\alpha}}\,\sup _{Q}|u|
$$
and, if $u$ is non-negative,
$$
\sup_{Q_-} u\le C_5\, \inf_{Q_+} u\,, \quad\text{where}
$$
$Q'=(r^2\,,\,3r^2)\times B(\xi, r)$, $Q_-=(r^2\,,\,2r^2)\times B(\xi,r)$
and $Q_+=(3r^2\,,\,4r^2)\times B(\xi,r)$.

\medskip

\item 
For any fixed integer $k \ge 0$ and $\epsilon\in (0\,,\,1)$
there is a constant $C_{k,\epsilon}$ such that for all $t>0$ and all
$\xi,\zeta\in \GM$, with $\alpha$ as above,
$$
|\partial_t^k h(t,\xi,\zeta)|
\le
\frac{C_{\epsilon,k}}{t^k\,V(\xi,\min\{1,\sqrt{t}\})}
\exp\left(
-\frac{\rho(\xi,\zeta)^2}{4(1+\epsilon)t}\right).
$$
Moreover,
$$
|\partial_t^k h(t,\xi,\zeta)|\le \frac{C_{\epsilon, k}}{\min\{1,t\}^k}\,
h\bigl((1+\epsilon)t,\xi,\zeta\bigr)
$$
and, for all $\zeta'$ with $\rho(\zeta,\zeta')\le \min\{1,\sqrt{t}\}$,
$$
|\partial_t^k h(t,\xi,\zeta)-\partial_t^k h(t,\xi,\zeta')|
\le \frac{C_{\epsilon,k}\,\rho(\zeta,\zeta')^\alpha}
{(\min\{1,\sqrt{t}\})^{\alpha+k/2}}
\,h\bigl((1+\epsilon)t,\xi,\zeta\bigr).
$$
\end{enumerate}
\end{thm}

Concerning the growth of the volume of large balls, we point out that the
hypothesis that the  volume doubling property holds locally uniformly
as in Theorem \ref{th-UR} implies that 
$$
V(\xi,r)\le e^{Cr/r_0}V(\xi,r_0) \quad\text{for all}\;\;r \ge r_0\,,
$$
see \cite[Lemma 5.2.7]{SalB}. We collect three of the main features.

\begin{cor} \label{cor-Htcons}
Under the hypotheses of {\em Theorem \ref{th-UR}}, the following properties
hold for the heat semigroup $\{H_t = e^{t\Delta}: t > 0\}$.
\begin{enumerate}
\item It is conservative (stochastically complete), that is, 
$e^{t\Delta}\mathbf 1=\mathbf 1$.
\item It sends ${\mathcal L}^\infty(\GM)$ into $\mathcal C_b(\GM)$.
\item It sends $\mathcal C_0(\GM)$ into itself.
\end{enumerate}
\end{cor}

The next corollary concerns global non-negative solutions of the heat equation.

\begin{cor} \label{cor-PosUniq}
Under the hypotheses of {\em Theorem \ref{th-UR}}, there exists a constant 
$C$ such that any non-negative weak solution $u$ of the heat equation on
$(0\,,\,T)\times \GM$ satisfies
$$
u(s,\xi)\le u(t,\zeta)\exp\Bigl(C\bigr(1+ t/s +
\rho(\xi,\zeta)^2/(t-s)\bigr)\Bigr) \quad\text{for all}\;\; \xi,\zeta\in \GM,
\; 0<s<t<T.
$$

Moreover, uniqueness of the positive Cauchy problem holds on
$(0\,,\,T)\times \GM$ for the heat equation.
More precisely, if $u: (0\,,\,T)\times \GM$ is non-negative and is a weak
solution of the heat equation on $(0\,,\,T)\times\GM$ then there exist a
non-negative Borel measure $\si$ on $\GM$ and $a>0$ such that
$$
\int_{\GM} e^{-a\,\rho(\xi_0,\xi)^2}\,d\si(\xi)<\infty
$$
for some (equivalently, any) $\xi_0\in \GM$, and
$$
u(t,\xi)=\int_{\GM} h(t,\xi,\zeta) \,d\si(\zeta) \quad\text{for all}\;\;
(t,\xi)\in (0,T)\times \GM\,.
$$
In particular, if $u$ is a non-negative weak solution of the heat equation
in $(0\,,\,T)\times \GM$ and there is some 
$u_0\in  {\mathcal L}^1_{\mbox{\em \tiny loc}}(\GM)$ such that
$$
\lim_{t\rightarrow 0}\int_{\GM} u(t,\cdot)\,f\,d\mu=\int_{\GM} u_0\, f\,d\mu
\quad\text{for all}\;\;f\in \mathcal C^\infty_c(\GM),
$$
then $\displaystyle u(t,\xi)=\int_{\GM} h(t,\xi,\cdot) u_0\, d\mu\,$.
\end{cor}

\begin{proof} $\;$ See {\sc Ancona and Taylor} \cite{AT},
{\sc Aronson}~\cite{Ar1}, \cite{Ar2}, {\sc Grigor'yan}~\cite[Th. 6.2]{Gri-surv}, 
\cite[Sec. 3]{Stu1} and \cite[Sec. 5.5.2]{SalB}.
\end{proof}

Next, we give some relatively simple sufficient conditions which imply that
the hypotheses of Theorem \ref{th-UR} are satisfied.

\begin{pro} Assume the following.
\begin{itemize} 
\item 
The manifold $(\Ms,g)$ is complete and satisfies the doubling property and 
${\mathcal L}^2$-Poincar\'e inequality at all scales, that is, there are positive 
constants $D_{\Ms}$ and $P_{\Ms}$ such that, for every $x_0\in \Ms,\;r>0$,
$$
V_{\Ms}(x_0,r)\le D_{\Ms}\,V_{\Ms}(x_0,2r),
$$
where $V_{\Ms}(x_0,r)$ is the Riemannian volume of the geodesic ball 
$B=B_M(x_0,r)$ of  radius $r$ around $x_0$ in $\Ms$, and
$$
\int_{B}|f-f_B|^2\,dx 
\le P_{\Ms}\,r^2 \int_B |\nabla_{\Ms} f|^2\,dx\quad\text{for all}
\;\;f\in  \mathcal W^1(B),
$$
where $f_B$ is the average of $f$ over $B$, and $dx$ is the volume element
of $\Ms$.

\item 
There are finite positive constants $c_0$ and $C_0$ such that
$$
\int_{e^-}^{e^+}\sqrt{\phi_e(s)}\, ds\ge c_0 \quad\text{for every}\;e \in E\,,
\quad\text{and}\quad \deg(v) \le C_0 \quad\text{for every}\;v\in V\,.
$$
Moreover, for any finite interval
$I\subset \Gamma^1$ with $\int_I\sqrt{\phi(s)}\,ds \le c_0\,$,
$$
\frac{\max_I\phi}{\min_I\phi }\le C_0 \AND 
\frac{\max_I\psi}{\min_I\psi}\le C_0.
$$
\end{itemize}

Under these hypotheses, $(\GM,\rho)$ is complete, and there are constants 
$D,P,r_0$ such that the properties (i) and (ii) of \emph{Theorem \ref{th-UR}}
hold.
\end{pro}

\begin{proof} Completeness follows clearly from Lemma \ref{lem-GMcomp}.
Moreover, under the above hypotheses on $\phi,\psi$, for any fixed
$r_0\,$, the functions $\phi$ and $\psi$ behave like constant functions 
(that is, there is $c=c(r_0) > 0$ such that $c \le \phi,\psi \le 1/c$) on
any ball of radius $r_0$ in $\GM$.
This means that the geometry of $\GM$ in such a ball $B$ is
comparable to the product of a piece of $\Gamma_1$ scaled by a constant factor
$\phi_B$ (corresponding to the size of $\phi$ in the ball in question)
and $(\Ms,\phi_B \,g)$.
The uniform local doubling property thus follows from
the global doubling property on $M$ and the fact that $\phi$ and $\psi$ are
approximately constant in $B$. The uniform local Poincar\'e inequality
follows by the argument used in the proof of Theorem \ref{th-locDP} that 
can now be carried through up to a uniformly fixed scale.
\end{proof}

\begin{exas} (a) Let $(\Gamma^1,l)$ be a metric graph as above with
$\min_{e\in E}\{l_e\}>0$.
Suppose that $\psi\in\mathcal S^\infty(\Gamma^1)$ has the property that for any
interval $I\subset \Gamma^1$ of length $1$, one has $\max_I\psi/\min_I\psi\le
C$ for some positive $C$, and that $\max_V \deg(v) < \infty\,$.

Then the weighted $1$-complex $(\Gamma^1,l,\psi(s)\, ds)$ satisfies the
hypotheses of Theorem \ref{th-UR}.
More generally, for any $k=0,1,2,\dots$, the strip complex $\GM$ with 
$\Ms=\R^k$, $\phi\equiv 1$ and $\psi$ as above, satisfies the
hypotheses of Theorem \ref{th-UR}.

\smallskip

\noindent
(b) The treebolic space $\HT$ equipped with any one of the forms
$\bigl(\mathcal E_{\alpha,\beta}\,,\mathcal W^1_{\alpha,\beta}({\HT})\bigr)$
satisfies the hypotheses of Theorem
\ref{th-UR}. This follows from the local result and the fact that
there is a transitive group of measure adapted isometries for any one
of these structures.
\end{exas}

\section{Smoothness of weak solutions}\label{smoothweak}

Throughout this section, we keep the setting and notation of  Section 
\ref{heat}. 

\subsection*{A. Harmonic functions}
By the general theory of Dirichlet forms, there is a Hunt  process with 
continuous sample paths defined for every starting point $\xi\in\GM$ 
associated with the semigroup $H_t=e^{t\Delta}: {\mathcal L}^2(\GM)\to 
{\mathcal L}^2(\GM)$. 
In general, since our semigroup is not always conservative, we must add an 
isolated point $\infty$ to $\GM$.

The distribution  $(\mathbb P_\xi )_{\xi\in \GM}$ of this process
on $\boldsymbol{\Omega}=\mathcal C ([0,\infty]\mapsto \GM\cup\{\infty\})$
is determined by the one-dimensional distributions
$$
\mathbb P_\xi(X_t\in U)=  \int_U h(t,\xi,\zeta) \, d\mu(\zeta)
=H_t\mathbf 1_U(\xi)
$$
for any open subset $U\subset \GM$, where $\xi$ is the starting point. 
The life time of the process is
$$
\tau_\infty
=\sup\{ t\ge 0: X_t\in \GM\}\,,
$$
and $H_t$ is conservative if and only if
$\mathbb P_\xi(\tau_\infty<\infty)=0$ for some (equivalently, all) 
$\xi\in \GM$.

For any relatively  compact  open set $U$, define the exit time
$$
\tau_U=\inf \{t>0: X_t\in U^c\}
$$
and, for $\xi\in U$,  the exit distribution
$$
\pi_U(\xi,B)=\mathbb E_\xi(X_{\tau_U}\in B).
$$
Since the process has continuous paths, for $\xi\in U$, the measure
$\pi_U(\xi,\cdot)$ is supported on the boundary $\partial U$ of $U$.
More generally, we set
$$
\pi_U(\xi,f)=\mathbb E_\xi\bigl(f(X_{\tau_U})\bigr)
$$
for any bounded Borel measurable function $f$ defined everywhere on 
$\partial U$.

The \emph{Green potential} of a continuous function $\varphi \ge 0$ with support in
$U$ can be written as
$$
G_{U}\varphi(\xi) 
= \mathbb{E}_{\xi}\left( \int_0^{\tau_{U}} \varphi(X_t) \,dt\right) \le +\infty
\,.
$$

\begin{dfn} A bounded Borel function $u$ in an open set $\Omega\subset \GM$
is \emph{$\mathbb P$-harmonic} (that is, harmonic
with respect to the process $X=(X_t)_{t\ge 0}$ with law $\mathbb P$)
if, for any open relatively compact set $B$ with $\overline{B}\subset \Omega$, 
we have
$$
\pi_B(\xi,u)=u(\xi) \quad \text{for all}\; \xi \in B\,.
$$
\end{dfn}

Since the associated semigroup $\{H_t :t>0 \}$ sends bounded measurable 
functions to bounded continuous functions it follows that any harmonic 
function is continuous; see e.g. {\sc Dynkin}~\cite[Vol. II]{Dynk}. The following 
result is important for our purpose.

\begin{thm}\label{weak-strong-harm}
 Let $\Omega\subset \GM$ be an open set.
\begin{itemize}
\item[(i)] If $u$ is a weak solution of $\Delta u=0$ in $\Omega$
then the continuous version of $u$ is $\mathbb P$-harmonic in $\Omega$.
\item[(ii)] If $u$ is $\mathbb P$-harmonic in $\Omega$
 then $u$ is a weak solution of $\Delta u=0$ in $\Omega$.
\end{itemize}
\end{thm}
\begin{proof} Part (i) is true in great generality, see 
\cite[Theorem 4.3.2]{FOT} (recall that, in our case, weak solutions are 
continuous). 

We now prove Part (ii). Without loss of generality, we can assume that $\Omega$
is relatively compact and $u\ge \epsilon>0$ in $\Omega$. 

Consider a fixed open set $V$ with $\overline V \subset \Omega$ (i.e.,
$V$ is relatively compact in $\Omega$). 

Let $\varphi$ be a non-negative continuous function (not identically $0$) 
with support in $U$, and let $w=G_{\Omega}\varphi\in \mathcal W^1_{0}(\Omega)$ be its 
Green potential in $\Omega$. 

Since $u$ is bounded from above in $U$ and the 
potential $w$ is bounded from below in $U$, there exists $t>0$
such that the excessive function $h=\min\{t\cdot w,u\}$ coincides with $u$ in
$U$, because $w|_{\partial \Omega} =0$. 
Moreover, $h$ coincides with $t\cdot w$ near the boundary of $\Omega$.
Since $h \leq t \cdot w$,  the function $h=G_{\Omega}\nu$ is the Green 
potential of a  measure $\nu$ with compact support in $\Omega$
and energy integral which is computed as
$$
\mathcal{E}(h,h)=\int_{\Omega} h\, d\nu < \infty\,.
$$
See {\sc Blumenthal and Getoor}~\cite[Ch. VI, Theorem 2.10]{BG} and 
{\sc Silverstein}~\cite[Ch. 1, Sec. 3]{Silv}.
In particular, $h\in \mathcal W^1_{0}(\Omega)$, and
since  $u$ coincides with $h$ in $V$, we see that $u$ is in
$\mathcal W^1_{\mbox{\tiny loc}}(\Omega)$.

Next, $u$ is represented inside any open set $V$ with 
$\overline V \subset \Omega$ 
as $u=\pi_V(\cdot,u)$.
Since $u$ is in $\mathcal W^1(V)$, the function $\pi_V(\cdot,u)$ coincides
with the Hilbert projection of $u$ on the linear subspace of weakly harmonic
functions in $V$. See \cite[Theorem 4.3.2]{FOT}.
\end{proof}

\subsection*{B. The bifurcation conditions}
The aim of this section is to prove that weak solutions of $\Delta u = 0$
are actually very regular in each strip and up to the bifurcation manifolds
although their various derivatives are typically not continuous across those
bifurcation manifolds. This will allow us to see that weak solutions verify
in a strong sense a particular bifurcation condition (or Kirchhoff's law)
along each bifurcation manifold. This bifurcation law is a crucial ingredient
in the analysis of our Dirichlet forms. It captures the influence of the
jumps of the functions $\phi$ and $\psi$ across bifurcation manifolds
and is crucial for an understanding of the domain of the infinitesimal 
generator.

Let us start by observing that, in any
open strip $S^o_e$, the infinitesimal generator $\Delta$ of our heat semigroup
is simply the weighted Riemannian Laplacian
$$
\Delta f= \frac{1}{\psi} \,\mbox{div}\bigl(\psi \, \mbox{grad}(f)\bigr),
$$
where $\mbox{ div }$ and $\mbox{ grad }$ refer, respectively, to the divergence
and gradient on the manifold
$$
\Bigl((e^-\,,\,e^+)\times M\,,\;\; \phi(s)\bigl((ds)^2+ g(\cdot,\cdot)\bigr)
\Bigr).
$$
More concretely, this means that for any $f$ in the domain of $\Delta$
and such that
$f\in \mathcal C^\infty(S^o_e)$,
$$
\Delta f 
= \frac{1}{\psi}\bigl[ \partial_s^2 +\Delta_M +\eta\, \partial_s\bigr]f\,,
\quad\text{where}\quad \eta=\partial_s \ln (\phi^{(n-1)/2} \,\psi).
$$
To be able to distinguish between the infinitesimal generator
and its expression in the interior of a strip, we make the
following definition.
\begin{dfn} \label{dfn-A}
For any $\xi\in \GM^o$ and any function $f$ which coincides with a smooth 
function in a neighbourhood of $\xi$, set
$$
\Lap f(\xi) = \frac{1}{\psi(\xi)}
\bigl[ \partial_s^2 +\Delta_M +\eta(\xi)\,\partial_s\bigr]f(\xi)\,,
\text{where}\quad \eta=\partial_s \ln (\phi^{(n-1)/2} \,\psi).
$$
In particular, $\Lap $ (as well as any of its integer powers $\Lap ^k$)
is a well defined continuous operator from
$\mathcal S^\infty(\GM^o)$ to ${\mathcal L}^\infty_{\mbox{\tiny loc}}(\GM)$.
\end{dfn}
In addition to the ``differential operator'' $\Lap $, there is another crucial
ingredient needed in order to describe harmonic functions on $\GM$ properly.
Namely, harmonic functions must satisfy a bifurcation condition
(or Kirchhoff law) along  each bifurcation manifold $\Ms_v\,$. To express
this bifurcation condition, we introduce the following notation.

\begin{dfn} 
Given $v\in V$ and $e\in E_v\,$, let $\nnu_{v,e}$ be the
outwards pointing normal unit vector relative to $S^o_e$ along $\Ms_v\,$.
\end{dfn}

We start by writing down  Green's formulas for a domain $\Omega$ with piecewise
smooth boundary contained in one strip $S_e$ and for smooth functions $f$, $h$
on $\Omega$. Then Green's formulas read as follows.

\begin{equation}\label{Green1}
\int_\Omega f\, \Lap h \,d\mu + 
\int_\Omega ( \nabla f,\nabla h )\, d\mu
=\int_{\partial \Omega}  (\nnu,\nabla h)\, f\,d\mu'
\end{equation}
and
\begin{equation}\label{Green2}
\int_\Omega (f \, \Lap h -  h \,\Lap f )\,d\mu
= \int_{\partial \Omega} \bigl((\nnu,\nabla h) \,f
          -(\nnu,\nabla f)\,h\bigr)\,d\mu'\,,
\end{equation}
where $\nnu$ is the outward unit normal vector to $\Omega$
and $\mu'$ is the induced measure on $\partial \Omega$. This measure has 
density $\psi_e(s)$ with respect to the Riemannian hypersurface measure on
$\Bigl(S_e\,,\,\phi\cdot\bigl((ds)^2+g(\cdot,\cdot)\bigr)\Bigr)$.

Let $u$ be a weak solution of $\Delta u = 0$
in a general domain $\Omega\subset \GM$ and let $U$
be a domain in a bifurcation manifold $\Ms=\Ms_v$ such that the closure of $U$
is contained in $\Omega$. 

Fix a strip $S=S_e^o$ attached to $\Ms_v\,$, and  
consider the outward unit normal derivative 
relative to $S^o_e$ along $\Ms_v\,$.

If $(U;x_1,\dots, x_n)$ is a local coordinate chart in $\Ms_v$ and
$(s,x_1,\dots,x_n)$ denotes the corresponding coordinate chart  in
$S^o_e=(e^-,e^+)\times U$ then that derivative is given by

\begin{equation}\label{nor-der}
(\nnu_{v,e}\,,\nabla) = \pm \phi_e(v)^{-1/2}\,\partial_s\,,\quad\text{if}\;
v=e^{\pm}.
\end{equation}
(The two signs have to coincide.) Note that it is crucial here to use the 
notation $\phi_e(v)$ since $\phi$ is not necessarily defined at $v$ and the 
values of the edge-wise extensions $\phi_e(v)$ of $\phi$ to the vertex $v$
may be distinct for different $e\in E_v\,$.

Suppose for the sake of simplicity that $v=e^-$.
Given $u$ as above, we want to define
$$
\delta = (\nnu_{v,e}\,,\nabla u)|_U = 
-\phi_e(v)^{-1/2}\,\partial_s u(v,\cdot)
$$
as a distribution on $U$. For $\ep>0$ (small enough), let 
$L_{\ep}= \{(s,x) \in S : s=s_{\ep}\}$ be the ``horizontal manifold''
in $S$ where $s_{\ep}$ is the point at distance $\ep$ from $e^-=v$
in the interval $I_e\,$. 
Let $U_{\ep}=\{(s,x): s=s_{\ep}\,, x\in U\}$. We assume that $\ep$ is so small 
that the closure of $U_{\ep}$ is contained in $\Omega$. For $0<\ep'<\ep$ fixed 
small enough we let $R_{\ep',\ep}$ be the rectangle with $U_{\ep'}$ and 
$U_{\ep}$  as horizontal sides.

Because $u$ is smooth inside the strip $S$, for any sufficiently 
small $\ep>0$ and any smooth function $\theta$ on $\Ms$ with compact support in 
$U$,
$$
\delta_{\ep}(\theta)= -\phi(s_{\ep})^{(n-1)/2}\,\psi(s_{\ep})
\int_{U} \partial_s u(s_{\ep},x) \, \theta(x) \,dx
$$
is well defined, and $\theta\mapsto \delta_{\ep}(\theta)$ is a distribution.

Now, we can compute $\delta_\epsilon(\phi)-\delta_{\epsilon'}(\phi)$ by setting
$\Theta(s,x)=\theta(x)$ and writing
$$
\int_{R_{\ep',\ep}} \bigl(\Theta \, \Lap   u - u \, \Lap \Theta\bigr) \,d\mu 
=\int _{\partial R_{\ep',\ep}}
\bigl((\nnu_{v,e}\,,\nabla u) \,\Theta
    -( \nnu_{v,e}\,,\nabla \Theta) \,u\bigr) \,d\mu' .
$$
Recall that $\theta$ is a smooth function with compact support in $U$.
It follows that $\Theta$ and $\nabla \Theta$ vanish on the vertical
components of $\partial R_{\epsilon',\epsilon}\,$. In addition, since $\Theta$
is independent of $s$, $\langle \nnu_{v,e}\,, \nabla \Theta\rangle$ vanishes on
the horizontal components of $\partial R_{\ep',\ep}\,$. Furthermore
$\Lap  u=0$ in $R_{\ep',\ep}$. Hence
$$
-\int_{R_{\ep',\ep}} u \, \Lap  \Theta \, d\mu 
= \int _{U_{\ep}}  (\nnu_{v,e}\,,\nabla u)\, \Theta\,  d\mu' 
 +\int _{U_{\ep'}} (\nnu_{v,e}\,,\nabla u)\, \Theta \,  d\mu',
$$
whence
$$
\left|\int_{R_{\ep',\ep}} u(s,x) 
\frac{1}{\phi(s)}\Delta_{\Ms} \theta (x)\, d\mu\,\right|
= |\delta_{\ep}(\theta)-\delta_{\ep'}(\theta)|.
$$
Since $u$ and $\Delta_{\Ms}\theta$ are uniformly bounded in a domain
containing all  rectangles $R_{\ep',\ep}$ with sufficiently small
$0<\ep'<\ep$, it follows that
$$
\lim_{\ep \rightarrow 0} \delta_{\ep}(\phi)=\delta(\phi)
$$
exists. If (as usually) $u_e$ denotes the restriction of $u$ to $S^o_e\,$, this
defines
$$
\delta= \bigl(\nnu_{v,e}\,, \nabla u_e(v,\cdot)\bigr)
$$
as a distribution on $U$.

In this way, we obtain $\deg(v)$ distributions $\delta_{v,e}\,$, one for each 
edge $e\in E_v\,$. Each $\delta_{v,e}$ corresponds to the unit outward normal 
derivative $\bigl(\nnu_{v,e}\,,\nabla u_{e}(v,\cdot)\bigr)$ in $S^o_e$ along 
$U\subset \Ms_v$. Now, the fact that $u$ is a weak solution of $\Delta u=0$
in $\Omega$ implies that 
\begin{equation}\label{dist-Kirch}
\sum_{e\in E_v} \psi_e(v)\, \delta_{v,e}=
\sum_{e\in E_v}\psi_e(v) \bigl(\nnu_{v,e}\,,\nabla u_e(v,\cdot)\bigr) =0
\mbox{ as distributions on } U\subset \Ms_v\,.
\end{equation}
We refer to this as the \emph{bifurcation condition} along $\Ms_v$ or
\emph{Kirchhoff's law,} in the sense of distributions.

For later purpose, it is useful to observe that the argument developed
above for weak solutions  of $\Delta u=0$
also works for weak solutions of the Poisson equation
$$
\Delta u= f
$$
in an open set $\Omega$ with a function $f$ that is H\"older continuous in  
$\Omega$. To be precise, we require here that $u\in \mathcal W^1(\Omega)$
and that for any
$h\in \mathcal W^1_0(\Omega)$,
$$
\mathcal E(u,h)=-\int f \,h \,d\mu\,.
$$
Note that by classical results, such a function $u$ has continuous
partial derivatives up to second order
and satisfies  $\Lap  u=f$ in the intersection of $\Omega$
with each open strip $S^o_e$. By an argument similar to the one used above
for weak solutions, the function $u$ must also satisfy the
bifurcation condition (\ref{dist-Kirch}) in the sense of distributions.

For instance, the function $u(\xi)=h(t,\zeta,\xi)$
is a weak solution of $\Delta u= f$ on $\HT$ with
$f(\xi)=\partial_t h(t,\zeta,\xi)$. Hence it satisfies
(\ref{dist-Kirch}) in the sense of distributions along each of
the bifurcation manifolds $\Ms_v$ in $\GM$.

\subsection*{C. Smoothness of harmonic functions}
The aim of this section is to show that weak solutions of
$\Delta u=0$ in an open set are smooth in the strip complex sense,
that is, they belong locally to $\mathcal C^\infty(\GM)$.
Since $\Delta $ is a non-degenerate elliptic operator
in each open strip, we know that harmonic functions are smooth there
(in the usual sense of having continuous partial derivatives of all orders).
The problem is to obtain smoothness up to the bifurcation manifolds
in each strip separately. Recall here that smoothness on $\GM$
does not imply continuity of the derivatives across
bifurcation manifolds.

\begin{thm} \label{th-Hsm}
Fix an open set $\Omega\subset \GM$.
For each $e\in E$, set $\Omega_e=\Omega\cap S^o_e\,$
and, if $u\in \mathcal C(\Omega)$, $u_e=u|_{\Omega_e}$.
A function $u$ is a weak solution of $\Delta u=0$ in $\Omega$ if and only if
it has  the following properties (more precisely, the continuous
version of $u$ has the following properties):
\begin{itemize}
\item $u\in \mathcal C^\infty(\Omega)\,$.
\item For any $e\in E$, one has $\;\Lap  u_e=0\;$ on $\Omega_e\,$. 
\item For any $v\in V$, one has 
$\;\sum\limits_{e\in E_v} \psi_e(v) \, (\nnu_{v,e}\,,\nabla u_e)=0\;$
along $\Ms_v\cap \Omega$.
\end{itemize}
\end{thm}

\begin{rmk} The first and third conditions are the crucial ones, since we
already know that the second condition  must hold by the local
ellipticity of our Laplacian in each open strip.
Concerning the first condition, we already know that
weak solutions are continuous (more precisely, have a continuous
representative) so the important part of the statement is that they
belong locally to $\mathcal S^\infty(\GM)$. We already observed
in \eqref{dist-Kirch} that the third condition must hold in the sense
of distributions but, if $u\in \mathcal C^\infty(\Omega)$, this is equivalent
to a classical pointwise statement as given by the theorem.
\end{rmk}

\begin{proof}[Proof of Theorem \ref{th-Hsm}] The proof goes through four steps 
and needs two auxiliary propositions.
\renewcommand{\qedsymbol}{}\end{proof}\begin{proof}
[Step 1: change of function] 
It will be useful to consider the functions
$$
w_e(\xi)= \beta_e(s) u(\xi)\,,\quad\text{where}\quad
\beta_e=\sqrt{\phi_e^{(n-1)/2}\psi_e} \AND  \xi=(s,x) \in I_e\times \Ms.
$$
Recall that $u$ satisfies 
$$
\Lap u=
\phi^{-1}\bigl[ \partial_s^2 +\Delta_{\Ms} + \eta \,\partial_s\bigr]u=0\,,
\quad\text{where}\quad\eta=\partial_s  \ln (\phi^{(n-1)/2} \,\psi)
$$
in each set $\Omega_e=\Omega\cap S_e^o$ and the bifurcation equation
$$
\sum_{e\in E_v} \psi_e(v)\,(\nnu_{v,e}\,, \nabla u_e)=0
$$
on each bifurcation manifold $\Ms_v\,$, where this is understood in the sense
of distributions. Observe that
$$
\frac{2\partial_s \beta_e}{\beta_e}
= \partial _s \ln (\phi_e^{(n-1)/2}\,\psi_e)=\eta_e\,.
$$
This implies that the functions $w_e\,$, $e\in E$, satisfy
$$
(\partial_s^2+\Delta_{\Ms})w_e = \phi_e \, \beta_e \,\Lap u_e
+ (\partial_s^2 \beta_e) u_e
= \frac{\partial_s^2 \beta_e}{\beta_e} \, w_e
$$
in each open strip $S^0_e$ and the bifurcation equation
$$
\sum_{e\in E_v} \psi_e(v)(\nnu_{v,e}\,, \nabla w_e)= 
-\left(\frac{1}{\phi_e(v)^{1/2}\, \beta_e(v)}
\sum_{e\in E_v} \ep_{v,e}\,\psi_e(v)\,
|\partial_s \beta_e (v)|\, w_e\right) \quad \text{along}\; \Ms_v\,,
$$
where
\begin{equation}\label{epsilon}
\epsilon_{v,e}= \begin{cases} 1\,, &\text{if}\; v=e^+\\
                             -1\,, &\text{if}\; v=e^-\,
	        \end{cases}
\end{equation}
\renewcommand{\qedsymbol}{}\end{proof}

\begin{proof}[Step 2: folding] 
As smoothness is a local property,
we can assume that $\Omega$ is a small neighbourhood
of a point $\xi_0=(v,x_0)$ on a fixed bifurcation manifold  $\Ms_v$ and
that $\Omega_e=\Omega\cap S_e^o$ is of the form $(v\,,\,r_e)\times U$ where
$r_e \in I_e^o = (e^-\,,\,e^+)$, and all intervals $(v\,,\,r_e)$ in $\Gamma^1$ 
have the same (small) length $l$. This provides us with an obvious way to 
identify all the different $\Omega_e$ with a  fixed set
$$
\Omega_+=(0\,,\,l)\times U \subset (0\,,\,\infty)\times \Ms.
$$
Using this identification, we can consider each  $w_e$
as a function defined on $\Omega_+$, namely,
$$
\Omega_+\ni (s,x)\mapsto w_e\bigl(s_{(v,e)},x\bigr)
$$
where $s_{(v,e)}$ is the point on $I_e\,$, $e \in E_v\,$, at distance
$s$ from $v$.
Now Theorem \ref{th-Hsm} will be an immediate consequence of the next 
result. 
\renewcommand{\qedsymbol}{}\end{proof}

In the following proposition, $E_v$ can be viewed as an arbitrary finite 
set of parameters whose elements are denoted by $e$.

\begin{pro} \label{pro-Hsm}
Let $U$ be a relatively compact domain in $\Ms$.
Let
$$
\Omega_+=(0\,,\,l)\times U \subset (0\,,\,\infty)\times \Ms
$$
and $I=\{0\}\times U$ be the bottom of $\Omega_+\,$. 
For all $e, e'  \in E_v\,$, let $\delta_e>0$, $\tilde\delta_e\in \R$
and $c_{e,e'}>0$ be fixed numbers.
Let $\gamma_e$, $e\in E_v\,$, be functions in $\mathcal C^\infty([0\,,\,l])$.

Assume that $w_e$, $e\in E_v$, are functions defined on $\Omega_+$
that belong to $\mathcal C^\infty(\Omega_+)$ and satisfy the following
hypotheses.
\begin{itemize}
\item  For each $e\in E_v\,$, the function $w_e$ is in 
$\mathcal C^\alpha(\overline{\Omega_+})$ for some $\alpha \in (0\,,\,1)$, and 
$$
w_e\big|_{I} = c_{e,e'} \, w_{e'}\big|_{I}
\quad\text{for all}\; \,e,e'\in E_v\,,
$$
\item  $\; [\partial_s^2+\Delta_{\Ms}] w_e=  \gamma_e\, w_e\;$ in $\Omega_+\,$.
\item  The partial derivatives $\partial_s w_e(0,\cdot)$, $e\in E_v\,$,
whose existence in the sense of distributions in $U$ is guaranteed by
the first two hypotheses, satisfy
$$
\sum _e \delta_e \, \partial_s w_e(0,\cdot)
=  \sum_e \tilde\delta_e\,w_e(0,\cdot)
$$
in the sense of distributions in $U$.
\end{itemize}
Then $w_e\in \mathcal C^\infty ([0\,,\,l)\times U)$ for each $e\in E_v\,$,,
i.e., it is smooth up to the bottom $I$ of $\Omega_+\,$.
\end{pro}

In order to prove this proposition, set
$$
W(s,x)=  \sum_{e\in E_v} \delta_e \,w_e(s,x)\,.
$$
The function $W$ is continuous in 
$\overline{\Omega_+}=[0\,,\,l]\times \overline{U}\,$.
Moreover, it satisfies
\begin{equation}\label{sysW}
\left\{\begin{array}{l}[\partial_s^2+\Delta_M]W = W_1\;
\mbox{ in }\; \Omega_+\,,\\[3pt]
\partial_s W(0,\cdot)= W_2\;  \mbox{ on }\; U\,,
\end{array}\right.
\end{equation}
where
\begin{equation} \label{sysW'}\left\{ \begin{array}{l}\displaystyle
W_1= \sum_{e} \delta_e \,\gamma_e \, w_e \in 
 \mathcal C^\alpha(\overline{\Omega_+}) \AND
\\[3pt]
W_2=\frac{1}{\delta} \sum_e \tilde\delta_e \, w_e(0,\cdot)
\in \mathcal C^\alpha(\overline{U})\,,\quad\text{with}\;
\delta=\sum_{e\in E_v} \delta_e.\end{array}\right.
\end{equation}

At this point, the proof of Proposition \ref{pro-Hsm} requires another
auxiliary result, as follows.

\begin{proof}[Step 3: improved regularity]
\renewcommand{\qedsymbol}{}\end{proof}

\begin{pro}\label{pro-Ck1} With notation as in Proposition \ref{pro-Hsm},
fix $\alpha\in (0\,,\,1)$ and a nonnegative integer $k$.
Also fix $h_1\in \mathcal C^{k+\alpha}(\overline{\Omega_+})$ and
$h_2\in \mathcal C^{k+\alpha}(\overline{U})$.
Let $f$ be a smooth function in $\Omega_+$  which belongs to
$\mathcal C^{k+\alpha}(\overline{\Omega_+})$ and
satisfies
$$
\left\{\begin{array}{l}[\partial_s^2+\Delta_{\Ms}] f= h_1 \; \mbox{ in }\;\Omega_+
\\[3pt]
\partial_s f = h_2 \;\mbox{ in }\; I\;
\mbox{ (in the sense of distributions when $k=0$)}.\end{array}\right.
$$
Then $f$ belongs to $\mathcal C^{k+1+\alpha}(\overline{\Omega'_+})$ for
every set $\Omega'_+= (0,l')\times U'$, where $0<l'<l$ and $U'$ is open 
and relatively compact in $U$.
\end{pro}

\begin{proof}[Proof of the proposition] 
Without loss of generality (because of well-known basic extension theorems,
see e.g. {\sc Seeley}~\cite{See}),
we can assume that $h_1=h|_{\Omega_+}$ is the restriction to $\Omega_+$
of a function $h \in \mathcal C^{k+\alpha}(\R \times \Ms)$ with compact 
support. Let $B$ be a ball in $\R \times M$ containing the support of $h$.
Let $H=G_Bh$ be the Green potential of $h$ relative to this ball $B$ and with
respect to the operator $\partial_s^2+\Delta_{\Ms}$. Then  
$H\in \mathcal C^{k+2+\alpha}_{\mbox{\tiny loc}}(B)$, and within $\Omega_+$ 
we have
$$
[\partial_s^2+\Delta_M](f+H)=0.
$$
Obviously, on the boundary $I$, the function $f+H$ satisfies
$$
\partial_s (f+H)\big|_{I}=h_2+ \partial_sH\big|_{I}\,.
$$
Note that $f+H\in \mathcal C^{k+\alpha}(\overline{\Omega_+})$
and $h_2+\partial _s H\big|_{I}\in 
\mathcal C^{k+\alpha}(\overline{I})$.
Thus, replacing $f$ by $f+H$, we are led to study the solutions
$f\in \mathcal C^{k+\alpha}(\overline{\Omega_+})$ of
$$
\left\{\begin{array}{l}[\partial_s^2+\Delta_{\Ms}] f=0\;
\mbox{ in }\;\Omega_+\\[3pt]
\partial_s f = h\; \mbox{ on }\; I\,,\end{array}\right.
$$
where $h\in \mathcal C^{k+\alpha}(\overline{U})$.
Indeed, to prove Proposition \ref{pro-Ck1}, it suffices to show  that such
$f$ must be in $\mathcal C^{k+1+\alpha} (\overline{\Omega'_+})$.

Recall that $I$ is the bottom of $\Omega_+\,$. Let 
$\Omega'_0 = \{0\} \times U'$. Identifying $\{0\} \times\Ms$ with $\Ms$,
there exists a function $f_1\in \mathcal C_c^{k+\alpha}(\Ms)$ and a set 
$J$ open in $I$ with $\overline {I'} \subset J \subset I$ such that 
$f|_J = f_1\big|_J\,$. Thus, if we decompose
$f|_I=f_1+f_2$ on $I$, then $f_2 = 0$ on $J$.

Let $(s,x)\mapsto F_1(s,x)$ be the harmonic function on $(0,\infty)\times \Ms$  
which coincides with $f_1$ on $\{0\} \times\Ms$, that is, the
Poisson integral given formally by $F_1(s,x)=e^{-s\sqrt{-\Delta_{\Ms}}}f_1(x)$.
Then, in $\Omega_+$, we have $f=F_1+F_2$ where $F_2$ is harmonic in
$\Omega_+$ with boundary values $0$ on $J$. In particular, $F_2$ has bounded
continuous derivatives of all orders up to  $J$.  Moreover, along $J$
we have in the sense of distributions on $J$
$$ 
h=\partial_s f(0,\cdot) =  -\sqrt{-\Delta_{\Ms}} \,f_1 + \partial_s F_2(0,\cdot).
$$
Write this as
$$
\bigl[\Id+\sqrt{-\Delta_{\Ms}}\,\bigr] f_1\big|_J  =
\bigl(- h + f_1 +\partial_sF_2(0,\cdot)\bigr)\big|_J\,, 
$$
again in the sense of distributions. (Here, $\Id$ is the identity operator.)
By hypothesis, the  right-hand side is in
$\mathcal C^{k+\alpha}_{\mbox{\tiny loc}}(J)$.
Let  $f_3\in \mathcal C^{k+\alpha}_c(J)$ be a function which
coincides with $\bigl(-h +f_1 +\partial_sF_2(0,\cdot)\bigr)\big|_J$
in a neighbourhood $J'$ of $\overline{I'}$ that is contained in $J$.
Let $f_4=\bigl[\Id+\sqrt{-\Delta_{\Ms}}\,\bigr]^{-1}f_3\,$. Then
$f_4\in \mathcal C^{k+1+\alpha}_{\mbox{\tiny loc}}({\Ms})\cap {\mathcal L}^2({\Ms})$, 
and the function  $f_1-f_4\in {\mathcal L}^2({\Ms})$ satisfies
$$
\bigl[\Id+\sqrt{-\Delta_{\Ms}}\,\bigr](f_1-f_4) =0 \;\;\mbox{ in }\; J'.
$$
In addition, the distribution $\bigl[\Id+\sqrt{-\Delta_{\Ms}}\,\bigr](f_1-f_4)=
\bigl(\bigl[\Id+\sqrt{-\Delta_{\Ms}}\,\bigr]f_1\bigr)-f_3$ can be represented by 
a function in ${\mathcal L}^2(\Ms)$ outside $I$ because $f_1$ is continuous with compact 
support in $I$. By the hypoellipticity of
$\bigl[\Id+\sqrt{-\Delta_{\Ms}}\,\bigr]$ (see Theorem \ref{th-hypsqrt} in the
Appendix) it follows that $f_1-f_4$ is in 
$\mathcal C^\infty_{\mbox{\tiny loc}}(J')$.
Hence $f_1$ is $\mathcal C^{k+1+\alpha} _{\mbox{\tiny loc}}(J')\,$:
it has the same smoothness as $f_4$ in $J'$. This implies
that the Poisson integral $F_1$ of $f_1$ is in 
$\mathcal C^{k+1+\alpha}(\overline {\Omega'_+})$. Hence $f=F_1+F_2$ is in
$\mathcal C^{k+1+\alpha}(\overline{\Omega'_+})$. This is the desired result.
\end{proof}

\begin{proof}[Step 4: final bootstrap] We now prove Proposition
\ref{pro-Hsm} by induction on the smoothness parameter $k$, using Proposition
\ref{pro-Ck1}. Assume we have proved that the functions $w_e$ in Proposition
\ref{pro-Hsm} are in $\mathcal C^{k+\alpha}(\overline{\Omega'_+})$ for some
integer $k$ and any $\Omega'=(0\,,\,l')\times U'$ relatively compact in 
$\Omega_+\,$. This implies that the functions $W_1, W_2$ of \eqref{sysW'}
are respectively in $\mathcal C^{k+\alpha}(\overline{\Omega'_+})$
and $\mathcal C^{k+\alpha}(\overline{U})$.  Hence we can apply
Proposition \ref{pro-Ck1} to the function $W$ of \eqref{sysW}. This gives
that $W\in \mathcal C^{k+1+\alpha}(\overline{\Omega^*_+})$
where $\Omega^*_+=(0\,,\,l^*)\times U^*$ with $l^*$ an arbitrary real in 
$(0\,,\,l')$ and $U^*$ an  arbitrary open relatively compact set in $U'$. 
Because $l'\in (0\,,\,l)$ and $U'$, relatively compact in $U$, are arbitrary,
we conclude that $W\in \mathcal C^{k+1+\alpha}(\overline{\Omega'_+})$ for any
$\Omega'=(0\,,\,l')\times U'$ relatively compact in $\Omega_+$.

The functions $w_e\,$, $e\in E_v\,$, are related on $\{0\}\times U$
by $$
w_e(0,x)= c_{e,e'}\,w_{e'}(0,x)
$$
and  thus are all equal on $\{0\}\times U$  to a fixed multiple of
$W(0,\cdot)\in \mathcal C^{k+1+\alpha}(U)$.
Each of the functions $w_e$ is solution of
$$
\left\{\begin{array}{l}\bigl[\partial_s^2+\Delta_{\Ms}\big] f=h_{e,1}\;
\mbox{ in }\; \Omega_+\\[3pt]
f(0,\cdot) = h_{e,2} \;\mbox{ on }\; U\,,\end{array}\right.
$$
where $h_{e,1}=\gamma_e\, w_e\in \mathcal C^{k+\alpha}(\overline{\Omega_+})$
and $h_{e,2} = w_e(0,\cdot)\in \mathcal C^{k+1+\alpha}(U)$.

Let $H_{e,1}$ be the Green potential (in a large ball in $\R\times M$)
of a compactly supported extension of $h_{e,1}$ that belongs to
$\mathcal C^{k+\alpha}(\R\times {\Ms})$. The function
$H_{e,1}$ is $\mathcal C^{k+2+\alpha}(\overline{\Omega_+})$,
and $w_e-H_{e,1}$ is solution of
$$
\left\{\begin{array}{l}\bigl[\partial_s^2+\Delta_M\bigr] f=0
\;\mbox{ in }\;\Omega_+\\[3pt]
f(0,\cdot)= h_{e,2} -H_{e,1}(0,\cdot) \;\mbox{ on }\; U\,,\end{array}\right.
$$
where $h_{e,2}-H_{e,1}(0,\cdot)\in \mathcal C^{k+1+\alpha}(U)$. It follows 
that $w_e-H_{e,1}$ is in $\mathcal C^{k+1+\alpha}(\overline{\Omega'_+})$.
This means that each of the functions $w_e$ is in
$\mathcal C^{k+1+\alpha}(\overline{\Omega'_+})$.
\end{proof}

Given an open connected set $\Omega$, consider the linear space $H(\Omega)$
of  all weak solutions of the Laplace equation $\Delta u=0$ in $\Omega$.
By the local H\"older regularity result and the fact that the notion of weak 
solutions and of $\mathbb P$-harmonic functions coincide, it follows that
$H(\Omega)$ equipped with the seminorms of the uniform convergence on compact
subsets of $\Omega$ is a complete seminormed vector space.

By Theorem \ref{th-Hsm}, any element $u$ of $H(\Omega)$ is in
$\mathcal C^\infty(\Omega)$. The closed graph theorem then
yields the following result.

\begin{cor} Let $\Omega$ be an open connected set in $\GM$ and
$\Omega_0$ relatively compact in $\Omega$.
Let $I\times U$ be a relatively compact coordinate chart in $\GM$
such that $K=\overline{I\times U}\subset \Omega_0\,$.
Fix $\kappa=(\kappa_0,\kappa_1,\dots,\kappa_n)$.
Then there exists a constant $C=C(\Omega _0,K,\kappa)$ such that
$$
\sup_{\xi\in K}| \partial_\xi^\kappa u(\xi)|
\le C \sup_{\Omega_0} |u| \quad\text{for all}\;u\in H(\Omega)\,.
$$
\end{cor}

\subsection*{D. Regularity of certain weak solutions of the heat equation}
Let $(t,\xi)\mapsto u(t,\xi)$ be a weak solution of the heat equation in
$(0\,,\,T)\times\Omega$, where $\Omega$ is an open set in $\GM$.
We already know that we can regard $u$ as a H\"older continuous function
on $(0\,,\,T)\times\Omega$. Our aim is to show that in some cases, including 
the case of the heat kernel, that $u(t,\cdot)\in \mathcal C^\infty(\Omega)$  
for each $t\in (0\,,\,T)$, and moreover, for any positive integer
$k$, $\partial_t^ku(t,\cdot)\in \mathcal C^\infty(\Omega)$.
(See Definition \ref{smoothOmega} for the definition of 
$\mathcal C^\infty(\Omega)$.) It is plausible that this result holds for 
any weak solution, but our proof below does not provide this stronger 
result. 

\begin{dfn} Fix  $k\in \{0,1,\dots,\infty\}$, $T>0$, $I=(0\,,\,T)$
and an open set $\Omega \subset \GM$. See \eqref{Floc} for the defintion of 
$\mathcal F_{\mbox{\tiny loc}}(I\times \Omega)$. We say that a weak solution  
$u \in \mathcal F_{\mbox{\tiny loc}}(I\times\Omega)$
of the heat equation in  $I\times \Omega$ is time regular to order $k$ if, 
for each $m\in\{0,1,\dots,k\}$, the distributional time derivative 
$\partial_t^m u$ exists and can be represented by a function 
$u_m\in \mathcal F_{\mbox{\tiny loc}}(I\times \Omega)$ which is a weak 
solution of the heat equation in $I\times \Omega$. 
When $u$ is time regular to infinite order we simply say that $u$ is a 
time regular weak solution in $I\times \Omega$.
\end{dfn}

\begin{exa} Fix $f\in {\mathcal L}^2(\GM)$. Then $u(t,\xi)=H_tf(\xi)=e^{t\Delta}f(\xi)$
is a time regular weak solution up to infinite order in 
$(0\,,\,\infty)\times \GM$. Fix $\zeta\in \GM$ and set 
$u(t,\xi)=h(t,\xi,\zeta)$. Then $u$ is again a time regular solution
to infinite order in $(0\,,\,\infty)\times \GM$.  Fix an open set 
$\Omega \subset \GM$ and consider the Dirichlet Laplacian $\Delta_\Omega$ in 
$\Omega$. This is the infinitesimal generator associated with the closure of the
form $\bigl(\int_\Omega |\nabla f|^2d\mu, \mathcal C_c^\infty(\Omega)\bigr)$.
Let $f\in {\mathcal L}^2(\Omega)$ and consider
$u(t,\xi)=e^{t\Delta_\Omega}f(\xi)$, $(t,\xi)\in (0,\infty)\times \Omega$.
This is a time regular weak solution up to infinite order in
$(0\,,\,\infty)\times \Omega$ and so is the corresponding Dirichlet heat
kernel in $\Omega$.
\end{exa}

\begin{thm} \label{th-Heatsm}
Fix $T>0$ and an open set $\Omega\subset \GM$.
For each $e\in E$, set $\Omega_e=\Omega\cap S^o_e$ and, if 
$u\in \mathcal C\bigl((0\,,\,,T)\times \Omega\bigr)$, set 
$u_e=u|_{(0,T)\times\Omega_e}$.
Any function $u$ which is a weak solution of $[\partial_t-\Delta] u=0$ in
$Q=(0\,,\,T)\times \Omega$ and is time regular to order $k$
has  the following properties:
\begin{itemize}
\item For any $m=0,1,2,\dots,k$,  the derivative $\partial^m_tu $ is a 
continuous function on $(0\,,\,T)\times \Omega$. Moreover, there is 
$\alpha \in (0\,,\,1)$ such that
$\partial_t^m u(t,\cdot)\in \mathcal C^{k-m+\alpha}(\Omega)$ for any 
$t\in (0\,,\,T)$.
\item For any $e\in E$, one has $[\partial_t-\Lap ] u_e=0$ on 
$(0\,,\,T)\times \Omega_e\,$. In particular, $u_e$ is smooth (in the usual 
sense) in the open set  $\Omega_e\,$.
\item For any $m\in \{0,1,\dots,k-1\}$ and
$v\in V$, 
$$
\sum_{e\in E_v} \psi_e(v) \, (\nnu_{v,e}\,,\nabla \partial^m_tu_e)=0\quad
\mbox{along }\; (0\,,\,T)\times (\Ms_v\cap \Omega).
$$
\end{itemize}
\end{thm}

\begin{proof} The proof goes through three steps and involves Proposition 
\ref{pro-Heatsm} 
below.
\renewcommand{\qedsymbol}{}\end{proof}

\begin{proof}[Step 1: change of function] As in the elliptic case,
we consider the functions
$$
w_e(t,\xi) = \beta_e(s) \, u_e(t,\xi),\quad\text{where}\quad 
\beta=\sqrt{\phi^{(n-1)/2}\,\psi}  \AND \xi=(s,x) \in I_e \times \Ms.
$$
Recall that $u$ satisfies 
$$
\Lap u=
\frac{1}{\phi}\bigl[ \partial_s^2 +\Delta_{\Ms} +\eta \partial_s\bigr]u
=\partial_t u\,,
\quad\text{where}\quad \eta=\partial_s \ln (\phi^{(n-1)/2} \,\psi)
$$
in each set $\Omega_e=\Omega\cap S_e^o$ and the bifurcation equation
$$
\sum_{e\in E_v} \psi_e(v)\,(\nnu_{v,e}, \nabla u_e)=0
$$
on each bifurcation manifold $M_v\,$, where this is understood in the sense
of distributions. As in the proof of Theorem \ref{th-Hsm},
this implies that the functions $w_e$ satisfy
$$
[\partial_s^2+\Delta_{\Ms}]w_e =
\frac{\partial_s^2 \beta_e}{\beta_e}\, w_e + \phi_e\, \partial_t w_e
$$
in each open strip $S^0_e$ and the bifurcation equation
$$
\sum_{e\in E_v} \psi_e(v)(\nnu_{v,e}\,, \nabla w_e)= 
-\left(\frac{1}{\phi_e(v)^{1/2}\, \beta_e(v)}
\sum_{e\in E_v} \epsilon_{v,e} \,\psi_e(v)\,
|\partial_s \beta_e (v)|\, w_e \right) \quad\text{along}\; M_v\,,
$$
where $\epsilon_{v,e}$ is as in \eqref{epsilon}.
\renewcommand{\qedsymbol}{}\end{proof}

\begin{proof}[Step 2: folding and improved regularity]
The following is analogous to Proposition \ref{pro-Hsm} except for
the role played by the function $\widetilde{w}_e$.
\begin{pro} \label{pro-Heatsm}
Let $U$ be a relatively compact domain in $M$. Let
$$
\Omega_+=(0\,,\,l)\times U \subset (0\,,\,\infty)\times \Ms
$$
and $I=\{0\}\times U$ be the bottom of $\Omega_+\,$. 
For all $e, e'  \in E_v\,$, let $\delta_e>0$, $\tilde\delta_e\in \R$
and $c_{e,e'}>0$ be fixed numbers.
Assume that $w_e, \widetilde{w}_e$, $e\in E_v\,$, are functions defined on 
$\Omega_+$ that belong to $\mathcal C^\infty(\Omega_+)$ and satisfy the 
following hypotheses.
\begin{itemize}
\item  For each $e\in E_v\,$, the functions $w_e,\widetilde{w}_e\,$ are in
$\mathcal C^{k+\alpha}(\overline{\Omega_+})$ for some integer $k$ and
$\alpha \in (0\,,\,1)$, and
$$ 
w_e\big|_{I} = c_{e,e'}\, w_{e'}\big|_{I} \in 
\mathcal C^{k+\alpha}(\overline{U})\quad\text{for all}\; \,e,e'\in E_v\,,
$$
\item  $[\partial_s^2+\Delta_{\Ms}] w_e = \widetilde{w}_e$ in $\Omega_+\,$.
\item  The partial derivatives $\partial_s w_e(0,\cdot)$, $e\in E_v\,$,
whose existence in the sense of distributions in $U$ is guaranteed by the 
first two hypotheses, satisfy
$$
\sum_e \delta_e\, \partial_s w_e(0,\cdot)
= \sum_e\tilde \delta_e\, w_e(0,\cdot)
$$
in the sense of distributions in $U$.
\end{itemize}
Then $w_e\in \mathcal C^{k+1+\alpha}([0\,,\,l)\times U)$ for each $e\in E_v\,$.
\end{pro}

The proof of this result follows exactly the same line as the proof
of Proposition \ref{pro-Hsm}, except for the very last step (bootstrap)
that cannot be performed in the present case because of the presence
of the functions $\widetilde{w}_e$ on the right-hand side of the
second condition.
This is why we only obtain improved smoothness from $\mathcal C^{k+\alpha}$
to $\mathcal C^{k+1+\alpha}$. 
\renewcommand{\qedsymbol}{}\end{proof}

\begin{proof}[Step 3: finite order bootstrap]
When applying  Proposition \ref{pro-Heatsm} to weak solutions of the heat 
equation, the function $\widetilde{w}_e$ has the form
$$\widetilde{w}_e= 
\frac{\partial_s^2\beta_e}{\beta_e} \,w_e+ \phi_e \, \partial_t w_e\,.
$$
In order to apply  Proposition \ref{pro-Heatsm} repeatedly,
we need to improve not only the smoothness of $w_e$ but also the
smoothness of $\partial_tw_e\,$. For instance, in order to apply Proposition 
\ref{pro-Heatsm} and obtain $\mathcal C^{1+\alpha}$-regularity of $w_e\,$, 
we need first to prove that $\partial_tw_e$ is H\"older continuous.
Observe that this property immediately follows if we know that
the original weak solution $u_e$ of the heat equation is such that
that $\partial_tu_e $ is also a weak solution of the heat equation.

Assume now that $u$ and all its time derivatives $\partial_t^m$
up to order $k$ are weak solutions of the heat equation in
$(0\,,\,T)\times \Omega$. Then all the partial derivatives $\partial^m_t u$,
$m\in \{0,\dots, k\}$ are H\"older continuous
and we can apply Proposition \ref{pro-Heatsm}
simultaneously to all the functions $\partial^m_tw_e$, where $e\in E_v$ and
$m\in \{0,1,2,k-1\}$, to conclude that these functions are in
$\mathcal C^{1+\alpha}$. Using this conclusion, and applying
Proposition \ref{pro-Heatsm} to $\partial^m_tw_e$, where $e\in E_v$ and
$m\in \{0,1,2,k-2\}$, we conclude that these functions are in
$\mathcal C^{2+\alpha}$.
Proceeding by finite induction, Theorem \ref{th-Heatsm} follows.
\end{proof}

\begin{dfn}
Fix $T>0$ and an open set  $\Omega\subset \GM$ and set 
$Q=(0\,,\,T)\times\Omega$. Let $R_k(Q)$ be the vector space of all weak 
solutions in $(0,T)\times\Omega$  that are time regular
to order $k$ in $(0\,,\,T)\times \Omega$, equipped with the seminorms
\begin{eqnarray*}
N_{k,Q'}(u)
&=&\sup_{Q'}\sup_{m\in \{0,\dots,k\}}\bigl|\partial^m_tu(t,\xi)\bigr|\\
&&+ \sup_{v\in \mathcal F_c(Q')}
       \left|\int_{Q'} v\,\partial^{k+1}_t\,u\, d\mu\, dt\right|
+ \sup_{v\in \mathcal F_c(Q')}
\left|\int_{Q'} (\nabla v,\nabla \partial_t^ku) \,d\mu \,dt\right|\,,
\end{eqnarray*}
where $Q'=I'\times \Omega'$ is relatively compact in $(0\,,\,T)\times \Omega\,$.
\end{dfn}
The first term in the seminorm $N_{Q'}$ controls the sup-norms (hence
the ${\mathcal L}^2$-norms) in $Q'$ of the time derivatives up to order $k$. Since these
functions are weak solutions, this yields a control of the ${\mathcal L}^2$-norms
of $|\nabla \partial^m_t u|$ for $m$ up to $k-1$. The last two terms
provide the additional control needed to insure that the seminormed
space $R_k(Q)$ is complete (a limit in this topology
of a sequence of weak solutions that are all time
regular up to order $k$ is, itself, such a solution).

\begin{cor} Let $T>0$, $(a\,,\,b)$ a relatively compact interval in
$(0\,,\,T)$ and $[a'\,,\,b']$ be a compact interval in $(a\,,\,b)$.
Let $\Omega$ be an open connected set in $\GM$ and
$\Omega'$ be a subset that is relatively compact in $\Omega$.
Set $Q=(0\,,\,T)\times \Omega$, $Q'=I'\times\Omega'$.
Let $I\times U$ be a relatively compact coordinate chart in $\GM$
such that
$K=\overline{I\times U}\subset \Omega'$.
Fix integers $k$, $\kappa_*$,
$\kappa=(\kappa_0,\kappa_1,\dots,\kappa_n)$ with
$\kappa_*+\sum_0^n\kappa_i\le k$.
Then there exists a constant $C=C(a,a',b,b',\Omega',K,k)$
such that if $u\in R_k(Q)$ is a weak solution of the heat
equation in $Q$, time regular to  order $k$, then we have
$$
\sup
\bigl\{\bigl|\partial_t^{\kappa_*}\partial _\xi^{\kappa}u(t,\xi)\bigr|:
(t,\xi)\in [a'\,,\,b']\times K \bigr\} \le C \,N_{k,Q'}(u).
$$
\end{cor}

Applying this to the heat kernel which is a time regular weak solution
to infinite order, we obtain the following important result.

\begin{thm} \label{th-HKsmth}
For any fixed $\zeta\in \GM$, and integer $k$,
the function  $\xi \mapsto \partial_t^kh(t,\xi,\zeta)$ is in
$\mathcal C^\infty(\GM)$.
\begin{itemize}
\item Fix a relatively compact coordinate chart
$I\times U $ and $\kappa=(\kappa_0,\kappa_1,\dots,\kappa_n)$. Then, for fixed 
$\xi\in I\times U$, the function
$$
(t,\zeta) \mapsto u(t,\zeta) =\partial_t^k\partial^\kappa_\xi h(t,\xi,\zeta)
$$ 
is in $\mathcal C^\infty(\GM)$. It is a weak solution of the heat equation,
and it satisfies the bifurcation condition
$$
\sum_{e\in E_v} \psi_e(v)\,(\nnu_{v,e}\,,\nabla u)=0
$$
(in the classical sense) along each bifurcation manifold $M_v\,$, $v\in V$.
\item Fix a compact time interval $[a\,,\,b]\subset (0\,,\,\infty)$ and a 
relatively compact coordinate chart $I\times U$ in $\GM$ with 
$\xi_0\in I\times U$. Fix also integers $k$ and $\kappa_0, \dots, \kappa_n$ and
set $\kappa=(\kappa_0,\dots,\kappa_n)$. Then there exists a constant
$C=C(a,b,I,U,k,\kappa)$ such that
$$
\sup \bigl\{ \bigl| \partial ^k_t\partial^\kappa_\xi h(t,\xi,\zeta)\bigr|:
(t,\xi)\in [a\,,\,b]\times I\times U \bigr\}
\le C \, h(2b,\xi_0,\zeta)
\quad \text{for all}\;\zeta\in \GM\,.
$$
\end{itemize}
\end{thm}

\section{Projections}\label{projs}

Recall the following simple version of transformation of phase space.
See \cite[Vol. II, Thm. 10.13]{Dynk}.

Let $X$ be a separable metrisable space equipped with a Radon measure 
$\mu$ with full support and with a symmetric Markov semigroup $\{ H_t : t>0\}$
of operators on ${\mathcal L}^2(X)={\mathcal L}^2(X,\mu)$. Denote also by 
$H_t$ the extension of that
operator from ${\mathcal L}^2(X)\cap {\mathcal L}^\infty(X)$ to 
${\mathcal L}^\infty(X)$. Assume that $(H_t)$ admits a transition function 
$h(t,x,\dot)$,
that is, for any $f\in \mathcal L^\infty (X)$ and for all $t>0$ we have
$H_t f(x)= \int_X f(y)\, h(t,x,dy)$ for $\mu$-almost every $x$.
Let $\bigl((X_t)_{t\ge 0}\,,\mathbb P_x\bigr)$ be the associated Markov process.
In the applications of interest to us here, $X=\GM$ and the process
is the one associated with our Dirichlet form.

Let $G$ be a locally compact group  acting properly and continuously
on $X$, and let $\uX$ be the topological quotient space and 
$\pi: X \to \uX$ the quotient map.
Assume that $H_t$ commutes with the action of $G$, that is,
$[H_tf](gx)= H_t f_g(x)$ for all bounded measurable functions $f$ on $X$,
where $f_g(x)=f(gx)$. Then $H_t$ induces  a semigroup of contractions
$\uH_t: {\mathcal L}^\infty(\uX)\rightarrow {\mathcal L}^\infty(\uX)$ defined by
$$
\uH_t f(\ux)= [H_t f\circ \pi](x)\,,\quad\text{where}\;\; \ux=Gx.
$$
Moreover, the formula $\uX_t=\pi(X_t)$, $t>0$,  defines a Markov process
on $\uX$ with law $\mathbb P_{\ux}$ satisfying
$\mathbb P_{\underline{x}}(\uX_t\in A)= \uH_t\mathbf 1_A(\ux)=
\mathbb P_x[X_t\in \pi^{-1}(A)]$, where $\pi(x)=\ux$.
Note that in general there is no obvious natural way to project the 
${\mathcal L}^2$-structure onto  $\uX$. In particular, in this abstract setting and 
unless either $X$ or $G$ is compact, there is \emph{a priori} no natural 
reference measure on $\uX$.

For the purpose of the next theorem, we say that a semigroup $\{P_t : t > 0\}$
defined on ${\mathcal L}^\infty(X)$ is a \emph{Markov semigroup} if it admits a 
transition function $p_t(x,f)$ as defined in \cite[Vol. I, Ch. 2]{Dynk}.
By \cite[Vol. I, Thm. 2.1]{Dynk}, this is equivalent to say that
$\{P_t : t > 0\}$ can be viewed as a semigroup of contractions on the
space $\mathcal B(X)$ of all bounded measurable functions on $X$ 
(not classes of functions!) 
that preserves positivity and such that $P_0f(x_0)=0$ if $f(x_0)=0$.
As for any $t > 0$ and $x \in X$, $p_t(x,\cdot)$ is a Borel measure on $X$,
the action of $P_t$ on ${\mathcal L}^\infty(X)$ is determined by 
its action on $\mathcal C_c(X)$.

 \begin{thm} \label{th-proj}
Let $\GM$ and $\GMQ$ be two strip complexes. Assume that there
is a locally compact group $G$ that acts continuously and properly on $\GM$
and such that the quotient of $\GM$ by $G$ is $\GMQ$ (as topological spaces). 
Let $\pi$ be 
the quotient map. Assume that $\GM$  is equipped with the data $(l,\phi,\psi)$ 
that induce a geometry, measure  and a Dirichlet form as discussed in the
preceding sections. Let $\{ H_t=e^{t\Delta} : t>0\}$
be the heat semigroup on $\GM$ associated with $(l,\phi,\psi)$.

Let a Markov semigroup  $\{ H_{0,t} : t>0\}$ acting on 
${\mathcal L}^\infty(\GMQ)$ be given
that satisfies $\lim_{t\rightarrow 0}H_{0,t}\phi=\phi$ for all
$\phi\in \mathcal C_c(\GMQ)$. Assume the following hypotheses.

\begin{enumerate}
\item  $\;(\GM,\rho)$ is complete and satisfies the volume condition
$$
\int_1^\infty \frac{r\,dr}{\ln V(\xi_0,r)}=\infty\,.
$$
\item $\; H_t$ commutes with the action of $G$ on $\GM$.
\item For any bounded function $\phi_0\in \mathcal C_c(\GMQ)$, the function
$u_0:(0,\infty)\times \GMQ\rightarrow \R$
defined by $u_0(t,\xi)=H_{0,t}\phi_0(\xi)$ is such that
$u=u_0\circ \pi$ is a weak solution of the heat equation on
$(0\,,\,T)\times \GM$.
\end{enumerate}

Then the semigroup $\{ \uH_t : t>0\}$, defined on
${\mathcal L}^\infty(\GMQ)$ by
$$
\uH_t f(\uxi)= H_t[f\circ \pi](\xi),\quad\text{where}
\;\;\pi(\xi)=\uxi,
$$ coincides with $H_{0,t}\,$. Consequently, if $(X_t,\mathbb P_\xi)$ and
$(X_{0,t},\mathbb P_{0,\xi_0})$ are the Markov process associated with
$\{ H_t : t>0\}$ and $\{H_{0,t} : t>0\}$ on $\GM$ and $\GMQ$, respectively,
then these processes are related by
$$
\mathbb P_{\xi_0} [X_{0,t}\in B]=\mathbb P_\xi [\pi(X_t)\in B],\quad\text{where}
\;\; \xi_0=\pi(\xi),
$$
for any measurable set $B\subset \GMQ\,$.
\end{thm}

\begin{proof} Let $\phi_0\in \mathcal C_c(\GMQ)$. Define
$$ 
f_{0,t}=H_{0,t}\phi_0\,,\quad \phi=\phi_0\circ \pi\,,\AND f_t=H_t \phi.
$$
It suffices to show that
$$ 
f_t=f_{0,t}\circ \pi.
$$
Since $\phi_0\in \mathcal C_c(\GMQ)$ and $\phi$ is a bounded, uniformly 
continuous function, it is clear that
$$
\forall\,\xi\in \GM,\;\;
\lim_{t\rightarrow 0} f_t (\xi) = \lim_{t\rightarrow 0}f_{0,t}\circ \pi(\xi) 
=\phi(\xi) \quad\text{for all}\;\xi\in \GM.
$$
We claim that both $u(t,\xi)=f_t(\xi)$ and
$\widetilde{u}(t,\xi)= f_{0,t}\circ \pi(\xi)$ are weak solutions of the heat 
equation on $(0\,,\,\infty)\times \GM$.
If we can prove this claim, the desired conclusion will follow from Theorem 
\ref{th-BdUniq}, that is, from the uniqueness property for the bounded Cauchy 
problem, because $\uH_t$ and  $H_{0,t}$  are determined on $L^\infty(\GMQ)$ 
by their action on $\mathcal C_c(\GMQ)$. Note that Theorem \ref{th-BdUniq} requires completeness of $\GM$ and
the volume growth condition that we are assuming here. 

By hypothesis,  $(t,\xi)\mapsto \widetilde{u}(t,\xi)=f_{0,t}\circ \pi (\xi)$
is a weak solution on $\GM$. This yields one half of the claim. To prove the 
other half, we use Theorem \ref{th-HKsmth} to see that the bounded function 
$f_t$ is a weak solution of the heat equation on $\GM$. Note that this indeed
requires some smoothness estimates on the heat kernel on $\GM$
since $f$ is not in ${\mathcal L}^2(\GM)$. Theorem \ref{th-HKsmth} is more than 
sufficient for this purpose. This yields the claim and completes the proof.
\end{proof}

\begin{rmks} {\bf (A)} Given that $\GMQ$ is the quotient of $\GM$ by a proper
continuous group action, Theorem \ref{th-proj} is based on three main 
hypotheses.
\begin{itemize}
\item Hypothesis (1) concerns $\GM$ and its meaning is quite clear: it implies
uniqueness for the bounded Cauchy problem for weak solution of the
heat equation.
\item Hypothesis (2) is also clear. It is satisfied whenever
the action of $G$ on $\GM$ is by measure-adapted isometries.
\item Hypothesis (3) is crucial and concerns the relation between the heat
equation on $\GM$ and a certain semigroup  on $\GMQ\,$. This hypothesis
captures a huge amount of information, and it is \emph{a priori} not entirely 
clear whether it is a reasonable hypothesis, or when it can actually be 
verified. We thus need study it in more detail.
\end{itemize}

\noindent 
{\bf (B)} It can occur that a group acts properly and continuously
on a strip complex $\GM$ equipped with data $\phi,\psi$ in an isometric,
measure adapted way, but that the quotient $\GMQ$ cannot be equipped with
corresponding data $\phi_0,\psi_0$ such that the quotient semigroup equals 
the semigroup on $(\GMQ,\phi_0,\psi_0)$. The problem comes from the function 
$\psi_0$ that defines the underlying  measure. Here is an example. 

Let $\Ms=\{0\}$ be trivial. Let $\Gamma$ be $\mathbb Z$ with edge lengths $1$, 
so that $\Gamma^1=\R$, equipped with $\phi\equiv 1$. Fix $\qq>1$ and let 
$\psi$ be defined by
$$
\psi(s)=\qq^{k-1}\,,\quad \mbox{if}\;\;s\in (2k\,,\,2k+2\bigr),
$$
so that $\psi$ is constant along pairs of edges sharing an odd integer
endpoint. Consider the obvious isometric group action by translation by an 
even integer. This is measure adapted (translation by  $2k$ changes the 
measure by a constant factor of $\qq^k$). The quotient of $\Gamma^1$ by this 
group action is the finite metric  graph $\Gamma^1_0$ with two vertices $a,b$ 
and two length $1$ edges $e,f$ joining $a$ to $b$. The vertices $a$ and $b$ 
correspond to even and odd integers, respectively. The problem comes from the 
following fact.

Assume that there is a function $\psi_0$ on $\Gamma^1_0$ so that
the projected semigroup  coincides with the semigroup on $(\Gamma_0,\psi_0)$.
On one hand, inspection shows that $\psi_0$ must be continous when passing 
through $a$ and it must have a jump of size $\qq$ when going through $b$.
On the other hand, $\psi_0$ must be constant over edges. These two
conditions are, of course, incompatible.
\end{rmks}

To prepare for the next proposition we make the following observations.
Let $\GM=\Gamma^1\times \Ms$ and $\GMQ =\Gamma^1_0\times \Ms_0$
be two strip complexes and  $G$ 
be a locally compact group that acts continuously and properly by 
isometries on $\GM$ with  quotient $\GMQ$ (as a topological space). 
Let $\pi$ be the quotient map. According to our definition 
(Definition \ref{def-iso}), isometries must 
send bifurcation manifolds to bifurcation manifolds and thus send 
$\GM^o$ to $\GM^o$. Hence the action of $G$ on $\GM$ induces an action of 
$G$ on the vertex set $V$ of $\Gamma$.

Observe further that for any $s\in \Gamma^1$ and $g\in G$, we must
have $g(\{s\}\times \Ms) =\{s'\}\times \Ms$ for some $s'\in \Gamma^1$
because for any $\tau,\tau'\in \Gamma^1$ and $x,y\in \Ms$, 
$\rho\bigl((\tau,x),(\tau',x)\bigr)=\rho\bigl((\tau,y),(\tau',y)\bigr)$. 
Indeed, this distance is equal to  the minimum of the  integral of 
$\sqrt{\phi}$ along any 
path in $\Gamma^1$ from $\tau$ to $\tau'$. Hence, the action
of $G$ on $\GM$ induces an action of $G$ on $\Gamma^1$. Moreover, 
topologically, the quotient of $\Gamma^1$ by this action is $\Gamma^1_0\,$. 
However, in general, it is not true that the quotient of $V$ by the action
of $G$ is $V_0$ because it might be the case 
that additional vertices and bifurcation manifolds are needed 
to turn $\GM/G$ into the  strip complex $\GMQ$.  This is best 
explained by two examples: 

(1) Take $\Gamma^1$ be the natural graph of $\mathbb Z$
($\equiv \mathbb R$ with the integers marked as vertices), $\Ms=\{0\}$, and 
$G= \mathbb Z$ acting by translation. Then the quotient is the circle with 
one marked point. This is not a strip complex (as a strip complex is 
required to have no loop) and we need to choose a second marked 
point to turn it into a strip complex.

(2)  Take $\Gamma^1$ as in (1) and $G=\{e,\sigma\}$ where $e$ is identity
and $\sigma$ is the reflexion with respect to $-1/2$. The quotient
is a half line with marking at $1/2$ and at the positive integers. 
To turn this into a strip complex, we need to add a vertex at the origin of 
the half line. 

Fortunately, this difficulty  (in the two examples above and in the general 
case) is solved by adding ``dummy''  middle vertices and corresponding
bifurcation manifolds in every strip $S^o_e\in \GM$, $e\in E$.
This yields a new strip complex $\GM'$ (isometric with $\GM$ as metric spaces,
and equivalent with $\GM$ for all analytic purposes)
with the same manifold $\Ms$ but the new graph 
$\Gamma'$ obtained by subdividing each edge of $\Gamma$ into two new edges 
with a new vertex in the middle. Furthermore, the action of $G$ on $\GM'$ 
(resp. $(\Gamma')^1$) is such that if $M_v$ and $M_w$ are two bifurcation 
manifolds (resp. $v,w$ are two vertices)  in the same orbit under the action 
of $G$ then the pair $\{v,w\}$ cannot be the pair of extremities 
of an edge $e$ in $E'$. It follows that  $(\Gamma')^1/G$ is naturally a 
metric graph with vertex set $V'_0=V/G$ and with no loops. Therefore, there is 
no loss of generality in assuming that $\Gamma^1_0=(\Gamma')^1/G$.

Consequently,  without loss of generality, we can assume that $\pi$ 
induces a natural graph homomorphism of 
$\Gamma$ onto $\Gamma_0\,$. The latter 
will also be denoted by $\pi$, so that we can speak about the vertices and
edges $\pi(v)$ and $\pi(e)$ of $\Gamma_0\,$, where $v \in V$ and $e \in E$,
respectively.

Consider a pair of  open strips $S^o\subset \GM$, $S^o_0\subset \GMQ$
with $\pi(S^o)=S^o_0\,$. Let  
$$
G_{S^o}=\{g\in G: g(S^o)=S^o\}/\{g\in G: g|_{S^o}=\mbox{\small \em id}\}
$$
be the effective quotient for the action of $G$ on $S^o$.  Since
any $g\in G$ such that $g\xi\in S^o$ for some $\xi\in S^o$ must send $S^o$ to
$S^o$, it follows that  $\pi(S^o)=S^o_0$ is also the (topological)
quotient of $S^o$ by the action of $G_{S^o}$ (see, e.g., 
{\sc Bourbaiki}~\cite[I.23]{BTG}),
and for any function $u_0$ on $\GMQ\,$, we have
\begin{equation}\label{local-strip}
u_0\circ \pi|_{S^o}= u_0|_{S^o_0}\circ \pi^{S^o}
\end{equation}
where  $\pi^{S^o}$ is the projection map from $S^o$ to $S_0^o$.

Note that $G_{S^o}$ acts by isometries on the manifold $S^o\,$.
In what follows we will assume that $G_{S^o}$ is a Lie subgroup of the
group of isometries of $S^o$ and that
$$
\pi^{S^o}: 
\Bigl(S^o=I\times \Ms\,,\, 
\phi\,\bigl((ds)^2+g(\cdot,\cdot)\bigr)\Bigr)  \rightarrow
\Bigl(S_0^o=I_{0}\times \Ms_0\,,\, 
\phi_0\,\bigl((d\tau)^2+g_0(\cdot,\cdot)\bigr)\Bigr)
$$
is a Riemannian submersion. This implies that the action of $G_{S^o}$
on $S^0$ is free. Moreover, $\pi^{S^o}$ sends any set of the form
$\{s\}\times \Ms$ to some set of the form $\{\tau\}\times \Ms_0$ and,
for any $ f_0\in \mathcal C^\infty(S^o_0)$ and any $(s,x)\in S^o$
with $\pi^{S^o}(s,x)= (\tau,x_0)$, we have
\begin{equation} \label{partialssubm}
\frac{1}{\phi(s)}\bigl|\partial_s f_0\circ \pi^{S^o} (s,x)\bigr|^2 =
\frac{1}{\phi_0(\tau)}\bigl|\partial_\tau f_0 (\tau,x_0)\bigr|^2
\end{equation}
and
\begin{eqnarray}\label{Lapsubm}
\lefteqn{
\frac{1}{\phi(s)}\Bigl[\partial^2_s +\Delta_{Ms} +
[\partial_s \log \phi(s)^{(n-1)/2}]\,\partial_s\Bigr] f_0 \circ \pi^{S^o}(s,x)}
							   \hspace{.5in}&&\\
&&\nonumber=
\frac{1}{\phi_0(\tau)}\Bigl[\partial^2_\tau +\Delta_{\Ms_0}+
[\partial_\tau \log \phi_0(\tau)^{(n-1)/2}]\,\partial_\tau\Bigr]  
f_0 (\tau,x_0).
\end{eqnarray}
This follows from the fundamental property of a Riemannian submersion and the 
fact that the expressions in  \eqref{Lapsubm} are the Laplace operators
of the relevant Riemannian metrics. Observe  that the weight functions
$\psi$ and $\psi_0$ do not appear in this formula.

\begin{pro}\label{pro-goodact}
Let $\GM$ and $\GMQ$ be two strip complexes. Assume that there
is a locally compact group $G$ that acts continuously and properly on $\GM$
and such that the quotient of $\GM$ by a $G$ is $\GMQ\,$.
Let $\pi$ be the quotient map. Assume that $\GM$  and $\GMQ$ are equipped 
with the data  $(\phi,\psi)$ and $(\phi_0,\psi_0)$, respectively, 
that induce a geometry, measure  and a respective Dirichlet form as discussed 
above. Assume furthermore that the following hypotheses are satisfied.
\begin{enumerate} 
\item  $G$ acts on $\GM$ by isometries and $\Gamma_0$ is the quotient 
of $\Gamma$ under the induced  action of $G$.
\item  For any  edge $e\in E$,  the group  $G_{S^o_e}$ is a Lie subgroup 
of the isometry group of $S^o_e\,$, the projection map $\pi^{S^o_e}$ is a 
Riemannian submersion of $S^o_e$ onto $S^o_0=\pi(S^o_e)\subset \GMQ\,$, and
\item 
there exists a constant $A(e)\in (0\,,\,\infty)$ such that
$$
\psi_e(s)=A(e)\,\psi_0|_{S^o_0}(\tau)
$$ 
for any $s,\tau$ such that $\pi^{S^o}(s,x)=(\tau,x_0)$ 
for some $x\in M$ and $x_0\in M_0\,$.
\item For any pair of vertices $v\in V$ and $v_0\in V_0$ such that
$\pi(M_v)=M_{0,v_0}\,$, there exists a constant $a(v)\in(0\,,\,\infty)$ such 
that
$$
\sum_{e\in E_{v}\,:\, \pi(e)=e_0} \psi_e(v)
= a(v)\, \psi_{0,e_0}(v_0) \quad\text{for all}\; e_0\in E_{v_0}\,.
$$
\end{enumerate}
Then, for any $T>0$ and any function
$u_0\in\mathcal C^\infty\bigl((0,T)\times \GMQ\bigr)$ which is a time regular
weak solution of the heat equation  on $(0\,,\,T)\times \GMQ\,$, the
function $u=u\circ \pi$ is a time regular weak solution of the heat equation 
on $(0\,,\,T)\times \GM$.
\end{pro}

\begin{proof}
Because of \eqref{local-strip} and assumption (2), $u=u_0\circ \pi$ and its 
time derivatives $\partial^k_tu$ are in $\mathcal C^\infty(\GM).$ For such 
a function, being a weak solution of the heat equation means:
$$ 
\left\{\begin{array}{l}
\partial_t u= \Lap u = \dfrac{1}{\phi(s)}
\bigl[\partial_s^2+\Delta_{\Ms} +\eta \partial_s\bigr]u=0, \quad\text{where}\; \; 
\eta=\partial_s \ln(\phi^{(n-1)/2} \psi)\,,
\\[5pt]
\displaystyle\sum_{e\in E_{v}}\psi_e(v)\,(\nnu_{v,e}\,,\nabla u_{e}) =0 \quad
\mbox{along}\; M_{v}\, \;\text{for all}\;v\in V.
\end{array}\right.
$$
That $u$ satisfies the first of those two identities follows by careful 
inspection using  \eqref{local-strip}, assumption (2), \eqref{Lapsubm} and 
assumption (3). The  second line identity similarly from \eqref{local-strip},
assumption (2), \eqref{partialssubm} and assumption (4).
\end{proof}

\begin{exa}
Let $\GM$ be a strip complex equipped with the data $\phi$ and $\psi$. Assume 
that the isometry group $G$ of $(\Ms,g)$ acts transitively on $\Ms$. This 
group also acts on $\GM$ in an obvious way, and this action
is measure adapted (in fact, measure preserving) and isometric.
The quotient of $\GM$ by this action is the $1$-dimensional complex
$\Gamma^1$. For each open strip $S^o$, $G_{S^o}$ is isomorphic with $G$ itself,
and assumption (2) of Proposition \ref{pro-goodact} is obviously satisfied.
Assumptions (3)--(4) of Proposition \ref{pro-goodact} are satisfied
if we equip $\Gamma^1$ with the data $\phi, A_0\,\psi,$ where $A_0$ is 
any fixed positive constant.

The same applies if $G$ is a subgroup of the isometry group that acts freely 
and properly on $\Ms$ with quotient $\Ms_0$. Then there exists a unique 
Riemannian structure on $\Ms_0$ that makes the quotient map a Riemannian 
submersion. The quotient of $\GM$ under the natural action of $G$ is $\GMQ$ 
with $\Gamma_0=\Gamma$. For each open strip $S^o\,$, the group $G_{S^o}$ is 
again isomorphic to $G$ itself, and assumption (2) of Proposition 
\ref{pro-goodact} is obviously satisfied. Assumptions (3)--(4) of Proposition 
\ref{pro-goodact} are satisfied if we equip $\GMQ$ with the data $\phi$ and 
$A_0\, \psi$ for any fixed positive constant $A_0\,$.
\end{exa}

\begin{exa}\label{ex-goodact} Let $\GM$ be a strip complex.
Assume that $G$ is a subgroup of the group
of automorphisms of the non-oriented version of the graph $\Gamma$. 
By adding dummy vertices in the middle of edges if necessary, we can
assume that the quotient  $\Gamma_0=\Gamma/G$ has no loops.
For any Riemannian manifold $\Ms=\Ms_0$, the group $G$ has a natural action on
$\GM$ with quotient $\GMQ\,=\Gamma_0\times \Ms_0=\Gamma_0\times \Ms$. 
Let $\pi$ be the quotient map
from $\Gamma^1$ to $\Gamma^1_0\,$. In particular, $\pi$ maps the 
edge set of $\Gamma$ onto the edge set of $\Gamma_0\,$. 
Fix data $\phi_0$ and $\psi_0$ on $\GMQ$ 
and equip $\GM$ with $\phi=\phi_0\circ \pi$. Then $G$ acts on
$\GM$ by isometries. Next, we consider the conditions (3)--(4) of Proposition
\ref{pro-goodact}. Condition (3) involves  numbers $A(e)>0$, $e\in E$, 
such that $\psi_e= A(e)\,\psi_{0,e_0}\circ \pi|_{S^o_e}$\,. 
Given that condition (3) is satisfied, condition (4)
requires that 
$$
\sum_{ e\in E_v\,:\,\pi(e)=e_0}A(e)=a(v) \quad\text{for all}\; v\in V\,,\;
\,e_0 \in E_{\pi(v)}\,.
$$

Let us examine some special cases. 

\smallskip

\noindent
{\bf (A)} First, assume that for any vertex $v$ of $\Gamma$, we have 
$\deg_{\Gamma}(v) = \deg_{\Gamma_0}\bigl(\pi(v)\bigr)$. Then the restriction
of $\pi$ from $E_v$ to $E_{\pi(v)}$ is bijective, or in other words, $\pi$ is a
graph covering. In this case, the above condition means that $A(e)=A(e')$ if 
the edges $e$ and $e'$ have a common end vertex. Since our graphs are 
connected, this actually  implies that $A(e)=A$ is a constant, that is, 
$\psi=A \cdot \psi_0\circ \pi$.

\smallskip

\noindent
{\bf (B)} Second, consider the specific example where $\Gamma= \T_2$ is the 
regular tree  with degree $3\,,$ drawn with respect to a reference end
$\varpi$ as in Figure~2. The graph $\Gamma_0$ is the two-way-infinite path, 
which we denote by $\Z$ (which is, more precisely, the vertex set of  
$\Gamma_0\,$, while the associated $1$-complex is $\R$). The group $G$ is the 
group of all graph automorphisms of the tree that fix every horocycle,
and the projection is $\pi = \hor$, the Busemann function with respect to
$\varpi$. Here, the projection is obviously not a graph covering.
For simplicity, we assume that all edges have length $1$ and that 
$\phi,\phi_0\equiv 1$. Furthermore, we assume that $\psi_0$ is constant on 
each edge of $\Gamma_0=\mathbb Z$. 
Recall that in this specific example, every vertex $v$ has one neighbouring 
vertex $v^-$ in the ``preceding'' horocycle and is itself the predecessor
of its ``forward'' neighbours $w_1, w_2$ that satisfy $w_i^-=v$. 
(This notation should not be mixed up
with the one for the endpoints $e^-$ and $e^+$ of an edge $e$.) 
If $\hor(v)=k$ then $e_v= [v^-,v]$ is the only edge in $E_v$ that projects
onto the edge $[k-1,k]$ of $\Z$. Therefore $A(e_v)=a(v)$. On the other hand,
both edges $e_{w_1}$ and $e_{w_2}$ project onto the edge $[k,k+1]$ of $\Z$.
Therefore the above condition can be rewritten in terms of the positive
function  $v \mapsto a(v)$. In order to be feasible, it is necessary and 
sufficient that it satisfies $a(w_1)+a(w_2)= a(v)$ for any vertex $v$ of
$\T$, where the $w_i$ are its forward neighbours.
Because $\T_2$ is a tree, we can construct infinitely many functions that 
satisfy  this property, and hence there are infinitely many functions
$\psi$, constant on open edges, so that conditions (3) and (4) are satisfied,
whenever the function $\psi_0$ is chosen to have constant value,
say $b_k\,$, on each open strip $(k-1,k) \times \Ms$ of $\GMQ$.
One solution for $\psi$ is given by
$$
\psi\big|_{S_e^o} \equiv 2^{-k}b_{k}\,,\quad\text{when}\;\;
\pi(e)=[k-1,k].
$$
This is the only solution for which the corresponding group action
is measure adapted.

\smallskip

\noindent
{\bf (C)} Consider the situation described in Theorem \ref{thm-projections}
concerning various projections of $\HT(\pp,\qq)$. The hypotheses (1) and (2)
of Theorem \ref{th-proj} are verified,  and  hypothesis (3) is also satisfied 
because of Proposition \ref{pro-goodact}.  Hence Theorem \ref{thm-projections}
 follows from Theorem \ref{th-proj} and 
Proposition \ref{pro-goodact}. Note that Proposition \ref{pro-goodact} makes 
heavy use of the results
of Section 5. Further related uniqueness theorems are given in the next two 
sections.
\end{exa}

\section{Uniqueness of the heat semigroup}\label{uniq}

Throughout this section, we use the basic setting of a strip complex with
data, distance, measure, Dirichlet form, Laplacian and heat semigroup as already
specified at the beginning of Section \ref{heat}. 
Our aim is to show that, in some strong sense,
there is only one semigroup of operators whose generator
coincides with the Laplacian $\Delta$ on a certain space of
smooth compactly supported functions. This property is important
in many applications.  We will discuss two different uniqueness results:
one concerns uniqueness on $\mathcal C_0(\GM)$, whereas the other
concerns uniqueness on ${\mathcal L}^2(\GM)$.

\subsection*{A. A candidate for a core of the infinitesimal generator}
In this section we introduce a very specific space, 
$\mathcal D^\infty_c\,$,  of compactly supported
smooth functions on $\GM$ that is a good candidate to be a core 
for the generator of the heat semigroup, either on $\mathcal C_0(\GM)$ or on 
$\mathcal L^2(\GM)$. In some cases, we will be able to show that 
$\mathcal D^\infty_c$ is indeed a core. Please note that
the spaces $\mathcal D^\infty$ and $\mathcal D^\infty_c$ introduced below   
depends on the fixed data $(l,\phi,\psi)$ on $\GM$.

\begin{dfn} The space $\mathcal D^\infty$ is
the space of all functions $f$ in $\mathcal C^\infty(\GM)$ such that
\begin{enumerate}
\item For any integer $k=0,1,\dots, $ any $v\in V$ and $e,e'\in E_v$,
$$
\mbox{Tr}^{S_e}_{\Ms_v}(\Lap^kf)=\mbox{Tr}^{S_{e'}}_{\Ms_v}(\Lap^kf).
$$
This means that the functions $\Lap^kf$, originally only defined
and continuous on $\GM^o$, are in fact
continuous functions on $\GM$ (after proper extension by continuity)
and thus in $\mathcal C^\infty(\GM)$.
\item For any integer $k=0,1,\dots, $ and $v\in V$
$$
\sum_{e\in E_v} \psi_e(v)\, (\nnu_{v,e}\,,\nabla \Lap^kf_e)=0 \quad\mbox{along}\; 
M_v\,. 
$$
This means that each function $\Lap^kf\in \mathcal C^\infty(\GM)$ satisfies the
bifurcation condition along any bifurcation manifold $\Ms_v\,$, $v\in V$.
\end{enumerate}
The space $\mathcal D^{\infty}_c$ is the subspace of all compactly supported 
functions in $\mathcal D^\infty$.
\end{dfn}

\begin{rmk}\label{xcont}
Fix a coordinate chart $(U;x_1,\dots,x_n)$ in $\Ms$.
Observe  that any function $f$ in $\mathcal C^\infty(\GM)$
viewed as a function of $(s,x)\in \Gamma^1\times U$
actually has continuous partial derivatives of all orders
$\partial^\kappa_xf(s,x)$ in the $x$ direction, but not in the $s$ direction 
in general. It follows that the continuity condition on $\Lap f$ reduces
to the continuity of 
$$
\partial^2_s f+\eta(s)\, \partial_s f
$$
across any bifurcation manifold $\Ms_v\,$. The bifurcation condition implies 
that, typically, the function $\partial_s f$ is not continuous across 
bifurcation manifolds. It follows that, typically, $\partial_s^2 f$ is not 
continuous and neither are $\partial^k_sf$, $k\ge 3$. An important 
consequence of this is that $\mathcal D^\infty$ and $\mathcal D^\infty_c$ 
are not algebras under pointwise multiplication.
\end{rmk}
\begin{rmk} Note that $\mathcal D^{\infty}_c$
is a subspace of $\mbox{Dom}(\Delta ^k)$ for every $k \ge 1$, and
$$
\Delta^k =\Lap^k \quad\mbox{on}\;\;\mathcal D^\infty_{c}\,.
$$
\end{rmk}

\begin{lem} \label{lem-C0d} The space $\mathcal D^{\infty}_{c}$ is dense in 
$\mathcal C_0(\GM)$ for the uniform topology.
\end{lem}

\begin{proof} This important result is an immediate corollary of Lemma 
\ref{lem-dens0} since we have 
$\mathcal C^\infty_{c,c}(\GM) \subset \mathcal D^\infty_{c}\,$. Indeed,
$\mathcal C^\infty_{c,c}(\GM)$ is the subspace of those functions $f$
in $\mathcal C^\infty_c(\GM)$ whose partial derivative $\partial_s f$ along
$\Gamma^1$ vanishes in a neighbourhood of any bifurcation manifold.
The desired inclusion thus follows from Remark \ref{xcont} above.
\end{proof}

The following is a simple corollary of Theorem \ref{th-HKsmth}.

\begin{thm} \label{th-HKD}
For every fixed $t>0$, $\zeta\in \GM$, $k=0,1,\dots$,  every relatively 
compact coordinate chart $I\times U\ni \zeta$ in $\GM$ and
$\kappa=(\kappa_0,\kappa_1,\dots\kappa_n)$,  the function
$$
\xi\mapsto \partial^k_t\partial^\kappa_\zeta h(t,\zeta,\xi)
$$
belongs to $\mathcal D^\infty$.
\end{thm}

\subsection*{B. Uniqueness of the heat semigroup on $\mathcal C_0(\GM)$}

Consider the operator $(\Lap ,\mathcal D^\infty_c)$ as a linear, densely defined  
operator on $\mathcal C_0(\GM)$. Recall that indeed, $\mathcal C^\infty_c$ is 
dense in $\mathcal C_0(\GM)$ for the uniform topology, see Lemma \ref{lem-C0d}.
We claim that $(\Lap ,\mathcal D^\infty_c)$ satisfies the positive 
maximum principle. That is, if $\xi_0\in \GM$ and $f\in \mathcal D^\infty_c$
are such that $\max_{\GM}\{ f\}=f(\xi_0)\ge 0$, then $\Lap f(\xi_0)\le 0$.

Indeed, if $\xi_0$ is not on a bifurcation manifold, this follows from the 
usual maximum principle. If $\xi_0=(v_0,x_0)$ is on a bifurcation manifold, 
let $(U;x_1,\dots,x_n)$ be a local coordinate chart in $\Ms$ around $x_0\,$.
Since $f\in \mathcal D^\infty_c$ is maximal at $\xi_0\,$,
the first order partial derivatives at $\xi_0$ along $\Ms$ must be $0$
and we must have $\partial^2_{x_i}f(\xi_0)\le 0$, $i=1,\dots, n$.
It follows that $\Delta_{\Ms}f(\xi_0)\le 0$.

Moreover, in any strip $S_e$ containing $\xi_0=(v_0,x_0)$,
the outward normal derivatives $\bigl(\nnu_{v,e}\,,\nabla f_e(\xi_0)\bigr)$ 
must be greater or equal to $0$. Hence, the bifurcation condition implies that
$\bigl(\nnu_{v,e}\,,\nabla f_e(\xi_0)\bigr)=0$. It follows that in any strip 
$S_e=I_e\times \Ms$ around $\xi_0\,$, we must have 
$\partial_s^2f_e(\xi_0)\le 0$. Hence
$$
\Lap f(\xi_0)= \frac{1}{\phi(\xi_0)}[\partial_s^2f(\xi_0)+\Delta_{\Ms}f(\xi_0)]
\le 0.
$$

Without further assumption on $\GM$, we do not know how to show that 
$(\Lap ,\mathcal D^\infty_c)$ admits an  extension that is the infinitesimal 
generator of a contraction semigroup on $\mathcal C_0(\GM)$. The difficulty 
lies in proving that the range $(\lambda \,\Id-\Lap )\mathcal D^\infty_c$ is dense 
in $\mathcal C_0(\GM)$ for some $\lambda>0$, that is, that 
$(\Lap ,\mathcal D^\infty_c)$ is closable in $\mathcal C_0(\GM)$.
However, by the results of {\sc van Casteren and Okitaloshima}~\cite{VC}, 
\cite{OVC}, we have the following  
\cite[Theorem 3.6 and Proposition 3.7]{OVC}: if $(\Lap ,\mathcal D^\infty_c)$ is 
closable, then its closure is the only linear extension of 
$(\Lap ,\mathcal D^\infty_c)$ that is the infinitesimal generator of a Feller 
semigroup (that is, a strongly continuous semigroup of contractions on
$\mathcal C_0(\GM)$ preserving positivity). This, together with Theorem
\ref{th-UR}, yields the following result.

\begin{thm}\label{th-UC0}
Let $\GM$ be a strip complex equipped with a geometry and measure
as above. Let  $h(t,\xi,\zeta)$  be the heat kernel associated with the Dirichlet
form $\bigl(\mathcal E,\mathcal \mathcal W^1_0(\GM)\bigr)$, where
$(t,\xi,\zeta)\in (0,\infty)\times \GM\times \GM\,$.
Assume that $(\GM,\rho)$ is complete and that there are constants $D,P,r_0$
such that \emph{(i)} and \emph{(ii)} hold.
\begin{itemize}
\item[(i)] For any $\xi\in \GM$ and $r\in (0\,,\,r_0)$, we have the doubling property
$V(\xi,r)\le D \,V(\xi,2r)$.
\item[(ii)] For any $\xi\in \GM$  and $r\in (0,r_0)$, setting $B=B(\xi,r)$,
$$
\int_B|f-f_B|^2\, d\mu\le P\, r^2\int_B|\nabla f|^2\,d\mu \quad
\text{for every}\;\;f\in \mathcal W^1(B)\,, \quad\text{where}\;\;
f_B=\frac{1}{\mu(B)}\int_B f\,d\mu\,.
$$
\end{itemize}
Then the densely defined linear operator
$(\Lap,\mathcal D^\infty_c)$ on $\mathcal C_0(\GM)$  is closable and its closure 
$\bigl(\olA,\mbox{\em Dom}(\olA)\bigr)$ is the infinitesimal generator of the 
Feller semigroup defined by
$$ 
\mathcal C_0(\GM)\ni  f\mapsto  e^{t\olA}f\,,\; t > 0\,,\quad\text{where}\quad
e^{t\olA}f(\xi)= \int_{\GM}h(t,\xi,\zeta)\,f(\zeta)\,d\mu(\zeta).
$$
Moreover, if $\bigl(\widetilde{\Lap },\mbox{\em Dom}(\widetilde{\Lap })\bigr)$
is an extension of $(\Lap ,\mathcal D^\infty_c)$
and is the infinitesimal generator of a Feller semigroup then
$\bigl(\widetilde{\Lap },\mbox{\em Dom}(\widetilde{\Lap })\bigr)
=\bigl(\olA,\mbox{\em Dom}(\olA)\bigr)$.
\end{thm}

\begin{rmk} It follows from the results in \cite{VC} and \cite{OVC} that,
under the hypotheses of Theorem \ref{th-UC0},
the martingale problem for the operator $(\Lap ,\mathcal D^\infty_c)$
is uniquely solvable (for any starting point $\xi\in \GM$).
See \cite[Theorem 3.6]{OVC}.
\end{rmk}

\subsection*{C. Uniqueness of the heat semigroup on ${\mathcal L}^2(\GM)$}

Let us observe that, because of the possibility to impose various boundary 
conditions, uniqueness on ${\mathcal L}^2(\GM)$ cannot hold unless we make the assumption
that $(\GM,\rho)$ is complete.

\begin{dfn}We say that a continuous function 
$\rho_0:\GM\rightarrow (0,\infty)$
is a strip-adapted exhaustion function 
if it has the following properties.
\begin{itemize}
\item The function $\rho_0$ belongs to $\mathcal C^\infty(\GM)$.
\item For any edge $e\in E$ and any $x\in \Ms$ the function
$s\mapsto \partial_s \rho_e(s,x)$ has compact support in $(e^-\,,\,e^+)$.
\item The function $\rho_0$ tends to infinity at infinity.
\item The functions  $|\nabla \rho_0|$ and $|\Lap \rho_0|$ are bounded on $\GM$.

\end{itemize}
\end{dfn}

Note that a strip-adapted exhaustion function is a continuous smooth function 
on $\GM$ which is locally constant in the direction of $\Gamma^1$ near 
each bifurcation manifold.
The existence of such exhaustion functions is a 
non-trivial matter that
will be discussed in Section \ref{adapted}. 

\begin{dfn}
A sequence of continuous compactly supported functions $\varrho_n$ is called
a strip-adapted approximation of $\uno$ 
if the following holds.
\begin{itemize}
\item Each $\varrho_n$ belongs to $\mathcal C^\infty_{c,c}(\GM)\,$.
\item Each $\varrho_n$ takes values in $[0\,,\,1]$,  and
$\;\displaystyle  \lim_{n\rightarrow \infty}\varrho_n(\xi)=1 \;$
for all $\xi\in\GM$.
\item There is a constant $C$ such that
$|\nabla \varrho_n|\le C$, $\;|\Delta \varrho_n|\le C$, and
for all $\xi\in\GM$,
$$
\lim_{n\rightarrow \infty}|\nabla \varrho_n(\xi)|
=\lim_{n\rightarrow \infty}|\Delta\varrho_n(\xi)|=0.
$$
\end{itemize}
\end{dfn}

\begin{rmks}\label{rmks-adapted} (a) If a strip-adapted exhaustion  function $\rho_0$ exists then
a strip-adapted approximation of $\uno$ is easily obtained by setting
$\varrho_n(\xi)=\theta\bigl(\rho_0(\xi)/n\bigr)$, where $\theta$ is a smooth, 
compactly supported function of \emph{one} variable taking value in $[0\,,\,1]$ 
and such that $\theta \equiv 1$ in a neighbourhood of $0$.
\\[5pt]
(b) Let $\varrho_n$ be
a strip-adapted approximation of $\uno$. Then 
$\varrho_n f\in \mathcal D^\infty_c$ for any  $f\in \mathcal D^\infty$. 
Compare this with the fact that, in general, $\varrho\in \mathcal D^\infty_c$ 
and $f\in \mathcal D^\infty$ does not imply $\varrho f\in 
\mathcal D^\infty_c$.
\end{rmks}

\begin{thm} \label{th-SA}
The operator $(\Lap , \mathcal D^{\infty}_{c})$ is symmetric on 
${\mathcal L}^2(\GM)$.
If $(\GM,\rho)$ is complete and there exists a strip-adapted approximation of 
$\uno$ then the symmetric 
operator $(\Lap , \mathcal D^{\infty}_{c})$ is essentially self-adjoint on 
${\mathcal L}^2(\GM)$, and its unique self-adjoint extension
is $\bigl(\Delta, \mbox{\em Dom}(\Delta)\bigr)$.
\end{thm}

\begin{rmk} When considering Theorem \ref{th-SA}, the reader should recall that
the relevant underlying data include the graph $\Gamma=(V,E)$, the
Riemannian manifold $(M,g)$, the function 
$\phi\in \mathcal C^\infty(\Gamma^1)$ which is part of the definition of the 
geometry on $\GM$ and plays a crucial role on whether $(\GM,\rho)$ is complete 
or not, as well as the function 
$\psi\in \mathcal S^\infty(\Gamma^o)$
which appears in the Definition \ref{dfn-mu} of the underlying measure $\mu$. 
Indeed, ${\mathcal L}^2(\GM)$ is the ${\mathcal L}^2$-space relative to that specific measure
$\mu$. It is interesting to observe how these different parameters enter
the definition of $\Delta$ and that of $\Lap $. Concerning $\Lap $, the functions
$\phi$ and $\psi$ appear in the formula defining $\Lap $ on each open strip.
However, the possible jump discontinuities of $\phi$ and / or $\psi$ only 
appear in the definition of $\mathcal D^\infty_c$ via the bifurcation condition.
This clearly shows that one cannot replace $\mathcal D^\infty_c$ by 
$\mathcal C^\infty_{c,c}(\GM)$ in Theorem \ref{th-SA} because then the role of 
the possible jumps of the functions $\phi$ and $\psi$ is lost.
\end{rmk}

The proof of Theorem \ref{th-SA} requires a number of lemmas. The symmetry of
$(\Lap , \mathcal D^{\infty}_{c})$ on ${\mathcal L}^2(\GM)$ follows from the various 
definitions by inspection. Let $\bigl(\Lap^*, \mbox{Dom}(\Lap^*)\bigr)$ be the
adjoint of $(\Lap , \mathcal D^{\infty}_{c})$.

\begin{lem} For any function $f\in \mathcal D^\infty\cap \mbox{\em Dom}(\Lap^*)$,
one has $\Lap^* f= \Lap f$ in ${\mathcal L}^2(\GM)$.
\end{lem}

\begin{proof} By definition, for any $f\in \mbox{Dom}(\Lap^*)$ and 
$h\in \mathcal D^{\infty}_{c}\,$, we have
$$
\langle \Lap^* f,h\rangle = \langle f,\Lap h\rangle\,,
$$
where $\langle\cdot,\cdot\rangle$ is the inner product on ${\mathcal L}^2(\GM)$.
But for $f\in \mathcal D^\infty$ and $h\in \mathcal D^{\infty}_{c}\,$,
Green's formula in each strip and the bifurcation conditions imposed on $f$ and
$h$ show that
$$
\langle f,\Lap h\rangle =\langle \Lap f,h\rangle.
$$
This proves the desired result.
\end{proof}

\begin{lem}\label{lem-A1}
Let $f\in  \mbox{\em Dom}(\Lap^*)$, $h \in \mathcal D^\infty$, and suppose that 
$h,\Lap h\in {\mathcal L}^2(\GM)$.  Then
$$
\langle \Lap^* f,h\rangle= \langle f,\Lap h\rangle.
$$
\end{lem}

\begin{proof} Consider the sequence $h_n=\varrho_n h$, where $\varrho_n$ is 
a strip-adapted approximation of $\uno$. Then $h_n\in \mathcal D^\infty_c$
and $h_n \to h$ as well as  $\Lap h_n \to \Lap h$ in ${\mathcal L}^2(\GM)$. Hence the desired 
equality follows from the fact that 
$\langle \Lap^*f,h_n\rangle= \langle f,\Lap h_n\rangle.$
\end{proof}

\begin{lem} \label{lem-A2} For any function $f\in \mbox{\em Dom}(\Lap^*)$ and 
$t>0$, the function $f_t=e^{t\Delta}f$ is in $\mbox{\em Dom}(\Lap^*)$
and 
$$
\Lap^* f_t=e^{t\Delta}\Lap^*f.
$$
\end{lem}

\begin{proof} For any $h$ in  ${\mathcal L}^2(\GM)$, the function $h_t=e^{t\Delta}h$ is a 
global weak solution of the heat equation that is time regular
to infinite order. By Theorem \ref{th-Heatsm} this implies
that $h_t\in \mathcal D^\infty$. Obviously, $h_t$ and  $\Lap h_t=\Delta h_t$
are also in ${\mathcal L}^2(\GM)$. Now, for $f\in \mbox{Dom}(\Lap^*)$ and
$h\in \mathcal D^{\infty}_c$, we have
$$
\langle e^{t\Delta} \Lap^*f, h\rangle=
\langle \Lap^* f, e^{t\Delta} h\rangle.
$$
Since $h_t=e^{t\Delta} h$ is in $\mathcal D^\infty$ and both $h_t$ and 
$\Lap h_t$ are in ${\mathcal L}^2(\GM)$, Lemma \ref{lem-A1} gives

\begin{eqnarray*}
\langle e^{t\Delta} \Lap^*f, h\rangle
&=&
\langle \Lap^*f, h_t\rangle 
=\langle f, \Lap h_t\rangle =\langle f, \Delta e^{t\Delta}h\rangle\\
&=&\langle f, e^{t\Delta} \Delta h\rangle 
= \langle e^{t\Delta} f, \Lap h\rangle = \langle f_t, \Lap h\rangle\,.
\end{eqnarray*}
This proves that $\Lap^* f_t=e^{t\Delta} \Lap^*f$ as desired.
\end{proof}

The next lemma will complete the proof of Theorem \ref{th-SA}.

\begin{lem} $\mathcal D^{\infty}_{c}$ is dense in $\mbox{\em Dom}(\Lap^*)$ in the 
graph norm.
\end{lem}

\begin{proof} 
Approximate $f\in \mbox{Dom}(\Lap^*)$ by
$f_t=e^{t\Delta}f$, where $t\rightarrow 0$. Then $f_t$ converges to $f$ in
${\mathcal L}^2(\GM)$ and, by Lemma \ref{lem-A2}, $\Lap^*f_t$ also converges to $\Lap^*f$ in 
${\mathcal L}^2(\GM)$. This shows that $\mathcal D^{\infty}\cap \mbox{Dom}(\Lap^*)$
is dense in $\mbox{Dom}(\Lap^*)$ in the graph norm. Now, we use multiplication 
by the strip-adapted sequence $\varrho_n$  that approximates $\uno$ and set 
$h_n=f_{1/n}\,\varrho_n$ to obtain the desired conclusion.
\end{proof}

\begin{rmk}\label{graphcase} Assume that $M=\{0\}$ is a singleton so that
$(\GM,\phi,\psi)$ reduces to the metric graph $\Gamma^1$ equipped with
the data $\phi,\psi$.  Assume that $(\Gamma^1,\rho)$ is complete.
In this case, the symmetric operator
$(\Lap ,\mathcal D^\infty_c)$ is always essentially self-adjoint on 
${\mathcal L}^2(\Gamma^1,\mu)$. This is proved in \cite{BSCg} following the 
argument used for complete Riemannian 
manifolds by {\sc Strichartz}~\cite{Str}. It is not clear that this argument 
can be adapted to the case when $M\neq \{0\}$. The difficulty lies in showing 
that any solution $f\in \mbox{Dom}(\Lap^*)$ of the equation $\Lap^*f=\lambda f$
is in fact in $\mathcal W^1_{\mbox{\tiny loc}}(\GM)$. On a manifold, this follows
from local ellipticity. On a graph, it can be checked by an adhoc argument
using very much the $1$-dimensional nature of the underlying space. See
\cite{BSCg}.
\end{rmk}

\section{Strip-adapted approximations of 
$\uno$}\label{adapted}

Unfortunately, the existence of a strip-adapted approximation of $\uno$
is a difficult question in full 
generality. Even in the case of complete Riemannian manifolds,
an adapted approximation of $\uno$ is not known to exist in general.
The proof of the essential self-adjointness of the Laplacian
(see, e.g., \cite{Str}) on a complete Riemannian manifold has to
avoid the use of an adapted approximation of $\uno$ and, instead, makes use
of the fact that the adjoint is an elliptic operator in the sense of 
distributions. See Remark \ref{graphcase} regarding the graph case. 
Whether or not that approach can be made to work
in the present setting is not clear, the main question being whether or not
one can prove that
$$
\mbox{Dom}(\Lap^*)\subset \mathcal W^1_{\mbox{\tiny loc}}(\GM).
$$
This appears to be a rather subtle question although one would conjecture that
the answer is ``yes''.

In this section we construct strip-adapted exhaustion functions
(or strip-adapted approximations of $\uno$) 
in a number of different special cases. 
We start with some simple-minded constructions.

\begin{pro} Assume that $(\Ms,g)$ is a complete Riemannian manifold
which admits an adapted approximation $(\varrho_{\Ms,n})$ of $\,\uno$.
Assume that the underlying metric graph $\Gamma$ satisfies
$$
l_*=\inf_E \{l_e\}>0\,,
$$
that is, edge lengths are bounded below.
Assume that $\GM$ is equipped with its bare strip complex structure,
that is, $\phi\equiv 1$ and $\psi\equiv 1$.
Then $\GM$ admits a strip-adapted approximation of $1$.
\end{pro}

\begin{proof} Let us first construct an edge-adapted 
exhaustion $s\mapsto\rho_1(s)$
on the  one dimensional complex $\Gamma^1$. (Here, the strips are the edges, so
that we use ``edge-adapted'' instead of ``strip-adapted''.)
Fix $\epsilon\in (0\,,\,l_*/8)$.
On $\Gamma^1$, consider a function $\alpha\in \mathcal C^\infty(\Gamma^1)$
with the property that for each edge $e$, the restriction $\alpha_e$
of $\alpha$ to $(e^-\,,\,e^+)$ has compact support in
$(e^-+\epsilon\,,\,e^+-\epsilon)$, is equal to $1$
in $(e^-+2\epsilon\,,\,e^+-2\epsilon)$, and  satisfies
$\sup_{\Gamma^1}|\partial_s \alpha|\le C$. Such a function obviously exists
because of the hypothesis $l_*>0$. Fix an origin vertex $v_0$
and, minimizing over all paths of the form 
$\gamma:[0\,,\,a] \rightarrow \Gamma^1$ from  $v_0$ to $s\in \Gamma^1$, 
parametrized by arclength, set
$$ 
\rho_*(s)=\min_{\gamma} \lambda(\gamma)\quad \text{where}\quad
\lambda(\gamma) = \int_0^a \alpha\bigl(\gamma (\tau)\bigr)\,d\tau.
$$
Observe that the function $\rho_*$ tends to infinity at infinity and that it 
is constant in a neighbourhood of any vertex $v$. If we had
$\rho_*\in \mathcal C^\infty(\Gamma^1)$, it would thus be a good candidate
for an edge-adapted exhaustion function. 

However, this function is not smooth
at  points $s$ in the interior of an edge $(e^-\,,\,e^+)$ 
with the property that
there are two minimizing paths $\gamma_1$ and $\gamma_2\,$, one passing 
through $e^-$, the other through $e^+$ and such that $\rho_*$ is not constant
in a neighbourhood of $s$. Observe that in this case, $s$ is a point of local 
maximum for $\rho_*\,$, and $\rho_*(s)\ge \max\{\rho_*(e^-),\rho_*(e^+)\}$. 
It follows that such an edge is never used by minimizing paths except 
those ending within the edge itself. Thus, changing $\alpha$ along such an 
edge has no effect on the values of $\rho_*$ elsewhere. 
Assume without loss of generality that $\rho_*(e^-)\le \rho_*(e^+)$
and replace $\alpha_e$ by a smaller smooth function $\widetilde{\alpha}_e$ 
satisfying $|\partial_s \widetilde{\alpha}_e|\le C$  and such that  
$\int_{e^-}^{e^+} \widetilde{\alpha}_e(s)\,ds=\rho_*(e^+)-\rho_*(e^-)$.
We can do this along any of those ``bad'' edges.
Globally, this defines a new function $\widetilde{\alpha}$, and using the same
notation as above, we set
$$
\rho_1(s) =\min_{\gamma} \widetilde\lambda(\gamma), \quad\text{where}
\quad  \widetilde\lambda (\gamma)=
\int_0^a \widetilde{\alpha}\bigl(\gamma (\tau)\bigr)\,d\tau\,.
$$
By construction, we have $\rho_1=\rho_*$ excepts on edges where $\alpha_e\neq
\widetilde{\alpha}_e$.  In particular, $\rho_1=\rho_*$ on vertices.
Moreover, $\partial_s\rho_1$ has compact support within every open edge.
Clearly, $\rho_1$ tends to infinity at infinity (along with $\rho_*$)
and satisfies $|\partial_s \rho_1|\le 1$ and $|\partial^2_s\rho_1|\le C$.
That is, $\rho_1$ is an edge-adapted exhaustion function.
As explained in Remark \ref{rmks-adapted}(a), this yields an 
edge-adapted approximation of $\uno$, 
say $\varrho_{1,n}\,$, on $\Gamma^1$. 
A strip-adapted approximation of $\uno$ 
on $\GM$ is obtained by setting
$$
\varrho_n(\xi)=\varrho_{1,n}(s)\,\varrho_{\Ms,n}(x),\;\;\xi=(s,x)\in \GM.
\eqno{\square}
$$
\renewcommand{\qedsymbol}{}\end{proof}
\vspace{-.7cm}

\begin{rmk} The conditions $\phi\equiv 1$, $\psi\equiv 1$ can be relaxed to
$$ 
\inf \phi > 0  \AND \sup |\partial_s \ln (\phi^{(n-1)/2}\psi)| < \infty\,.
$$
\end{rmk}

Our next result deals with the treebolic spaces $\HT(\pp,\qq)$. 

\begin{pro}  \label{pro-tbexhaust}
The treebolic space $\HT(\pp,\qq)$ equipped with
$\phi,\psi$ as in Example \ref{main-ex-HT} admits a 
strip-adapted exhaustion.
\end{pro}
\begin{proof} We will use freely  the notation introduced in Section 
\ref{geometry}. First we construct a smooth function 
$\eta: (0\,,\,\infty)\rightarrow (0\,,\,\infty)$ such that
$\eta\equiv 1$ on $\bigl(1-1/(8\qq)\,,\,1+\qq/8\bigr)$ and  
$\eta ( \qq^ky)=\qq^k\eta (y)$
for all $k\in \mathbb Z$. Obviously there is a $C>0$ such that 
this function satisfies 
$$   
C^{-1}\le y\,\eta(y)\le C\,,\quad \sup_{y>0}|\eta'(y)|\le C\,,
\quad\sup_{y>0} \{y\,|\eta''(y)|\}\le C.
$$ 

As a first step, consider the case $\pp=1$, $\qq>1$ where $\HT(1,\qq)$
is the upper half-space with the horizontal lines $\{z= x + \im y : y=\qq^k\}$
marked as bifurcation lines. Consider the function
$$
\delta(z)= \log \left(1 +\frac{1+x^2+ y^2}{y}\right).
$$
Away from the point $\im$, this is comparable with the hyperbolic distance 
between  $z$ and the point $\im$. Computing partial derivatives,
one easily checks that $y^2(|\partial_x\delta (z)|^2+\partial_y\delta(z)|^2)\le C_1$
and $y^2(|\partial_x^2\delta(z)|+|\partial_y^2\delta(z)|)\le C_1$ for some $C_1 > 0$. 
In particular, $\delta$ has bounded hyperbolic gradient and bounded hyperbolic 
Laplacian. Set
$$
\rho(z)= \delta(x + \im \eta(y)).
$$
Then it is not hard to check that $\rho$ is a 
strip-adapted exhaustion function on $\HT(1,\qq)$.
The role of $\eta$ is to make $\rho$ constant in $y$ along 
the lines $\{y=\qq^k\}$.

Let us now consider the general case $\HT(\pp,\qq)$, $\pp\ge 1,\qq>1$.
Recall that 
$\HT(\pp,\qq)=\{(z,w)\in \Hb \times \T_\pp:\hor(w)=\log_\qq y\}$. 
Hence, we can  consider $\rho$ as a function on $\HT(\pp,\qq)$
by setting $\rho(z,w)=\rho(z)$. This function satisfies all requirements for a
strip-adapted exhaustion function, except that it does not tend to $\infty$ 
along any fixed  horocyclic level $\hor(w)=\qq^k$.

To treat this difficulty, fix an end $\uf_0 \in \bd^*\T$. Let $V(\uf_0)$
be the set of all vertices $v\in \T^0_{\pp}$ such that $v\in \geo{\uf_0\,\om}$.
For any $v\in V(\uf_0)$, let $\mathbf T(v)$ be the set of those elements 
$w\in \T_{\pp}$ such that  $w\cf \uf_0=v$. This set 
$\mathbf T(v)$ is the maximal subtree of $\T_{\pp}$ containing $v$ and 
intersecting $\uf_0$ only at $v$. The tree $\T_\pp$ is the disjoint union
$$
\T_\pp=  \geo{\uf_0\,\om} \,\cup\left( \bigcup_{v\in  \T^0_{\pp}
\bigcap \geo{\uf_0\,\om}}\mathbf T(v)\setminus\{v\}\right),
$$
where (recall) $\geo{\uf_0\,\om}$ is the geodesic between $\uf_0$ and $\om$.
By construction, we have $\hor(w)\ge \hor(v)$ if $w\in \mathbf T(v)$.
Thus, for $(z,w)\in \HT(\pp,\qq)$ with $z=x + \im y$ and $w\in \mathbf T(v)$, 
we have $ y\ge \qq^{\hor(v)}$.

We define a function $\kappa$ on $\HT(\pp,\qq)$ by setting
$$
\kappa(z,w)=\begin{cases}0\,, &\text{if}\; w\in \geo{\uf_0\,\om},\\[3pt]
\log \bigl(\eta(\qq^{-\hor(v)}y)\bigr)\,,& \mbox{if}\; w\in \mathbf T(v).
\end{cases}
$$
This function $\kappa$ has the property that  it tends to infinity on 
$\HT(\pp,\qq)$ when its argument $(z,w)$ escapes to infinity along 
a fixed  horocycle $\{(z,w)\in \HT(\pp,\qq): \log_\qq(y)=\hor(w)=t\}$, 
$t\in \mathbb R$. This is because, as $(z,w)$ escapes to infinity
with $\log_\qq(y)=\hor(w)=t$, the vertex $v=v(w)\in \geo{\uf_0\,\omega}$ 
such that $w\in \mathbf T(v)$ must tend to $\om$ and thus $\hor(v)$ tends 
to $-\infty$.

Now, we set  
$$
\rho_1: \HT(\pp,\qq)\rightarrow (0\,,\,\infty)\,, \quad (z,w)\mapsto  
\rho_1(z,w)=\rho(z) +\kappa (z,w).
$$
From the construction, it is clear that  
$\rho_1$ is a strip-adapted exhaustion function.
\end{proof}

\section{Appendix: some results concerning
$\sqrt{-\Delta_{\Ms}}$}\label{appendix}

Let $(\Ms,g)$ be a Riemannian manifold (equipped with its Riemannian measure 
$dx$) and let $\Delta_{\Ms}$ be its Laplacian defined on 
$\mathcal C^\infty_c(\Ms)$. Abusing notation, we let $\Delta_{\Ms}$ denote also 
its Friedrichs extension. Let $h_{\Ms}(t,x,y)$ be the heat kernel (the smooth 
positive integral kernel of $e^{t\Delta_{\Ms}}$) and let 
$\sqrt{-\Delta_{\Ms}}$ be defined by spectral theory, that is,
$\sqrt{-\Delta_M}=\int_0^\infty \sqrt{\lambda} \,d E_\lambda\,$,
where $E_\lambda$ is a spectral resolution of $-\Delta_{\Ms}$.
The domain of $\sqrt{-\Delta_{\Ms}}$ is the Sobolev space 
$\mathcal W^1_0(\Ms)=\mathcal W^1_0\,$.

Let $\mathcal W^\alpha_0$ be the dual of $\mathcal W^{-\alpha}_0$ 
(under the identification of ${\mathcal L}^2(\Ms)$ with its own dual). 
Hence, for $\alpha>\beta>0$, we have
$$
\mathcal W^\alpha_0\subset \mathcal W^\beta_0\subset {\mathcal L}^2(\Ms)
\subset \mathcal W^{-\beta}_0\subset \mathcal W^{-\alpha}_0\,.
$$
The intersection $\mathcal W^\infty_0=\bigcap_{\alpha}\mathcal W^\alpha_0$ is 
dense in any $\mathcal W^\alpha_0\,$, and the operator  
$(\Id+\sqrt{-\Delta_{\Ms}}\,)^{\gamma}$, initially defined on 
$\mathcal W^\infty_0\,$, extends as a unitary operator from 
$\mathcal W^\alpha_0$ to $\mathcal W^{\alpha-\gamma}_0\,$. Moreover,
$$
(\Id+\sqrt{-\Delta_{\Ms}}\,)^{\alpha}
(\Id+\sqrt{-\Delta_{\Ms}}\,)^{\beta}= (\Id+\sqrt{-\Delta_{\Ms}}\,)^{\alpha+\beta},
\quad
\alpha,\beta\in \R,
$$
and $(\Id+\sqrt{-\Delta_{\Ms}}\,)^{0}=\Id$.
Because $\mathcal C^\infty_c(\Ms)\subset \mathcal W^\infty_0$ 
(with continuous embedding
when equipped with their natural families of seminorms), it is clear that
any $\mathcal W^\alpha$ can be understood as a space of distributions.

On ${\mathcal L}^2(\Ms)$, the operator $(\Id+\sqrt{-\Delta_{\Ms}}\,)^{-1}$ has an integral 
kernel given by
\begin{equation}\label{form-G}
G(x,y)=\frac{1}{\sqrt{\pi}}\int_0^\infty e^{-t} \int_0^\infty
\frac{e^{-u}}{\sqrt{u}}\, h_M(t^2/4u,x,y)\, du\,dt\,.
\end{equation}

It obviously satisfies
$$
\int_{\Ms}G(x,y)\,dy=\int_{\Ms}G(x,y)\,dx\le 1.
$$
It follows that, for any $f\in \mathcal C^\infty_c(\Ms)$, we have
\begin{eqnarray*}
\sqrt{-\Delta_{\Ms}}f &=&f +(-\Id+\sqrt{-\Delta_{\Ms}}\,)f\\
&=& f +(\Id+\sqrt{-\Delta_{\Ms}})^{-1}
(\Id+\sqrt{-\Delta_{\Ms}})(-\Id+\sqrt{-\Delta_{\Ms}}\,)f\\
&=& f+ (\Id+\sqrt{-\Delta_{\Ms}})^{-1}[-f-\Delta_{\Ms}f]\\
&\in& {\mathcal L}^1(\Ms)\cap {\mathcal L}^\infty(\Ms).
\end{eqnarray*}
Let now $f\in {\mathcal L}^1(\Ms)+{\mathcal L}^\infty(\Ms)$. 
The previous observation implies that we
can make sense of $\sqrt{-\Delta_{\Ms}} f$ explicitly
as a distribution on $\Ms$ by setting,
$$
\bigl[\sqrt{-\Delta_{\Ms} }f\bigr](h)
=\int_{\Ms} f \, \bigl[\sqrt{-\Delta_{\Ms}}\, h\bigr] \,dx
\quad \text{for}\;\;h\in \mathcal C^\infty_c(\Ms)\,.
$$

By (\ref{form-G}) and the local regularity of the heat kernel,
for any fixed precompact compact coordinate chart $(U;x_1,\dots,x_n)$ in $\Ms$
and any open set $\Omega\supset \overline{U}$, we have
\begin{equation} \label{Gsm1}
\forall y\in M\setminus \Omega,\;\;
\sup_{x\in U} \bigl|\partial_x^m G(x,y)\bigr| 
\le C_{U,\Omega,m} \,\inf_{x\in U}G(x,y)\quad\text{for all }\; 
y\in \Ms\setminus \Omega\,,
\end{equation}
where $m=(m_1,\dots,m_2)$ and 
$\partial_x^m f= \partial^{m_1}_{x_1}\dots\partial^{m_n}_{x_n}f$.
Furthermore, if $(U';y_1,\dots,y_n)$ is a relatively compact coordinate chart
with $U'\subset \Ms\setminus \overline{\Omega} $ then
\begin{equation} \label{Gsm2}
\sup_{x\in U}\sup_{y\in U'}
\bigl|\partial_x^m\partial_y^k G(x,y)\bigr| \le C_{U,U',m,k}\, .
\end{equation}

We need the following simple hypoellipticity type result. It is certainly
well-known but it does not seem very easy to find a precise reference.
(See e.g. {\sc Bogdan and Byczkowski}~\cite{BB}, where $(\Ms,g)$ is
Euclidean space.)
In particular, note that some care is needed because $\sqrt{-\Delta_{\Ms}}$
is not a local operator.

\begin{thm} \label{th-hypsqrt}
Let $f\in {\mathcal L}^2(\Ms)$ and let $F$ be the distribution
$F=(\Id+\sqrt{-\Delta_{\Ms}}\,)f$. Fix two open relatively compact sets
$\Omega \subset \Omega'\subset \Ms$ with $\overline{\Omega}\subset \Omega'$.
Assume that
\begin{itemize}
\item $F=0$ in $\Omega$, that is, $F(u)=0$ for all 
$u\in \mathcal C^\infty_c(\Omega)$, and
\item   $F|_{X\setminus \Omega'} \in {\mathcal L}^2(\Ms)$, that is,
there exists $h\in {\mathcal L}^2(\Ms)$ such that 
$$
F(u)=\int_{\Ms} h\,u\, dx \quad\text{for all }\;
u\in \mathcal C^\infty_c(\Ms\setminus \overline{\Omega'}).
$$
\end{itemize}
Then $f\in \mathcal C^\infty_{\mbox{\em \tiny loc}} (\Omega)$.
\end{thm}

\begin{proof} Without loss of generality, we can assume that $h=0$
in a neighbourhood of $\overline{\Omega}$. It then follows easily from
(\ref{Gsm1}) that
$$
(I+\sqrt{-\Delta_{\Ms}}\,)^{-1}h =Gh 
\in \mathcal C^\infty_{\mbox{\tiny loc}}(\Omega).
$$
Next, for any two open sets $\Omega_0,\Omega_1$ with
$\overline{\Omega_0}\subset \Omega_1$ and $\overline{\Omega_1}\subset \Omega$,
and any relatively compact neighbourhood $\Omega_2$ of $\overline{\Omega'}$, the
distribution $F-h$ is supported in $\Omega_2\setminus \overline{\Omega_1}\,$.
We can approximate this distribution
by functions in $F_j\in \mathcal C^\infty_c(\Ms)$ supported in
$\Omega_2\setminus \overline{\Omega_1}$ and such that there exist
a constant $C$, an integer $l$, and a finite covering of
$K=\overline{\Omega_2}\setminus \Omega_1$ by relatively compact charts
$(U^i,x^i_1,\dots,x^i_n)$, $i\in I$, such that for all $j$
$$ 
\int_{\Ms} F_j\, u \,dm\le
C \,\sup\left\{ \bigl|\partial_{x^i}^k  u(x)\bigr| : x\in U^i\,,\; i\in I\,,\;
k=(k_1,\dots,k_n) \;\text{with}\;  {\textstyle\sum} k_i\le l\right\}.
$$
It then  follows from (\ref{Gsm2}) that, given any local chart
$(U;x_1,\dots,x_n)$ contained in $\Omega_0$ and any integer $m$,
the functions $(\Id+\sqrt{-\Delta_M}\,)^{-1}F_j=G F_j$ satisfy
$$
\sup_j\,\sup\left\{
\bigl|\partial_x^mGF_j(x)\bigr| : x\in U\,,\; m=(m_1,\dots,m_n)\,,\; 
{\textstyle\sum} m_i\le m \right\} \le C.
$$
This implies that the limit distribution
$(\Id+\sqrt{-\Delta_{\Ms}}\,)^{-1}(F-h)=\lim_j GF_j$ can be represented by a 
smooth function in $\Omega_0\,$. Hence,
$$
f=(\Id+\sqrt{-\Delta_{\Ms}}\,)^{-1}F= (\Id+\sqrt{-\Delta_{\Ms}}\,)^{-1}h
+(\Id+\sqrt{-\Delta_{\Ms}}\,)^{-1}(F-h)
$$  
satisfies
$$
f|_\Omega= \bigl[(\Id+\sqrt{-\Delta_{\Ms}}\,)^{-1}F\bigr]\big|_{\Omega}
\in \mathcal C^\infty_{\mbox{\tiny loc }}(\Omega).
$$
This concludes the proof.\end{proof}

\end{document}